\documentclass[11pt, reqno]{amsart}
\usepackage{lmodern}
\usepackage{amsmath, amsthm, amssymb, amsfonts}
\usepackage[normalem]{ulem}
\usepackage{hyperref}
\usepackage[all,cmtip]{xy}
\usepackage{verbatim}
\usepackage{nccmath}
\usepackage{stmaryrd}
\usepackage{caption}
\usepackage{mathrsfs}
\usepackage{mathtools}
\usepackage{esvect}
\usepackage{cite}
\usepackage{bbm}
\usepackage{eucal}

\usepackage{mathbbol}
\usepackage{tabularx}
\usepackage[toc,page]{appendix}

\usepackage{tikz-cd}

\theoremstyle{plain}
\newtheorem{thm}{Theorem}[section]
\newtheorem{cor}[thm]{Corollary}

\newtheorem{lemma}[thm]{Lemma}
\newtheorem{prop}[thm]{Proposition}

\newtheorem*{ass}{Assumptions}

\newtheorem{thml}{Theorem}

\theoremstyle{definition}
\newtheorem{defn}[thm]{Definition}

\makeatletter
\newcommand\ackname{Acknowledgements}
\if@titlepage
\newenvironment{acknowledgements}{%
	\titlepage
	\null\vfil
	\@beginparpenalty\@lowpenalty
	\begin{center}%
		\bfseries \ackname
		\@endparpenalty\@M
\end{center}}%
{\par\vfil\null\endtitlepage}
\else

\fi
\makeatother

\theoremstyle{remark}
\newtheorem{rmk}[thm]{Remark}

\newcommand{\BC}{{\mathbb{C}}}

\newcommand{\BF}{{\mathbb{F}}}

\newcommand{\BL}{{\mathbb{L}}}

\newcommand{\BN}{{\mathbb{N}}}

\newcommand{\BP}{{\mathbb{P}}}
\newcommand{\BQ}{{\mathbb{Q}}}
\newcommand{\BR}{{\mathbb{R}}}

\newcommand{\BT}{{\mathbb{T}}}

\newcommand{\BZ}{{\mathbb{Z}}}

\newcommand{\CC}{{\mathcal C}}

\newcommand{\CF}{{\mathcal F}}
\newcommand{\CG}{{\mathcal G}}
\newcommand{\CH}{{\mathcal H}}

\newcommand{\CL}{{\mathcal L}}
\newcommand{\CM}{{\mathcal M}}
\newcommand{\CN}{{\mathcal N}}
\newcommand{\CO}{{\mathcal O}}

\newcommand{\CT}{{\mathcal T}}
\newcommand{\CU}{{\mathcal U}}

\newcommand{\FM}{{\mathfrak{M}}}

\newcommand{\ch}{{\mathrm{ch}}}
\DeclareMathOperator{\Hilb}{Hilb}
\newcommand{\rk}{{\mathrm{rk}}}

\newcommand{\Coh}{\mathfrak{Coh}}

\newcommand{\rCoh}{\mathfrak{Coh}_{r}}
\newcommand{\rpCoh}{\mathfrak{Coh}^{\sharp}_{r}}
\newcommand{\pCoh}{\mathfrak{Coh}^{\sharp}}
\newcommand{\Cohc}{\mathrm{Coh}}

\newcommand{\tr}{{\mathrm{tr}}}

\DeclareFontFamily{OT1}{rsfs}{}
\DeclareFontShape{OT1}{rsfs}{n}{it}{<-> rsfs10}{}
\DeclareMathAlphabet{\curly}{OT1}{rsfs}{n}{it}

\newcommand\Ext{\operatorname{Ext}}
\newcommand\Hom{\operatorname{Hom}}
\newcommand{\p}{\mathbb{P}}

\newcommand{\Mbar}{{\overline M}}

\newcommand{\td}{\mathrm{td}}

\newcommand{\Pic}{\mathop{\rm Pic}\nolimits}

\newcommand{\Sym}{{\mathrm{Sym}}}

\newcommand{\ev}{{\mathrm{ev}}}

\newcommand\Tor{\operatorname{Tor}}

\newcommand{\Eff}{\mathrm{Eff}}

\newcommand{\Aut}{\mathrm{Aut}}
\newcommand{\thickslash}{\mathbin{\!\!\pmb{\fatslash}}}

\newcommand{\FC}{\mathfrak{C}}

\begin{document}
	
\title[Quasimaps to moduli spaces of  sheaves]
{Quasimaps to moduli spaces of  sheaves}

	\author{Denis Nesterov}
\address{ETH Z\"urich, Departement Mathematik}
\email{denis.nesterov@math.ethz.ch}
\maketitle
	\begin{abstract}
		
		We develop a theory of quasimaps to a moduli space of sheaves $M$ on a surface $S$.  Under some  assumptions, we prove that moduli spaces of quasimaps are proper and carry a perfect obstruction theory. Moreover, they are naturally isomorphic to moduli spaces of sheaves on threefolds $S\times C$, where $C$ is a nodal curve. Using Zhou's  theory of entangled tails, we establish  a wall-crossing formula which therefore relates  the Gromov--Witten theory of  $M$ and the Donaldson--Thomas theory of $S\times C$ with relative insertions. We evaluate the wall-crossing formula for Hilbert schemes of points $S^{[n]}$, if $S$ is a del Pezzo surface. 
	
	\end{abstract}
	
\setcounter{tocdepth}{1}
\tableofcontents

	\section{Introduction} 
	
	\subsection {Overview}
	In \cite{CFKM}, Ciocan-Fontanine, Kim and Maulik  defined moduli spaces of stable quasimaps to GIT quotients, generalizing previously existing constructions \cite{CFKto, MDP, MM,TodaC}. Moduli spaces of stable quasimaps and stable maps are different compactifications of  moduli spaces of stable maps with smooth domains. There also exists a mixed theory of $\epsilon$-stable quasimaps that interpolates between the two.  In \cite{CFKmi, CFKse, CFK14}, wall-crossing formulas between quasimap invariants and Gromov--Witten invariants were conjectured. In \cite{YZ}, they were proved in full generality. In a nutshell, the difference between two theories is measured  by an $I$-function - a generating series of localised invariants associated to quasimaps from $\p^1$. 
	
	In this article, we develop a theory of quasimaps to moduli spaces of sheaves on surfaces. More specifically, let $M(\mathbf{v})$ be a moduli space of stable sheaves  in a class $\mathbf{v} \in K_{\mathrm{num}}(S)$ on a smooth projective complex surface $S$. The space $M(\mathbf{v})$ naturally embeds into the rigidified\footnote{We refer to Section \ref{rig} for the rigidification.} moduli space of all, not necessarily stable,  sheaves,  
	\[M(\mathbf{v}) \subset \rCoh(S, \mathbf{v}).\]
	A quasimap to $M(\mathbf{v})$  from a nodal curve $C$ is defined to be a map 
\[f\colon C \rightarrow	\rCoh(S, \mathbf{v}),\]
which maps to $M(\mathbf{v})$ at a general point of $C$. We construct  moduli spaces of such quasimaps, showing that they are proper and carry a perfect obstruction theory under some natural assumptions.

The most remarkable feature of quasimaps to moduli spaces of sheaves is that they can be given a completely sheaf-theoretic interpretation. A reason for that is simple: by construction, a map from a scheme $B$ to  $\rCoh(S, \mathbf{v})$ is given by a sheaf on $S\times B$ flat over $B$. In particular, for any quasimap from a curve $C$, we obtain a naturally associated sheaf on $S\times C$.
This correspondence is of fundamental importance. On one hand, we use it to prove all aforementioned properties of  moduli spaces of quasimaps. On the other hand, it connects quasimaps to Donaldson--Thomas theory. In fact, this holds on all levels: stability, obstruction theory, insertions, etc. In short,

\[\text{Quasimap theory of }M(\mathbf{v})=\text{Donaldson--Thomas theory of } S\times C,\] 
such that primary insertions on the left correspond to relative insertions on the right.

Following \cite{YZ}, we prove wall-crossing formulas in our set-up. They have exactly the same form as those of GIT quasimaps. With respect to the identification of theories above,  $I$-functions become 1-leg Vertex functions associated to the localised Donaldson--Thomas theory on $S\times \p^1$, introduced in \cite{MNOP1, MNOP}. Hence the wall-crossing formulas express the Donaldson--Thomas theory of $S\times C$ in terms of the Gromov--Witten theory of $M(\mathbf{v})$ and Vertex functions, and vice versa. 
\subsection{Results} \label{Results}
Throughout the article, $S$ is a smooth projective  surface over the field of complex numbers $\BC$. We fix a very ample line bundle $\CO_{S}(1)\in \Pic(S)$, a class $\mathbf{v}$ in the numerical $K$-group $K_{\mathrm{num}}(S)$, and another class $\mathbf{u}$ in the algebraic $K$-group $K_{0}(S)$.  We equip $K$-groups with the Euler norm $\chi$, which is defined as $\chi(F):=\sum_i(-1)^{i}h^i(S,F)$ for a coherent sheaf $F$ on $S$. 

\begin{ass} We make the following assumptions: 
	\
	\begin{itemize}
	\item  $h^{1}(S)=0,$ 
	\item $\rk(\mathbf{v})>0$,
	\item $\chi(\mathbf{v}\cdot \mathbf{u})=1$,
	\item  for $\mathbf{v}$ and $\CO_S(1)$,  all  semistable\footnote{Semistability of sheaves is defined with respect to reduced Hilbert polynomials.} sheaves are  stable.	
\end{itemize}
Section \ref{general} explains why these assumptions are made. 
\end{ass}


A quasimap  $f\colon (C,\mathbf{p}) \rightarrow	\rCoh(S, \mathbf{v})$ is said to be stable, if 
\begin{itemize}
	\item nodes and markings  $\mathbf{p}$ are mapped to $M(\mathbf{v})$,
\item $C$ does not have rational tails\footnote{Rational components with one special point, i.e.\ with one separating node or one marking.},
\item $|\Aut(f)|<\infty$. 
\end{itemize}
More generally, for a real number $\epsilon \in \BR_{>0}$, in Definition \ref{stabilityqm} we introduce $\epsilon$-stable quasimaps, which specialise to stable quasimaps, if $\epsilon \ll 1$, and to stable maps\footnote{If $g\neq 0$ or $n\neq 0$, then it is sufficient to take $\epsilon>1$.},  if $\epsilon>2$. The degree of quasimaps is defined by pulling back determinant line bundles on  $\rCoh(S, \mathbf{v})$.  
Let 
\[Q^\epsilon_{g,N}(M(\mathbf{v}), \beta)\]
be a moduli space of $\epsilon$-stable quasimaps of genus $g$ and degree $\beta$ with $N$ markings. The starting point of the theory of quasimaps to moduli spaces of sheaves is the following result. 

\begin{thml}[Theorem \ref{mapssheaves}] \label{thm1} Under our assumptions, there exists a moduli space $M_{\mathbf{v},\check{\beta}}^{\epsilon}(S\times C_{g,N})$ parametrizing sheaves on threefolds $S \times C$ for varying nodal curves $C$, and a natural isomorphism
\[ Q^\epsilon_{g,N}(M(\mathbf{v}), \beta) \cong  M^{\epsilon}_{\mathbf{v},\check{\beta}}(S\times C_{g,N}),\] 
such that the degree of quasimaps $\beta$ together with the class $\mathbf{v}$ determine the Chern character of sheaves on $S\times C$. 
\end{thml}
Remark \ref{modulisheaves} gives an explicit description of  sheaves $F$, such that  \[F \in M_{\mathbf{v},\check{\beta}}^{\epsilon}(S\times C_{g,N})(\BC), \quad \text{ for } \epsilon \ll 1.\]  
An identification similar to  Theorem \ref{thm1} also holds if one considers a moduli space of objects in a perverse heart of $\mathrm{D^b}(S)$ given by the t-structure associated to a torsion pair. A particularly useful torsion pair is provided by 0-dimensional objects and their right-orthogonal complement in $\Cohc(S)$. We denote the rigidified moduli stack of objects in the associated heart by $\rpCoh(S, \mathbf{v})$. In Section \ref{perverse}, this is discussed in the context of Hilbert schemes of points, which naturally embed into  $\rpCoh(S, \mathbf{v})$,
\[S^{[n]} \subset \rpCoh(S, \mathbf{v}),\]
it is proved that quasimaps to the pair above correspond to stable pairs on $S\times C$ in the sense of \cite{PT}, Theorem \ref{PT}. 

Theorem \ref{thm1} gives us access to tools from the theory of moduli spaces of sheaves, which are used in combination with tools from the theory of quasimaps to prove the following result. 
\begin{thml}[Theorem \ref{proper}, \ref{obsthe}, \ref{PT}] Under our assumptions, a moduli space $Q^\epsilon_{g,N}(M(\mathbf{v}), \beta)$ is a proper  Deligne--Mumford stack with an obstruction theory. If the virtual tangent complex\footnote{There is a naturally defined sheaf-theoretic obstruction theory on $M(\mathbf{v})$, we refer to Section \ref{obstruction} for more details.} of $M(\mathbf{v})$ is a locally free sheaf in degree 0, then the obstruction theory is perfect. Moreover, the identification of Theorem \ref{thm1} respects naturally defined obstruction theories. The same holds for Hilbert schemes of points inside the moduli stack of objects in the perverse heart. 
	\end{thml}


We define descedent quasimap invariants, using virtual fundamental classes,
	\[\langle \gamma_{1}\psi^{k_1}, \dots, \gamma_{N}\psi^{k_N} \rangle^{\epsilon}_{g,\beta}:= \int_{[Q^{\epsilon}_{g,N}(M(\mathbf{v}),\beta)]^{\mathrm{vir}}}\prod^{i=N}_{i=1}\ev^{*}_{i}(\gamma_{i})\psi_{i}^{k_{i}},\]
where $\gamma_{1}, \dots, \gamma_{N} \in H^{*}(M(\mathbf{v}))$ and $\psi_{1}, \dots, \psi_{N}$ are $\psi$-classes associated to markings of the curves. With respect to the identification of Theorem \ref{thm1}, primary quasimap insertions, i.e.\ $k_i=0$ for all $i$, correspond to relative Donaldson--Thomas insertions.  

Let us now discuss the wall-crossing formula, which relates descendent quasimap invariants for different values of $\epsilon \in \BR_{>0}$. Let
\[V(M(\mathbf{v}),\beta) \]
be the space of  quasimaps  $f \colon \p^1 \rightarrow \rCoh(S,\mathbf{v})$ subject to the following conditions:
\begin{itemize}
	\item $f$ is of degree $\beta$,
	\item  $f(\infty) \in M(\mathbf{v})\subset \rCoh(S,\mathbf{v})$. 
\end{itemize}
By evaluating a quasimap at $\infty \in \p^1$, we obtain an evaluation map 
\[\mathrm{ev} \colon V(M(\mathbf{v}),\beta) \rightarrow M(\mathbf{v}).\]
There is a $\BC^*$-action on $\p^1$, 
\[t[x:y]=[tx:y], \quad t\in \BC^{*},\]
which induces a $\BC^*$-action on $V(M(\mathbf{v}),\beta)$. The space $V(M(\mathbf{v}),\beta)$ is not proper, but its $\BC^*$-fixed locus is.  Using virtual localisation, we obtain a localised virtual fundamental cycle, 
\[[V(M(\mathbf{v}),\beta)]^{\mathrm{vir}}\in H_*(V(M(\mathbf{v}),\beta)^{\BC^*})[z^\pm],\]
where $z$ is the $\BC^*$-equivariant parameter. We then define  $I$-functions, which are also known as Vertex functions, as follows,
\begin{align*}
	I_\beta(z)&:= \ev_* [V(M(\mathbf{v}),\beta)]^{\mathrm{vir}} \in H^*(M(\mathbf{v}))[z^\pm], \\
	\mu_{\beta}(z)&:=[zI_\beta(z)]_{z^{\geq 0}}  \in H^*(M(\mathbf{v}))[z].
\end{align*}
We now state the wall-crossing formula for values of $\epsilon$  corresponding to stable quasimaps and stable maps, which we denote by $0^+$ and $\infty$, respectively. 

\begin{thml}[Corollary \ref{wallcrossingHilb2}] \label{thm3} Assuming $(g, N)\neq (0,1)$, we have
	\begin{multline*}
		\langle \lambda_1 \psi^{k_1}_1, \dots,\lambda_N \psi^{k_N}_N \rangle^{0^+}_{g,\beta} = \langle \lambda_1 \psi^{k_1}_1, \dots,\lambda_N \psi^{k_N}_N \rangle^{\infty}_{g,\beta} \\
		+\sum_{\underline{\beta}} \langle \lambda_1 \psi^{k_1}_1, \dots,\lambda_N \psi^{k_N}_N, \mu_{\beta_1}(-\psi_{N+1}), \dots, \mu_{\beta_k}(-\psi_{N+k}) \rangle^{\infty}_{g,\beta_0}/k!,
	\end{multline*}
	where $\underline{\beta}$ runs through all $(k+1)$-tuples of effective quasimap classes 
	\[\underline{\beta}=(\beta_0, \beta_{1}, \dots, \beta_{k}),\]
	such that $\beta=\sum^{i=k}_{i=0}\beta_i$, and $\beta_{i}\neq 0$ for all $i\in \{1, \dots, k\}$. 
\end{thml}

 The wall-crossing invariants of $S^{[n]}$ for a del Pezzo surface $S$ can be explicitly computed, allowing us to evaluate the wall-crossing formula.  More precisely, we can compute wall-crossing invariants for the pair
  \[ S^{[n]} \subset \rpCoh(S, \mathbf{v}).\]
To state the result, recall that  for $n>1$,
Nakajima--Grojnowski operators provide the following identification, 
\[H_2(S^{[n]},\BZ) \cong H_2(S,\BZ)\oplus \BZ.\]
Using this identification, we denote a curve class in $H_2(S^{[n]},\BZ)$ by  
\[(\gamma,m) \in H_2(S,\BZ)\oplus \BZ.\] We then define the following generating series, 
\[
\langle \gamma_1,\dots, \gamma_N \rangle_{g,\gamma}^{\sharp,\epsilon}:=\sum_{m\geq 0}  \langle \gamma_1,\dots,\gamma_N \rangle_{g,(\gamma,m)}^{\sharp,\epsilon}y^{m},
\]
where the superscript ``$\sharp$" indicates that we consider perverse quasimap invariants. 
\begin{thml}[Corollary \ref{Delpezzowall}]
 Assume $N>2$, then for Hilbert schemes of points $S^{[n]}$ on a del Pezzo surface $S$, we have
\[\langle \gamma_1,\dots, \gamma_N \rangle_{g,\gamma}^{\sharp,0^+}=(1+y)^{\mathrm{c}_1(S)\cdot \gamma}\cdot \langle \gamma_1,\dots,\gamma_N\rangle^{\sharp,\infty}_{g,\gamma}.\]
\end{thml}
Since perverse quasimaps  can be identified with stable pairs on $S\times \p^1$, the above result relates the quantum cohomology of $S^{[n]}$ to the ring whose structure constants are given by the Pandharipande--Thomas invariants of $S\times \p^1$. The change of variables as above was predicted\footnote{Communicated to the author by Georg Oberdieck.} by Davesh Maulik. 	

In fact,  the result holds true more generally for $(g,N)\neq (0,1)$. One just needs to insert $I$-functions from Proposition \ref{IfunctionFano} into Corollary \ref{wallcrossingHilb2}. However, the formula becomes more involved, because for some critical values of  $(g,N,\beta)$ the string and divisor equations fail. 
\subsection{Hilbert schemes of points} Assume  that $M(\mathbf{v})$ is a Hilbert scheme of points $S^{[n]}$ on $S$, i.e.\ 
\[ \ch(\mathbf{v})=(1,0,-n).\]
We can consider the following enumerative theories related to $S^{[n]}$:
\vspace{0.3cm}
\begin{itemize}
\item $\mathsf{GW}(S^{[n]})$ - Gromov--Witten theory of $S^{[n]}$,
\vspace{0.3cm}

\item $\mathsf{GW}_{\mathrm{orb}}(S^{(n)})$ - orbifold Gromov--Witten theory of an orbifold symmetric product $[S^{(n)}]$,
\vspace{0.3cm}

\item $\mathsf{GW}_{\mathrm{rel}}(S\times C_{g,N})$ - Gromov--Witten theory of $S\times C_{g,N}/\Mbar_{g,N}$ with relative insertions,\vspace{0.2cm} \vspace{0.001cm}

\item  $\mathsf{DT}_{\mathrm{rel}}(S\times C_{g,N})$ -  Donaldson--Thomas theory of $S\times C_{g,N}/\Mbar_{g,N}$ with relative insertions.
\end{itemize}
\vspace{0.3cm}
\noindent These theories are related by the following fundamental conjectures:
\vspace{0.3cm}

\begin{itemize}
\item \textsf{CRC} - analytic continuation and a change of variables relating theories  \textsf{GW} and $\mathsf{GW}_{\mathrm{orb}}$, provided by the Crepant resolution conjecture, proposed in \cite{YR}, refined in \cite{BG,CIT},
\vspace{0.3cm}

 \item \textsf{DT/GW} - analytic continuation and a change of variables relating theories $\mathsf{GW}_{\mathrm{rel}}$ and $\mathsf{DT}_{\mathrm{rel}}$, provided by the Donaldson--Thomas/Gromov--Witten correspondence, proposed in \cite{MNOP1,MNOP}.
 \vspace{0.3cm}
 
\end{itemize}

\noindent The present work and, in particular, Theorem \ref{thm3} provide a relation between theories \textsf{GW}  and $\mathsf{DT}_{\mathrm{rel}}$: 
\vspace{0.3cm}
\begin{itemize}
 \item \textsf{DT/Hilb} - wall-crossing formulas relating theories \textsf{GW}  and $\mathsf{DT}_{\mathrm{rel}}$. 
 \vspace{0.3cm}
\end{itemize}
In the follow-up work \cite{N22}, we derive a similar relation between theories $\mathsf{GW}_{\mathrm{orb}}$ and $\mathsf{GW}_{\mathrm{rel}}$, called Gromov--Witten/Hurwitz wall-crossing:
\vspace{0.3cm}
\begin{itemize}
\item  \textsf{GW/H}  - wall-crossing formulas relating theories $\mathsf{GW}_{\mathrm{orb}}$ and $\mathsf{GW}_{\mathrm{rel}}$. 
\vspace{0.3cm}
\end{itemize}

All four relations together can be represented by the Square in Figure \ref{square}.  For $S=\BC^2$, the Square was established in a series of articles: \cite{BP, BG, OP10,OP10c,OP10b}, the culmination of which were \cite{PT19} and \cite{PT19b}, where relations were shown to hold on the level of cohomological field theories.  Establishment of these relations for $\BC^2$ is fundamental for many developments in the field of modern enumerative geometry, like the proof of \textsf{DT/GW} for a quintic 3-fold in \cite{PP}. Our wall-crossings give a more geometric context to these results.

\begin{figure} [h!]
	\vspace{2.5cm}
	\scriptsize
	\[
	\begin{picture}(200,75)(-30,-50)
		\thicklines
		\put(25,-25){\line(1,0){30}}
		\put(95,-25){\line(1,0){30}}
		\put(25,-25){\line(0,1){40}}
		\put(25,30){\makebox(0,0){\textsf{DT/Hilb}}}
		\put(25,20){\makebox(0,0){\textsf{wall-crossing}}}
		\put(25,35){\line(0,1){40}}
		\put(125,-25){\line(0,1){40}}
		\put(125,30){\makebox(0,0){\textsf{GW/H}}}
		\put(125,20){\makebox(0,0){\textsf{wall-crossing}}}
		\put(25,75){\line(1,0){35}}
		\put(90,75){\line(1,0){35}}
		\put(125,35){\line(0,1){40}}
		\put(140,85){\makebox(0,0){$\mathsf{GW}_{\mathrm{orb}}(S^{(n)})$}}
		\put(15,85){\makebox(0,0){$\mathsf{GW}(S^{[n]})$}}
		\put(75,75){\makebox(0,0){\textsf{CRC}}}
		\put(5,-35){\makebox(0,0){$\mathsf{DT_{rel}}(S\times C_{g,N})$}}
		\put(150,-35){\makebox(0,0){$\mathsf{GW_{rel}}(S\times C_{g,N})$}}
		\put(75,-25){\makebox(0,0){\textsf{PT/GW}}}
	\end{picture}
	\]
	\caption{The Square}
	\label{square}
	\vspace{-0.2cm}
\end{figure}
\subsection{More applications}

This article is the first in a series of three articles. The two other are \cite{NK3,N22}, the latter was already mentioned in the previous section. In \cite{NK3}, we focus on quasimaps to moduli spaces of sheaves on a $K3$ surface.  Moduli spaces of sheaves on $K3$ surfaces require a special treatment due to existence of holomorphic symplectic forms and, consequently, vanishing of standard virtual fundamental classes.  

 Using these wall-crossings, we establish the following results. In \cite{NK3}, we prove:
\begin{itemize}
	\item  quantum cohomology of $S^{[n]}$ is determined by relative Pandharipande--Thomas theory of $S\times \p^1$, if $S$ is a $K3$ surface, conjectured in \cite{OP1},
	\item the Hilbert-schemes part of the Igusa cusp form conjecture, conjectured in \cite{OPa};
	\item  relative higher-rank/rank-one Donaldson--Thomas  correspondence for  $S \times \p^1$ and $S \times E$, if $S$ is a $K3$ surface and $E$ is an elliptic curve,
	\item  relative Donaldson--Thomas/Pandharipande--Thomas  correspondence for  $S \times \p^1$, if $S$ is a $K3$ surface.
\end{itemize}
In \cite{N22}:

\begin{itemize}
	
	\item 3-point genus-0 $\mathsf{CRC}$ in the sense of \cite{BG} for the pair $S^{[n]}$ and $[S^{(n)}]$, if $S$ is a toric del Pezzo surface,
	\item the geometric origin of $y=-e^{iu}$ in \textsf{PT/GW} through $\mathsf{CRC}$.
\end{itemize}

Moreover,  the quasimap wall-crossing played a crucial role in establishing a holomorphic anomaly equation for $K3^{[n]}$ in \cite{Ob22}.

From the perspective of mathematical physics, the quasimap wall-crossing is related to so-called \textit{dimensional reduction}. For example, it was used in \cite{KW}. In fact, our quasimap wall-crossing for moduli spaces of rank 0 sheaves is one of the algebro-geometric aspects of \cite{KW}. This is addressed in \cite{NHiggs}. For more on dimensional reduction in a mathematical context, we refer to \cite{GLSY}, and in a physical context - to \cite{BJSV}.

\subsection{Methods}
Despite the fact that the stack $\rCoh(S, \mathbf{v})$ can be covered by Zariski open quotient substacks, we cannot reduce our quasimap theory to the one of GIT quasimaps from \cite{CFKM}. There are several reasons for that. Firstly, the stack $\rCoh(S, \mathbf{v})$ is not a local complete intersection, which is one of the requirements of \cite{CFKM}. Secondly, the stack $\rCoh(S, \mathbf{v})$ is unbounded. As we increase the degree of quasimaps, they are allowed to map further away from the stable locus. Miraculously, stability of quasimaps is enough to guarantee boundedness of their moduli spaces.  However,  it is not enough to consider a fixed bounded quotient substack of $\rCoh(S, \mathbf{v})$, and it is unclear if these quotient substacks are quotients of affine schemes, which is another requirement of \cite{CFKM}. 

The general outline of our proofs and definitions is inspired by \cite{CFKM}. However, due to the aforementioned reasons, we cannot reduce our proofs to the GIT set-up.  Instead, we use  results from geometry of moduli spaces of sheaves in an essential way. For example, Langton's semistable reduction and the Mumford--Castelnuovo regularity  are indispensable tools in our endeavours. Overall, the theory requires juggling map-theoretic and sheaf-theoretic methods at the same time. 

Luckily, as soon as we establish foundational results, the constructions of \cite{YZ} immediately apply to our set-up. Hence the proof of  wall-crossing formulas in our set-up is exactly the same as in the GIT case.   

\subsection{Generalisations} \label{general}
Let us comment on the assumptions made in Section \ref{Results}. Firstly,  $h^{1}(S)=0$ is needed to ensure that cohomology of $S\times C$ does not contain algebraic classes of the form $\alpha \otimes \alpha' \in H^{\mathrm{odd}}(S)\otimes H^{\mathrm{odd}}(C)$. These classes are not seen by determinant line bundles on $\rCoh(S,\mathbf{v})$, which we use to define degrees of quasimaps. However, if we work with a fixed curve $C$, allowing it to sprout $\p^1$-bubbles, this should not be an issue, because in this case we can define the degree directly via the Chern character of sheaves on $S\times C$ (note that $h^1(\p^1)=0$). In general, one needs a more refined notion of degree of quasimaps to capture the geometry of such classes. On the other hand, all of our arguments should extend to the case $h^{1}(S)\neq0$, as they involve modifications of sheaves along fibers, which take place in $H^{\mathrm{even}}(S)\otimes H^{\mathrm{even}}(C)$.

The only place where we use the assumption $\rk(\mathbf{v})>0$ is  Theorem \ref{positive1}, which is needed to establish positivity of certain line bundles with respect to quasimaps in Proposition \ref{positive}. This is essential for proving that moduli spaces of quasimaps are quasi-compact in  Section \ref{Sectionalg}. The case of 1-dimensional sheaves on local curves is treated in \cite{NHiggs}. The author was not able to find an appropriate reference for results of \cite[Section 8]{HL}
in the context of torsion sheaves on an arbitrary surface.

The class $\mathbf{u} \in K_0(S)$, such that $\chi(\mathbf{v}\cdot \mathbf{u})=1$, is needed to construct a universal family on the rigidified stack $\rCoh(S, \mathbf{v})$. However, there always exists a finite gerbe over $\rCoh(S, \mathbf{v})$ with a universal family. In this case, we need to work with orbifold quasimaps.

Finally, we assume that all semistable sheaves are stable to guarantee that $M(\mathbf{v})$ is a projective scheme. The case of moduli spaces with strictly semistable sheaves is partly discussed in the context of Higgs bundles in \cite{NHiggs}. 

\subsection{Acknowledgments }First and foremost I would like to thank my advisor Georg Oberdieck. Not a single result of this work would be possible without his inexhaustible support and guidance. I am grateful to Daniel Huybrechts and Thorsten Beckmann  for multiple discussions on related topics and for reading preliminary versions of the present work. 
My understanding of quasimaps has improved greatly through discussions with Luca Battistella, and a great intellectual debt is owed to Yang Zhou for his theory of entangled tails, without which the wall-crossings would not be possible.  I also thank the Mathematical Institute of University of Bonn for stimulating academic environment and excellent working conditions. 

Finally, I want to express the words of gratitude to an anonymous referee for carefully reading the article and providing numerous comments which significantly improved the presentation and quality of proofs.

\subsection{Notation and conventions}
We work over the field of complex numbers $\BC$. We set 
$e_{\BC^*}(\BC_{\mathrm{std}})=z,$
where $\BC_{\mathrm{std}}$ is the weight 1 representation of $\BC^*$ on the vector space $\BC$. All functors are derived, unless stated otherwise. Cohomologies and homologies have rational coefficients, unless stated otherwise. All curves are assumed to be projective. 

Let $N$ be a semigroup and $\beta \in N$ be its generic element. By $\BQ[\![ q^\beta ]\!]$ we will denote the (completed) semigroup algebra 
$\BQ[\![ N]\!]$. In our case, $N$ will be various semigroups of effective curve classes.

After fixing an ample line bundle $\CO_S(1)$ on a surface $S$,  we define $\deg(F)$ to be the degree of a sheaf $F$ with respect to $\CO_S(1)$. By a  \textit{general fiber} of a sheaf $F$ on $S \times C$, we will mean a fiber of $F$ over some dense open subset of the curve $C$. 

\section{Stack of coherent sheaves}
\subsection{Preliminaries} \label{rig} Let $S$ be a smooth projective surface over the field of complex numbers $\BC$. We assume that the first Betti number of $S$ vanishes,  $h^{1}(S)=0.$  For an ample line bundle $\CO_S(1)$, semistability of sheaves is defined with respect to the reduced Hilbert polynomial, 
\[p(F,t)=\frac{\chi(S,F \otimes \CO_S(t))}{\rk(F)}.\]
The order of polynomials is given by the lexicographic order of their coefficients. A sheaf $F$ is semistable, if $p(F,t)\geq p(G,t)$ for all proper subsheaves $0\neq G \subset F$. It is stable, if the inequality is strict.

We fix  a very ample line bundle $\CO_{S}(1)\in \Pic(S)$, a class $\mathbf{v} \in K_{\mathrm{num}}(S)$, and another class\footnote{See Section \ref{sectiondet} for reasons why we need the class $\mathbf{u}$. } $\mathbf{u} \in K_{0}(S)$, such that:
\begin{itemize}
	\item $\rk(\mathbf{v})>0$,
	\item  $\chi(\mathbf{v}\cdot \mathbf{u})=1$,
	\item  for $\mathbf{v}$ and $\CO_S(1)$,  all semistable sheaves are  stable.
\end{itemize}

\subsection{Rigidification}
 Let 
\[\Coh(S,\mathbf{v}) \colon (Sch/ \BC)^{\circ} \rightarrow Grpd\]
be the stack of coherent sheaves on $S$ in the class $\mathbf{v}$, constructed, for example, in \cite[Section 08KA]{stacks-project}. There is a locus of sheaves stable with respect to $\CO_S(1)$, 
\[\CM(\mathbf{v}) \hookrightarrow \Coh(S,\mathbf{v}),\]
which is a $\BC^{*}$-gerbe over a projective scheme $M(\mathbf{v})$.  The $\BC^{*}$-automorphisms come from multiplication by scalars. In fact,  we can quotient out $\BC^{*}$-automorphisms of the entire stack $\Coh(S,\mathbf{v})$, as explained in \cite[Appendix C]{AGV}, thereby obtaining a \textit{rigidified} stack
\[\rCoh(S,\mathbf{v}):=\Coh(S,\mathbf{v}) \thickslash \BC^{*}.\]
A $B$-valued point of $\rCoh(S,\mathbf{v})$ can be represented by a pair $(\CG, \phi)$, where $\CG$ is a $\BC^*$-gerbe over $B$ and $\phi\colon \CG \rightarrow \Coh(S,\mathbf{v})$ is a $\BC^*$-equivariant map (here we will ignore 2-categorical technicalities, see \cite[Appendix C.2]{AGV} for more details). 
The moduli space $M(\mathbf{v})$ canonically embeds into the stack $\rCoh(S,\mathbf{v})$, giving rise to the following square
\[
\xymatrix{
	\CM(\mathbf{v}) \ar@{^{(}->}[r] \ar[d]^{\BC^{*}-\text{gerbe}} &
	\Coh(S,\mathbf{v}) \ar[d]^{\BC^{*}-\text{gerbe}}
	\\
	M(\mathbf{v}) \ar@{^{(}->}[r] & \rCoh(S,\mathbf{v})
}
\]


\subsection{Determinant line bundles} \label{sectiondet}
Let $\CF$ be the universal sheaf on $S\times \Coh(S,\mathbf{v})$, then we have a naturally defined determinant-line-bundle map,

\[\lambda\colon K_{0}(S) \rightarrow \Pic(\Coh(S,\mathbf{v})),\]
constructed as the following composition,
\begin{multline*}
	\lambda \colon K_{0}(S)\xrightarrow{p_{S}^{*}} K_{0}(S\times \Coh(S,\mathbf{v})) \xrightarrow{\cdot [\CF]}  K_{0}(S\times \Coh(S,\mathbf{v})) \\
	 \xrightarrow{p_{\Coh(S,\mathbf{v})*}} K_{0}(\Coh(S,\mathbf{v}))\xrightarrow{\det} \Pic(\Coh(S,\mathbf{v})).
	\end{multline*}
The construction of $\lambda$ requires some care, because  $ \Coh(S,\mathbf{v})$ is not of finite type. The correct approach is provided by Waldhausen's $K$-theory \cite{Wald}, which was used to construct a universal determinant map for perfect complexes in \cite[Section 3]{STV}. Perfectness of $\CF$ and smoothness of $p_{\Coh(S,\mathbf{v})}$ are crucial. 

In general, the weight of a line bundle $\lambda(u)$ with respect to the $\BC^*$-scaling is equal to the Euler characteristics $\chi(\mathbf{v}\cdot u)$,
\[w_{\BC^{*}}(\lambda(u))=\chi(\mathbf{v}\cdot u).\]
There are two types of classes that will be of interest to us. Firstly, a class $u \in K_{0}(S)$, such that $\chi(\mathbf{v} \cdot u)=1$, gives a trivilisation of the  $\BC^{*}$-gerbe 
\[\tau \colon \Coh(S,\mathbf{v}) \rightarrow \rCoh(S,\mathbf{v}).\]
More precisely, since $B\BC^*$ parametrises line bundles, a trivialisation is given by the map 
\begin{equation} \label{sec}
(\tau,\lambda(u))\colon \Coh(S,\mathbf{v})\xrightarrow{\sim} \rCoh(S,\mathbf{v})\times B\BC^*.
\end{equation}
Equivalently, by \cite[{Lemma 06QG}]{stacks-project}, the trivialisation is given by a section of $\tau$, 
\begin{equation*}
	s_{u}\colon \rCoh(S,\mathbf{v}) \rightarrow \Coh(S,\mathbf{v}),
\end{equation*}
which  is induced by the descend of the twisted universal family \[\CF\otimes p_{\Coh(S,\mathbf{v})}^*\lambda(u)^{-1}\] from $S\times \Coh(S,\mathbf{v})$ to $S\times \rCoh(S,\mathbf{v})$.  We will need this section to lift quasimaps from $\rCoh(S,\mathbf{v})$ to $\Coh(S,\mathbf{v})$, as the latter stack has a better modular interpretation. For this reason, we fixed the class $\mathbf{u} \in K_0(S)$. 
\\

On the other hand, for a class $u \in K_{0}(S)$, such that $\chi(\mathbf{v} \cdot u )=0$, the line bundle $\lambda(u)$ descends to $\rCoh(S,\mathbf{v})$. Consider
\[K_{\mathbf{v}}(S):=\mathbf{v}^{\perp}=\{u \in K_{0}(S) \mid  \chi(\mathbf{v} \cdot u)=0 \}  \subset K_{0}(S),\]
then  $\lambda$ restricted to $K_{\mathbf{v}}(S)$ descends to a map to $\Pic(\rCoh(S,\mathbf{v}))$, 
\[\lambda_{\mathbf{v}}\colon K_{\mathbf{v}}(S) \rightarrow \Pic(\rCoh(S,\mathbf{v})).\]
The class $\mathbf{v}$ will be frequently dropped from the notation in $\lambda_{\mathbf{v}}$, when it is clear what stack is considered.  We define 
\[\Pic_{\lambda}(\Coh(S,\mathbf{v})):=\mathrm{Im}(\lambda), \quad \Pic_{\lambda}(\rCoh(S,\mathbf{v})):=\mathrm{Im}(\lambda_{\mathbf{v}}).\]
There exists a particular class of elements in $K_{\mathbf{v}}(S)$, which deserves special mention and will be used extensively later,
\begin{align*}
	u_{i}&:=-\rk(\mathbf{v})\cdot h^{i}+\chi(\mathbf{v}  \cdot h^{i})\cdot [\CO_{\mathrm{pt}}], \\
\CL_{i}&:=\lambda(u_{i}),
\end{align*}
where $\CO_{\mathrm{pt}}$ is a structure sheaf of a point $\mathrm{pt} \in S$, and $h=[\CO_{H}]$ for a hyperplane $H \in |\CO_{S}(1)|$. 
The importance of these classes is due to the following theorem.

\begin{thm}\label{positive1}
	The line bundles $\CL_{1}$ and $\CL_{0}\otimes \CL_{1}^{m}$ are nef and ample, respectively, on $M(\mathbf{v})$ for all $m \gg 0$. 
\end{thm}
\textit{Proof.} See \cite[Chapter 8]{HL}, we use the assumption that $\rk(\mathbf{v})>0$. \qed

\section{Quasimaps} \label{quasim}
\subsection{Preliminaries}Let $(X, \mathfrak{X})$ be either a pair $(M(\mathbf{v}), \rCoh(S,\mathbf{v}))$ or a pair $(\CM(\mathbf{v}), \Coh(S,\mathbf{v}))$. 
\begin{defn} \label{quasimapsdef}
	A map $f\colon (C,\mathbf{p}) \rightarrow \mathfrak{X}$ is a \textit{quasimap} to $(X, \mathfrak{X})$ of genus $g$ and of degree $\beta \in \Hom(\Pic_{\lambda}(\mathfrak{X}), \mathbb{Z})$, if 
	\begin{itemize}
		\item $(C,\mathbf{p})$ is a connected  marked nodal curve of genus $g$,
		\item $\CL \cdot_{f} C:=\deg(f^*\CL)=\beta(\CL)$ for all $\CL \in \Pic_{\lambda}(\mathfrak{X})$,
		\item $| \{ p\in C\mid  f(p) \in \mathfrak{X}\setminus X \} | < \infty$.
	\end{itemize}
	We will refer to the set $\{ p\in C| \ f(p) \in \mathfrak{X}\setminus X \}$ as \textit{base} points. A quasimap $f$ is \textit{prestable}, if 
	\begin{itemize}
		\item $\{ p\in C| \ f(p) \in \mathfrak{X}\setminus X \} \cap \{\mathrm{nodes}, \mathbf{p} \}=\emptyset$. 
	\end{itemize} 
	We define 
	\[\Eff(X, \mathfrak{X}) \subset  \Hom(\Pic_{\lambda}(\mathfrak{X}), \mathbb{Z})\]
	to be the cone of classes of quasimaps. 
	
\end{defn}

	Let  
	\begin{equation*}
		\Lambda:=\bigoplus_{p\geq 0} H^{p,p}(S)
	\end{equation*}
	be the $(p,p)$-part of the cohomology of the surface $S$. 
	For a smooth connected curve $C$, we  have the K\"unneth  decomposition of the $(p,p)$-part of the cohomology of the threefold $S\times C$,
	\begin{equation}\label{decomposition}
		\bigoplus_{p\geq 0} H^{p,p}(S\times C)=\Lambda \otimes H^0(C,\BC) \oplus \Lambda \otimes H^2(C,\BC)=\Lambda \oplus \Lambda(-1),
	\end{equation}
where $(-1)$ denotes the Tate twist. 
\subsection{Sheaves associated to quasimaps} \label{sheavesass}
	Consider a quasimap
	\[f \colon C \rightarrow \Coh(S,\mathbf{v}).\]
	 By construction of $\Coh(S,\mathbf{v})$, the quasimap $f$ is given by a sheaf $F$ on $S \times C$ which is flat over $C$.  However, we are mainly interested in the rigidified stack  $\rCoh(S,\mathbf{v})$. To associate a sheaf on $S\times C$ to a quasimap
	 \[f \colon C \rightarrow \rCoh(S,\mathbf{v}),\]
	  we take a lift 
	  \[s_{\mathbf{u}} \circ f \colon C \rightarrow \Coh(S,\mathbf{v})\] by the section $s_{\mathbf{u}}$ constructed in Section \ref{sectiondet}.  We then take the sheaf $F$ associated to $s_{\mathbf{u}} \circ f$.  We fix a section $s_{\mathbf{u}}$  for the rest of the article by fixing the class $\mathbf{u} \in K_{0}(S)$, such that $\chi(\mathbf{v} \cdot \mathbf{u})=1$.   Sheaves associated to lifts $s_{\mathbf{u}} \circ f $ will be characterized in Section \ref{sectionrel}.
	  \begin{rmk} \label{liftsrem} In fact, a quasimap $f \colon C \rightarrow \rCoh(S,\mathbf{v})$ always admits a lift to some quasimap $f' \colon C \rightarrow \Coh(S,\mathbf{v})$. Indeed,  by \cite[Appendix C.2]{AGV} a map $C \rightarrow \rCoh(S,\mathbf{v})$ is given by a $\BC^*$-gerbe $\CG$ over $C$ with an $\BC^*$-equivariant map $\phi\colon \CG \rightarrow \Coh(S,\mathbf{v})$. It can be checked that 
	  	\[H_{\mathrm{fppf}}^2(C, \CO_C^*)=H^2_{\text{\'et}}(C, \CO_C^*)=0,\] 
	  	hence a $\BC^*$-gerbe is always trivial on $C$. This allows to lift the map. Different lifts of $f$ are related by tensoring an associated sheaf $F$ on $S\times C$ by  line bundles coming from $C$. However, the drawback of this construction is that it does not work in families, which we need to identify moduli spaces of quasimaps with moduli spaces of sheaves on $S\times C$. Hence we choose to work with the lift given by $s_{\mathbf{u}}$ from the beginning. 
	  	
	  	\end{rmk}
	
	Let us discuss how the degree of a quasimap is related to the Chern character of $F$. For a smooth curve $C$, the Chern character has two components with respect to the decomposition in (\ref{decomposition}),
	\[\ch(F)=(\ch(F)_{\mathrm{f}}, \ch(F)_{\mathrm{d}}) \in \Lambda \oplus \Lambda(-1),\]
	where the subscripts ``f" and ``d" stand for $\textit{fiber}$ and $\textit{degree}$, respectively.  As the notation suggests, \[\ch(F)_{\mathrm{f}}=\ch(\mathbf{v}),\]
	which can be seen by pulling back $\ch(F)$ to a fiber over $C$ and using the flatness of $F$.  On the other hand, the degree component of $F$ is related to the degree of a quasimap in the following way. 
	\begin{lemma} \label{degreechern}  Assume $C$ is smooth and let $f \colon C \rightarrow \Coh(S,\mathbf{v})$ be a quasimap of degree $\beta$, then
		
		\label{degree}
		\[\beta(\lambda(u))=\int_{S} \ch(u) \cdot \ch(F)_{\mathrm{d}} \cdot \td_S.\]
		\end{lemma}
		\textit{Proof.} By the functoriality of the determinant line bundle construction, we have 
	\[\beta(\lambda(u))= \deg(\lambda_{F}(u)),\]
	where $\lambda_{F}(u)$ is the determinant line bundle associated to the family $F$ on $S\times C$ and a class $u \in K_{0}(S)$. Using the Grothendieck--Riemann--Roch theorem and the projection formula, we obtain
	\begin{align*}
		\deg(\lambda_F(u))=\deg(p_{C*}(p^{*}_{S}u\cdot [F]))&=\int_{C} \ch(p_{C*}(p^{*}_{S}u\cdot [F]))\\
		&=\int_{S\times C}\ch(p^*_{S}u\cdot [F]) \cdot p^*_{S}\td_S\\
		&=\int_S \ch(u) \cdot p_{S*}\ch(F) \cdot \td_S\\
		&=\int_{S} \ch(u) \cdot \ch(F)_{\mathrm{d}} \cdot \td_S.
	\end{align*}
\qed 
\\

Now let $\beta_{\Lambda}\colon \Lambda \rightarrow \BQ$ be the descend of $(\beta \circ \lambda)_{\BQ}\colon K_{0}(S)_{\BQ} \rightarrow \BQ$ to $\Lambda$ via the Chern character,
\[
\begin{tikzcd}[row sep=small, column sep = small]
	& \Lambda \arrow[r, dashed,"\beta_{\Lambda}"] & \BQ \\
	& K_{0}(S)_{\BQ} \arrow{u}{\ch} \arrow{ur}[swap] {(\beta \circ \lambda)_{\BQ}}&
\end{tikzcd}
\]
which exists by the formula from Lemma \ref{degree}. The formula also shows that the descend $\beta_{\Lambda}$ and $\beta$ determine each other. We thereby obtain an expression of $\ch(F)_{\mathrm{d}}$ in terms of $\beta_{\Lambda}$,
\begin{equation*}\ch(F)_{\mathrm{d}}=\beta_{\Lambda}^{\vee}\cdot \td_{S}^{-1},
\end{equation*}
where $\beta_{\Lambda}^{\vee}$ is the dual of $\beta_{\Lambda}$ with respect to the cohomological intersection pairing on $\Lambda$. This motivates the following definition which allows us to pass between degrees of quasimaps and Chern characters of sheaves.

\begin{defn} \label{identify} We define a linear map 
	\[ \check{(...)}\colon \Eff(\CM(\mathbf{v}), \Coh(S,\mathbf{v})) \rightarrow \Lambda, \quad \beta \mapsto \check{\beta},  \] 
	where  $\check{\beta}:=\beta_{\Lambda}^{\vee}\cdot \td_{S}^{-1}$. 
\end{defn}

Assume now $C$ is a nodal curve and consider a sheaf $F$ on $S\times C$ flat over $C$. Let
\[\pi \colon \bigcup_i  S\times C_{i} \rightarrow S\times C\]
be the normalisation of $S\times C$, such that $S\times C_{i}$ are its irreducible components. We define 
\[F_i:= \pi^*F_{|S\times C_i}.\]
\begin{defn}
By using the natural identification (\ref{decomposition}), we define the Chern character of $F$ on a nodal curve $C$  as follows, 
\[\ch(F):=(\ch(\mathbf{v}), \sum_i \ch(F_i)_\mathrm{d})\in \Lambda \oplus \Lambda(-1).\]
\end{defn}
As a consequence of these definitions, we obtain a natural extension of Lemma \ref{degree} to the case of a singular curve $C$, 
\[\ch(F)_\mathrm{d}=\check{\beta}, \]
where $F$ is the sheaf associated to a quasimap of degree $\beta$.


\subsection{Positivity}\label{chpositive}
The aim of this section is to establish positivity of certain line bundles - Proposition \ref{positive}. We start with the following result, which is inspired by  \cite[Proposition 4.4]{BM}.
\begin{lemma}\label{key} Let $F$ be the sheaf on $S\times C$ associated to a map $f\colon C \rightarrow \Coh(S,\mathbf{v})$, then 
	\begin{align*}\CL_{1}\cdot_{f} C&=\mathrm{deg}(\mathbf{v})\mathrm{rk}(p_{S*}F)-\mathrm{rk}(\mathbf{v})\mathrm{deg}(p_{S*}F),\\
		\CL_{0}\cdot_{f} C&=\chi(\mathbf{v})\mathrm{rk}(p_{S*}F)- \mathrm{rk}(\mathbf{v})\chi(p_{S*}F),
	\end{align*}
	where $\deg(...) $ is the degree of a class with respect to $\CO_S(1)$.  
\end{lemma}
\textit{Proof.} 
Since $\chi(\mathbf{v} \cdot u_i )=0$, we get that 
\[\rk(p_{C*}(p^{*}_{S}u_i\cdot [F]))=0,\] 
hence 
\begin{align*}\CL_{i}\cdot_{f} C =\deg(p_{C*}(p^{*}_{S}u_i\cdot [F]))=\chi(p_{C*}(p^{*}_{S}u_i\cdot [F])) &=\chi(p^{*}_{S}u_i\cdot [F]) \\
	&=\chi(u_i\cdot p_{S*}[F]).
	\end{align*}
The claim then follows from the following computation for a class $u\in K_0(S)$,
\begin{align*}
	\chi( u_{1}\cdot u) &=-\rk(\mathbf{v})\chi(u \cdot h)+\chi(\mathbf{v}\cdot h)\chi([\CO_{\mathrm{pt}}]\cdot u) \nonumber \\
	&=-\rk(\mathbf{v})\left(\deg(u)-\frac{\rk(u)}{2}H^{2}-\frac{\rk(u)}{2}H\cdot \mathrm{c}_{1}(S) \right) \nonumber\\
	&+\left(\deg(\mathbf{v})-\frac{\rk(\mathbf{v})}{2}H^{2}-\frac{\rk(\mathbf{v})}{2}H\cdot \mathrm{c}_{1}(S)\right)\rk(u)\nonumber \\
	& =\rk(u)\deg(\mathbf{v})-\rk(\mathbf{v})\deg(u),  \nonumber \\
	\chi( u_{0}\cdot u) &= -\rk(\mathbf{v})\chi(u)+\chi(\mathbf{v})\chi([\CO_{\mathrm{pt}}]\cdot u) \nonumber \\
	&= \chi(\mathbf{v})\rk(u)-\rk(\mathbf{v})\chi(u),
\end{align*}
this finishes the proof. 
\qed
\\

In the following proposition, a quasimap is \textit{constant}, if it does not have base points and the induced map to $M(\mathbf{v})$ is constant. 

\begin{prop} \label{positive}
	Let $f\colon C \rightarrow \Coh(S,\mathbf{v})$ be a prestable quasimap of degree $\beta$.  Then there exists $m_{0} \in \BN$ which depends only on $\beta$, $\mathbf{v}$ and  $\CO_S(1)$, such that for all $m \geq m_{0}$, the quasimap is non-constant, if and only if
	\[ \CL_{0}\otimes \CL_{1}^{m} \cdot_{f} C>0.\]
	This also holds true for all subcurves $C'\subseteq C$ and the induced maps for the same choice of $m$.
\end{prop}  

We firstly establish a simpler version of the proposition, which will be necessary for its proof. 

\begin{lemma} 	\label{nonnegative}Let $f\colon C \rightarrow \Coh(S,\mathbf{v})$ be a prestable quasimap, then

\[\CL_1\cdot_f C\geq 0. \]

\end{lemma} 

\textit{Proof.}
 Let $F$ be the family of sheaves on $S\times C$ associated to $f$. Assume for simplicity that $f$ has one base point $b \in C$.  By Langton's semistable reduction \cite{La}, the sheaf $F$ can be modified at a point $b$ to a sheaf which is stable over $b$ and is isomorphic to $F$ away from $S\times \{b\} \subset S\times C$.  The modification is given by a finite sequence of short exact sequences,
\begin{align*}
	0 \rightarrow F^{1}\rightarrow &F^0 \rightarrow Q^{1} \rightarrow  0, \\
	&\vdots \\
	0\rightarrow F^{k} \rightarrow &F^{k-1} \rightarrow Q^{k} \rightarrow 0, 
\end{align*}
where $F^0=F$, the sheaf $F^{k}$ is stable over $b \in C$, and $Q^{i}$ is the maximal destabilising quotient sheaf of $F^{i-1}_{b}$, i.e.\ the quotient of $F^{i-1}_{b}$ by the maximal destabilising subsheaf $E^i$, 
\[0\rightarrow E^i \rightarrow F^{i-1}_{b} \rightarrow Q^i \rightarrow 0.\]
More precisely, if $F^{i-1}_b$ does not have torsion, we refer to \cite[Section 1.3]{HL} for the construction of the maximal destabilising subsheaf $E^i$.  If $F^{i-1}_{b}$ has torsion $T^i$, then $E^i$ is defined to be  the preimage of the maximal destabilizing subsheaf of the torsion free sheaf  $F^{i-1}_{b}/T^i$. If $ F^{i-1}_{b}/T^i$ is stable, then $Q^i=F^{i-1}_{b}/T^i$.  Note that by the maximality, $Q^i$ is  torsion free.  

   Since $Q^{i}$ is a destabilizing quotient of a sheaf in the class $\mathbf{v}$, we have
\begin{equation}\label{dest} \mathrm{deg}(\mathbf{v})\mathrm{rk}(Q^{i})-\mathrm{rk}(\mathbf{v})\mathrm{deg}(Q^{i})\geq 0.
\end{equation}
Applying the derived pushforward $p_{S*}$ to each sequence, we get distinguished triangles
\[p_{S*}(F^{i}) \rightarrow p_{S*}(F^{i-1}) \rightarrow Q^{i} \xrightarrow{}.\]
By Lemma \ref{key} we obtain that 
\begin{equation}\label{red}
	\CL_{1} \cdot_{f^{i-1}}C=\CL_{1} \cdot_{f^i} C+\mathrm{deg}(\mathbf{v})\mathrm{rk}(Q^{i})-\mathrm{rk}(\mathbf{v})\mathrm{deg}(Q^{i}),
\end{equation}
where $f^i$ is the quasimap associated to $F^i$. The line bundle $\CL_{1}$ is nef on $M(\mathbf{v})$ by Theorem \ref{positive1}, therefore
\begin{equation} \label{nef}
\CL_{1} \cdot_{f^k}C\geq 0,
\end{equation}
because $f^k$ does not have base points. The claim for a quasimap with one base point now follows from (\ref{dest}), (\ref{red}) and  (\ref{nef}).  By applying the semistable reduction to all base points at the same time, we extend the argument  to a quasimap with arbitrarily many base points. 
\qed
\\

Before proceeding to the proof of Proposition \ref{positive}, we derive an important consequence of Lemma \ref{nonnegative}.  Recall that the slope $\mu(F)$ of a sheaf $F$ is defined to be the ratio $\deg(F)/\rk(F)$. The maximal slope  $\mu_{\mathrm{max}}(F)$ is defined to be the slope $\mu(G)$ of the maximal destabilizing subsheaf $G$ of a sheaf $F$. 

\begin{cor} \label{bounds} Let $f\colon C \rightarrow \Coh(S,\mathbf{v})$ be a prestable quasimap of degree $\beta$, and let $F$ be the associated sheaf. Let $Q^i$ be a maximal destabilizing quotient sheaf that appears in the semistable reduction of $F$ in the proof of Lemma \ref{nonnegative}, then 
	\begin{align*}
	 -|\deg(\mathbf{v})| -\frac{\beta(\CL_{1})}{\rk(\mathbf{v})}\leq &\deg(Q^i)\leq |\deg(\mathbf{v})|, \\ 
	&\mu_{\mathrm{max}}(Q^i)< |\deg(\mathbf{v})|+\frac{\beta(\CL_{1})}{\rk(\mathbf{v})}.
	\end{align*}

	\end{cor}

 \textit{Proof.} Let us denote $\CL_1 \cdot_f C$ by $\beta(\CL_1)$. By Lemma \ref{nonnegative}, (\ref{dest}) and (\ref{red}), we have
  \[\beta(\CL_{1}) \geq \mathrm{deg}(\mathbf{v})\mathrm{rk}(Q^{i})-\mathrm{rk}(\mathbf{v})\mathrm{deg}(Q^{i})\geq0, \]
this bounds the degrees of $Q^i$  as follows,
\begin{equation}\label{chern}
	\frac{\deg(\mathbf{v})\rk(Q^i)-\beta(\CL_{1}) }{\rk(\mathbf{v})}\leq\deg(Q^i)\leq  \frac{\deg(\mathbf{v})\rk(Q^i)}{\rk(\mathbf{v})}. 
\end{equation}
We therefore get a uniform bound on $\deg(Q^i)$ depending on the sign of $\deg(\mathbf{v})$,  

\begin{equation} \label{boundschern}
\begin{aligned}
	-\frac{\beta(\CL_{1})}{\rk(\mathbf{v})}<\deg(Q^i)\leq \deg(\mathbf{v}),\quad &\text{if $\deg(\mathbf{v})\geq 0$, } \\
	\deg(\mathbf{v})-\frac{\beta(\CL_{1})}{\rk(\mathbf{v})}\leq \deg(Q^i) < 0, \quad& \text{if $\deg(\mathbf{v})< 0$},
\end{aligned}
\end{equation}
combing these two bounds together we get the first inequality in the statement of the corollary.

We now deal with $\mu_{\mathrm{max}}(Q^i)$. Let  $E^i$ be the maximal destabilizing subsheaf  of $F^{i-1}_{b}$ from the proof of Lemma \ref{nonnegative}. If $E^i$ is torsion, then 
\[\mu_{\mathrm{max}}(Q^i)=\mu(Q^i)\leq \mu(\mathbf{v}),\]
since $F^{i-1}_{b}/ E^i=Q^i$ and $Q^i$ is stable. 
 If  $E^i$ is non-torsion,  then  by the maximality we must have
\[\mu_{\mathrm{max}}(Q^i) \leq \mu(E^i).\]
 Indeed, if it was not the case, then by taking the preimage of the maximal destabilizing  subsheaf of $Q^i$, we would  obtain a subsheaf of $F^{i-1}_{b}$, whose slope is greater than $\mu(E^i)$, contradicting the maximality of $E^i$. Using the fact that $\deg(\mathbf{v})=\deg(Q^i)+\deg(E^i)$, and (\ref{chern}), we obtain a uniform upper bound on $\mu(E^i)$, which therefore bounds $\mu_{\mathrm{max}}(Q^i)$,
 
 \begin{equation}\label{boundsmax}
\begin{aligned}
	\mu_{\mathrm{max}}(Q^i) \leq \mu(E^i)< \deg(\mathbf{v})+\frac{\beta(\CL_{1})}{\rk(\mathbf{v})},\quad &\text{if $\deg(\mathbf{v})\geq 0$, } \\
	\mu_{\mathrm{max}}(Q^i) \leq \mu(E^i)\leq \frac{\beta(\CL_1)}{\rk(\mathbf v)}, \quad& \text{if $\deg(\mathbf{v})< 0$},
\end{aligned}
\end{equation}
combining these two bounds together with the bound for torsion $E^i$, we get the second inequality in the statement of the corollary. \qed 
 \\

\textit{Proof of Proposition \ref{positive}.} We now deal with the claim in the proposition. Assume $f$ is constant, then  $\CL_{0}\otimes \CL_{1}^{m} \cdot_{f} C=0$ for all $m$. 

 Conversely, assume $f$ is non-constant, and let $F$ be the associated sheaf. We apply the semistable  reduction to $F$ at all points at once, then by  Lemma \ref{key}, we obtain 
\begin{multline}\label{eq11}
	\CL_{0}\otimes\CL_{1}^{m} \cdot_{f} C=\CL_{0}\otimes\CL_{1}^{m} \cdot_{f^k} C \\
	+m\sum_{i} \mathrm{deg}(\mathbf{v})\mathrm{rk}(Q^i)-\mathrm{rk}(\mathbf{v})\mathrm{deg}(Q^i)+\sum_{i}\chi(\mathbf{v})\mathrm{rk}(Q^i)-\mathrm{rk}(\mathbf{v})\chi(Q^i),
\end{multline}
where $f^k$ is a stable map.  As in Lemma \ref{nonnegative}, $Q^i$ is destabilizing,  hence 
\begin{equation} \label{po}
	\mathrm{deg}(\mathbf{v})\mathrm{rk}(Q^i)-\mathrm{rk}(\mathbf{v})\mathrm{deg}(Q^i)\geq 0.
	\end{equation}
  We therefore have to analyse the terms
  \begin{equation}\label{terms}
  \chi(\mathbf{v})\mathrm{rk}(Q^i)-\mathrm{rk}(\mathbf{v})\chi(Q^i).
  \end{equation}
We will split our analysis, depending on whether (\ref{po}) is positive or zero.   If 
\[\mathrm{deg}(\mathbf{v})\mathrm{rk}(Q^i)-\mathrm{rk}(\mathbf{v})\mathrm{deg}(Q^i)=0,\]
then since $Q^i$ is a destabilizing quotient, it must be destabilizing in the second coefficient of the Hilbert polynomial (see Section \ref{Sectionslopes} for the definition of Hilbert polynomials), i.e.\  
\begin{equation} \label{pos1}
	\chi(\mathbf{v})\mathrm{rk}(Q^i)-\mathrm{rk}(\mathbf{v})\chi(Q^i)>0.
\end{equation}

Consider  now the case of $Q^i$, such that \[\mathrm{deg}(\mathbf{v})\mathrm{rk}(Q^i)-\mathrm{rk}(\mathbf{v})\mathrm{deg}(Q^i)>0.\]
By \cite[Corollary 3.3.3]{HL} and Serre's duality, dimensions of $H^0( Q^i)$ and $H^2(Q^i)$ can be bounded as follows,
\begin{equation} \label{b1}
	\begin{aligned}
h^0(Q^i)&\leq \frac{\deg(S)\rk(\mathbf{v})}{2}\left[\frac{\mu_{\mathrm{max}}(Q^i)}{\deg(S)}+\frac{\rk(\mathbf{v})-1}{2}+2\right]^2_+ ,\\
h^2(Q^i)&\leq \frac{\deg(S)\rk(\mathbf{v})}{2} \left[\frac{\mu_{\mathrm{max}}((Q^i)^\vee\otimes\omega_S)}{\deg(S)}+\frac{\rk(\mathbf{v})-1}{2}+2\right]^2_+,
\end{aligned}
\end{equation}
where $\deg(S)=\mathrm{c}_1(\CO_S(1))^2$,  and $[x]_+=\mathrm{max}\{0,x\}$.

By  Corollary \ref{bounds}, $\mu_{\mathrm{max}}(Q^i)$ is bounded from above. On the other hand, since $Q^i$ is torsion free, we can also bound $\mu_\mathrm{max}((Q^i)^\vee\otimes\omega_S)$ from above in terms of $\deg(Q^i)$ and $\mu_\mathrm{max}(Q^i)$. Firstly, 
\[ \mu_\mathrm{max}((Q^i)^\vee\otimes\omega_S)= \mu_\mathrm{max}((Q^i)^\vee)+\deg(\omega_S).\]
Now let $E \hookrightarrow (Q^i)^\vee$ be a saturated subsheaf\footnote{Saturated subsheaf is a subsheaf $G\subset F$, such that $F/G$ is torsion free; by the proof of \cite[Proposition 1.2.6]{HL} it is enough to check the slope on saturated subsheaves.}, then by dualizing we get a map $(Q^i)^{\vee\vee} \rightarrow E^\vee$, which is surjective away from a codimension 2 locus. We also have an embedding  $Q^i \hookrightarrow (Q^i)^{\vee\vee}$  which is an isomorphism away from a codimension 2 locus. For the purpose of calculating the slope, we therefore may assume that we have a surjection $Q^i \twoheadrightarrow E^\vee$. Let $K$ be its kernel, the slope of  $E^\vee$ can  be bounded as follows, 
\begin{align*}
\mu(E^\vee)=\frac{\deg(Q^i)}{\rk(E^\vee)}-\frac{\rk (K)}{\rk(E^\vee)}\mu(K) &\geq \frac{\deg(Q^i)}{\rk(E^\vee)}-\frac{\rk (K)}{\rk(E^\vee)}\mu_{\mathrm{max}}(Q^i).
\end{align*}
 This bounds the slope of $E$ from above,
\[ \mu(E)=-\mu(E^\vee)\leq \frac{\rk (K)}{\rk(E^\vee)}\mu_{\mathrm{max}}(Q^i)-\frac{\deg(Q^i)}{\rk(E^\vee)},\]
using Corollary \ref{bounds}, we obtain a uniform bound, 
\[ \mu(E) <\left(\rk(\mathbf{v})+1\right)\cdot \left(|\deg(\mathbf{v})|+\frac{\beta(\CL_{1})}{\rk(\mathbf{v})}\right).\]
We conclude that 
\begin{equation}\label{b2}
\mu_\mathrm{max}((Q^i)^\vee\otimes\omega_S)<\left(\rk(\mathbf{v})+1\right)\cdot \left(|\deg(\mathbf{v})|+\frac{\beta(\CL_{1})}{\rk(\mathbf{v})}\right)+\deg(\omega_S).
\end{equation}
Overall, using Corollary \ref{bounds}, (\ref{b1}) and (\ref{b2}), we obtain that $\chi(Q^i)$ can be bounded by an explicit constant $\chi_0$ that depends only on $\beta$, $\mathbf{v}$ and $\CO_S(1)$, assuming $S$ is fixed, 
\[\chi(Q^i) \leq h^0(Q^i)+h^2(Q^i)<\chi_0.\] 
 This allows us to uniformly bound the terms (\ref{terms}), 
\begin{align}\label{boundq}
	\chi(\mathbf{v})\mathrm{rk}(Q^i)-\mathrm{rk}(\mathbf{v})\chi(Q^i)>-\mathrm{rk}(\mathbf{v}) \cdot (|\chi(\mathbf{v})|+\chi_0 ).
\end{align}

Now let $m_{0} \in \BN$ be such that  $\CL_{0}\otimes \CL_{1}^{m_{0}}$ is ample on $M(\mathbf{v})$, which is possible by Theorem \ref{positive1}, and 
\begin{equation*}
\begin{aligned}
m_{0}\cdot(\mathrm{deg}(\mathbf{v})\mathrm{rk}(Q^i)-\mathrm{rk}(\mathbf{v})\mathrm{deg}(Q^i))-\mathrm{rk}(\mathbf{v}) \cdot (|\chi(\mathbf{v})|+\chi_0 )>0,\quad 
\end{aligned}
\end{equation*}
for all $Q^{i}$ satisfying \[\mathrm{deg}(\mathbf{v})\mathrm{rk}(Q^i)-\mathrm{rk}(\mathbf{v})\mathrm{deg}(Q^i)>0.\]   Such $m_0$ exists, because the quantity  \[\mathrm{deg}(\mathbf{v})\mathrm{rk}(Q^i)-\mathrm{rk}(\mathbf{v})\mathrm{deg}(Q^i)\] is an integer, and, in particular, it is at least 1, if it is positive.  More explicitly, we can choose  $m_0$ such that  $\CL_{0}\otimes \CL_{1}^{m_{0}}$ is ample on $M(\mathbf{v})$, and 
\[m_0 \geq \mathrm{rk}(\mathbf{v}) \cdot (|\chi(\mathbf{v})|+\chi_0 )+1.\]

By the definition of $m_0$ and the fact that $\mathrm{deg}(\mathbf{v})\mathrm{rk}(Q^i)-\mathrm{rk}(\mathbf{v})\mathrm{deg}(Q^i)$ is an integer,  the choice of $m_{0}$ depends only on $\mathbf{v}$, $\CO_S(1)$ (via the ampleness requirement)  and the constant $\chi_0$. The latter in turn depends on $\beta$, $\mathbf{v}$ and $\CO_S(1)$ by Corollary \ref{bounds} and bounds presented above,  assuming $S$ is fixed.  We conclude that $m_0$ depends only on $\beta$, $\mathbf{v}$ and $\CO_S(1)$.  

By (\ref{pos1}), the choice of $m_0$ and (\ref{boundq}), for all $m\geq m_0$ and for all $i$, we have 
\[ m\cdot (\mathrm{deg}(\mathbf{v})\mathrm{rk}(Q^i)-\mathrm{rk}(\mathbf{v})\mathrm{deg}(Q^i))+\chi(\mathbf{v})\mathrm{rk}(Q^i)-\mathrm{rk}(\mathbf{v})\chi(Q^i)>0.\]
Hence by (\ref{eq11}) and Theorem \ref{positive1}, for all $m\geq m_0$, we have
 \[ \CL_{0}\otimes \CL_{1}^{m} \cdot_{f} C>0\]
 for all non-constant prestable quasimaps of degree $\beta$. 
  By Lemma \ref{nonnegative} and Lemma \ref{bounds}, the same bounds hold for the restriction of $f$ to a subcurve  $C'\subseteq C$, since $\CL_{1}\cdot_{f_{|C'}}  C'\leq \CL_{1}\cdot_{f}  C$. Hence, for all $m\geq m_0$,  
  \[\CL_{0}\otimes \CL_{1}^{m}\cdot_{f_{|C'}}  C'>0,\]
   if and only if $f_{|C'}$ is not constant. \qed

\subsection{Stable quasimaps} \label{stabilisation}

For all $\beta \in \Eff(\CM(\mathbf{v}), \Coh(S,\mathbf{v}))$, we fix once and for ever a line bundle,
\[\CL_{\beta}:=\CL_0\otimes \CL_1^{m} \in \Pic_{\lambda}(\Coh(S,\mathbf{v})),\]
 such that it satisfies  conclusions of Proposition \ref{positive}.

Given a quasimap $f\colon C \rightarrow  \Coh(S,\mathbf{v})$ of degree $\beta$ and a point $p \in C$ in the regular locus of $C$. By Langton's semistable reduction, we can modify the quasimap $f$ at the point $p$ to obtain a quasimap  
\[f_p \colon C \rightarrow  \Coh(S,\mathbf{v}),\]
which maps to the stable locus $\CM(\mathbf{v})$ at $p$ and is isomorphic to $f$ away from $p$. If $p$ is not a base point, then $f_p=f$.   In other words, since $\CM(\mathbf{v})$ satisfies the existence part of the discrete valuation criterion of properness (while $M(\mathbf{v})$ is proper), and $C$ is spectrum of a discrete valuation ring at $p$, one can eliminate the indeterminacy of $f$ at $p$, if we view it as a rational map to $\CM(\mathbf{v})$. We refer to $f_p$ as the \textit{stabilisation} of $f$ at $p$.  

Langton's semistable reduction gives a canonical choice of such stabilisation, however, it is unique only for quasimaps to  $M(\mathbf{v})$. Otherwise, we can always tensor a sheaf with a line bundle from a curve which is trivial away from the point $p$. Nevertheless, any choice of stabilisation suffices for our purposes, as intersection numbers with $\CL_0$ and $\CL_1$ are independent of the choice, as they are invariant with respect to tensoring sheaves with line bundles from curves. 

\begin{defn} \label{length}
	Given a prestable quasimap $f \colon C \rightarrow  \Coh(S,\mathbf{v})$ of degree $\beta$, we define the \textit{length} of a point $p \in C$  with respect to $f$, 
	\[\ell(p):=\CL_{\beta} \cdot_{f} C-\CL_{\beta} \cdot_{f_{p}} C.\]
	By the proof of Proposition \ref{positive},  $\ell(p)\geq 0$; and $\ell(p)=0$, if and only if $p$ is not a base point. 
\end{defn}
In what follows, by $0^+$ we will denote a number $A \in \BR_{>0}$, such that $A \ll 1$. 
\begin{defn} \label{stabilityqm} Given $\epsilon \in \BR_{>0}\cup\{0^+, \infty\}$, a prestable quasimap $f\colon(C, \mathbf{p}) \rightarrow  \mathfrak{X}$ of degree $\beta$ is $\epsilon$-stable, if 
	\begin{itemize}
		\item[$\mathbf{(i)}$] $\omega_C(\mathbf{p}) \otimes f^*\CL_{\beta}^\epsilon$ is positive,
		\item[$\mathbf{(ii)}$] $\epsilon \ell(p)\leq 1$ for all $p \in C$.
	\end{itemize}
	We will refer to $0^+$-stable and $\infty$-stable quasimaps as stable quasimaps and stable maps respectively. 
	
\end{defn} A \textit{family} of $\epsilon$-stable quasimap over a base $B$ is a triple
\[(\CC,  \mathbf{p}, f),\]
 consisting of
 a family of marked nodal curves $(\CC, \mathbf{p})$ over $B$ and a map $f\colon (\CC,\mathbf{p}) \rightarrow \mathfrak{X}$, such that the geometric fibers of $f$ over $B$ are $\epsilon$-stable quasimaps.  An isomorphism of triples, 
  \[(g_1,g_2)\colon (\CC,  \mathbf{p}, f) \xrightarrow{\sim} (\CC',  \mathbf{p}', f'),\] is given by an isomorphism $g_1 \colon (\CC, \mathbf{p}) \xrightarrow{\sim} (\CC',  \mathbf{p}')$ together with an isomorphism $g_2 \colon f'\xrightarrow{\sim} f\circ \phi_1$. 
\begin{defn}
Let
\begin{align*}
	Q^{\epsilon}_{g,N}(M(\mathbf{v}), \beta) \colon (Sch/ \BC)^{\circ} &\rightarrow Grpd \\
	B&\mapsto \{\text{families of $\epsilon$-stable quasimaps over $B$}\}
\end{align*}
be the moduli space of $\epsilon$-stable quasimaps of genus $g$ and degree $\beta$ with $N$ marked points  to a pair  $(M(\mathbf{v}),\rCoh(S,\mathbf{v}))$. 
\end{defn}
\subsection{Relative moduli spaces of sheaves} \label{sectionrel}
There is a different modular interpretation of $Q^{\epsilon}_{g,N}(M(\mathbf{v}), \beta)$, it can be viewed as a space that parametrises sheaves $F$ on threefolds $S\times C$ associated to quasimaps $f \colon C \rightarrow \rCoh(S,\mathbf{v})$.  In this section we will make this precise.   

By the discussion in Section \ref{sectiondet}, the class $\mathbf{u}$ provides a section of $\tau$,
\begin{equation*}
	s_{\mathbf{u}}\colon \rCoh(S,\mathbf{v}) \rightarrow \Coh(S,\mathbf{v}),
\end{equation*}
which is induced by the descend of the twisted universal family \[\CF\otimes p_{\Coh(S,\mathbf{v})}^*\lambda(\mathbf{u})^{-1}\] from $S\times \Coh(S,\mathbf{v})$ to $S\times \rCoh(S,\mathbf{v})$.  Sheaves associated to quasimaps to $\Coh(S,\mathbf{v})$ obtained by this lift can be explicitly characterized.  
By the projection formula, there exists a canonical identification, 
\[\phi \colon \det(p_{\Coh *}(p_{S}^{*}\mathbf{u}\otimes \CF\otimes p_{\Coh}^*\lambda(\mathbf{u})^{-1})\xrightarrow{\sim } \lambda(\mathbf{u})\otimes \lambda(\mathbf{u})^{-1}\xrightarrow{\sim}\CO_{\Coh(S,\mathbf{v})},\]
which also descends to $\rCoh(S,\mathbf{v})$.  
We obtain that for the sheaf $F$ associated to the quasimap $s_{\mathbf{u}}\circ f$, there exists a canonical identification
\[\phi \colon \det(p_{C*}(p_{S}^{*}\mathbf{u}\otimes F))\xrightarrow{\sim }\CO_{C}.\]
This motivates the following definition. 

\begin{defn} Let $(C, \mathbf{p})$ a marked nodal curve. Let $F$ be a sheaf on $S\times C$ flat over $C$, such that $\ch(F)=(\ch(\mathbf{v}),\check{\beta})$. Given $\epsilon \in \BR_{>0}\cup\{0^+, \infty\}$, we say that $F$ is $\epsilon$-stable, if
	\begin{itemize}
	\item[$\mathbf{(i)}$] the associated quasimap $f \colon (C,\mathbf{p}) \rightarrow \Coh(S,\mathbf{v})$ is $\epsilon$-stable,
	\item[$\mathbf{(ii)}$] $\det(p_{C*}(p_{S}^{*}\mathbf{u}\otimes F))\cong \CO_{C}$. 
	\end{itemize}
	\end{defn}
 A \textit{family} of $\epsilon$-stable sheaves over a base $B$ is a quadruple
 \[(\CC,  \mathbf{p}, F, \phi),\] 
consisting of a family of marked nodal curves $(\CC, \mathbf{p})$ over $B$,  a sheaf $F$ on $S\times \CC$ flat over $B$, and an identification $\phi\colon \det(p_{\CC*}(p_{S}^{*}\mathbf{u}\otimes F))\xrightarrow{\sim }\CO_{\CC}$, such that the geometric fibers of $F$ over $B$ are $\epsilon$-stable sheaves. The notion of isomorphisms for these quadruples  is defined  as for quasimaps - it involves isomorphisms of all parts of quadruples. 
\begin{defn} \label{sheavesqms} Let  
	\begin{align*}
		M_{\mathbf{v},\check{\beta}}^{\epsilon}(S\times C_{g,N}) \colon (Sch/ \BC)^{\circ} &\rightarrow Grpd \\
		B&\mapsto \{\text{families of $\epsilon$-stable sheaves over $B$}\}
	\end{align*}
be the moduli space of $\epsilon$-stable sheaves with the Chern character $(\ch(\mathbf{v}),\check{\beta})$, such  that curves are of genus $g$ with $N$ markings. 
	\end{defn}
\begin{rmk} \label{modulisheaves} If $\epsilon=0^+$, then $\BC$-valued points of $M_{\mathbf v, \check{\beta}}^{0^+}(S\times C_{g,N})$ are triples  $(C, \mathbf{p}, F)$, such that: 
	\begin{itemize}
		\item $(C, \mathbf{p})$ is a prestable nodal curve (no rational tails),
		\item a sheaf $F$ on $S\times C$ is flat over $C$,
		\item$ \ch(F)=(\ch(\mathbf v), \check \beta)$,
		\item $F_p$ is stable for a general $p\in C$,
		\item  $F_p$ is stable, if $p$ is a node or a marking,
		\item the group $\{g \in \mathrm{Aut}(C, \mathbf{p}) \mid g^*F\cong F\}$ is finite,
	
		\item $\det(p_{C*}(p_{S}^{*}\mathbf{u}\otimes F))\cong\CO_{C}$.
	\end{itemize}
In the appendix, we will show that slope stability of a general fiber of $F$ implies slope stability of $F$ for a suitable polarisation, Corollary \ref{stability}. The converse is shown in Corollary \ref{converse},  if $\rk(\mathbf{v})=2$. 

For a general value of $\epsilon \in \BR_{>0}$, $\epsilon$-stability controls stability of fibers of a sheaf $F$ on $S\times C$. We refer to Corollary \ref{Hilb} for the precise statement in the case $\rk(\mathbf{v})=1$. In the case of $\rk(\mathbf{v})\geq 2$,  an expression of $\epsilon$-stability in purely sheaf-theoretic terms seems more difficult to state.

Among other applications, these spaces provide a way of defining higher-rank Donaldson--Thomas theory relative to vertical divisors on threefolds of the form $S\times C$, which can be used to derive higher-rank degeneration formulas. 

	Our determinant-line-bundle condition is natural for families. The standard determinant-line-bundle condition would involve a choice of a line bundle for families which might not even exist. For a fixed smooth curve, the two determinant-line-bundles conditions are not far from each other, as is shown in Lemma \ref{det}.
\end{rmk}

By design, there exist two natural transformations between 2-functors $Q^{\epsilon}_{g,N}(M(\mathbf{v}), \beta)$ and  $M_{\mathbf{v},\check{\beta}}^{\epsilon}(S\times C_{g,N})$.  In one direction the transformation  is given by composing quasimaps with the section $s_{\mathbf{u}}$ and then associating sheaves on threefolds, 
 \begin{align*} p \colon Q^{\epsilon}_{g,N}(M(\mathbf{v}), \beta) &\rightarrow M_{\mathbf{v},\check{\beta}}^{\epsilon}(S\times C_{g,N}), \\
 	(\CC, \mathbf{p}, f) \mapsto (\CC,\mathbf{p}, &s_{\mathbf{u}}\circ f) \mapsto (\CC, \mathbf{p}, F, \phi ),
 	\end{align*}
 where for the second association we used Lemma \ref{flatness}.  In the opposite direction the transformation is given by associating quasimaps with the target $\Coh(S,\mathbf{v})$ to sheaves on threefolds and then composing with the  projection $\tau$, 
 \begin{align*} q \colon  M_{\mathbf{v},\check{\beta}}^{\epsilon}(S\times C_{g,N})  &\rightarrow  Q^{\epsilon}_{g,N}(M(\mathbf{v}), \beta), \\
	(\CC, \mathbf{p}, F, \phi ) \mapsto (\CC,&\mathbf{p}, f) \mapsto (\CC, \mathbf{p}, \tau \circ f),
	\end{align*}
where for the first association we again used  Lemma \ref{flatness}. 
\begin{thm} \label{mapssheaves} The transformation $p$ is an equivalence of  2-functors, 
	\[ p \colon Q^{\epsilon}_{g,N}(M(\mathbf{v}), \beta) \xrightarrow{\sim} M_{\mathbf{v},\check{\beta}}^{\epsilon}(S\times C_{g,N}),\]
	such that the inverse is $q$. 
	\end{thm}
\textit{Proof.} By (\ref{sec}), we have an isomorphism
\[(\tau,\lambda(\mathbf{u}))\colon \Coh(S,\mathbf{v})\xrightarrow{\sim} \rCoh(S,\mathbf{v})\times B\BC^*,\]
which identifies maps from a base scheme $B$ to  $\rCoh(S,\mathbf{v})$, i.e.\ the $B$-valued points,  with maps from $B$ to $\Coh(S,\mathbf{v})$  together with  a trivialisation $\phi\colon \lambda(\mathbf{u})_{|B}\xrightarrow{\sim} \CO_{B}$. By construction of $\Coh(S,\mathbf{v})$ and $\lambda(\mathbf{u})$,  maps from $B$ to $\Coh(S,\mathbf{v})$ are given by sheaves $F$ on $S\times B$ and $\lambda(\mathbf{u})_{|B} =\det(p_{B*}(p_{S}^{*}\mathbf{u}\otimes F))$. Applying this to families of curves $\CC$ over a base scheme $B'$, and using Lemma \ref{flatness} to match flatness of sheaves over $\CC$ with flatness over $B'$, we obtain that $p$ and $q$ must be inverses of each other. 
\qed 

\begin{lemma} \label{flatness}Let $\CC$ be a family of nodal curves over a base $B$ and $F$ be a sheaf on $S\times \CC$ flat over $B$. Then the fiber sheaf $F_b$ on $S\times \CC_b$ is flat over $\CC_b$ for all closed points $b \in B$, if and only if $F$ is flat over $\CC$. 
	\end{lemma}

\textit{Proof}. By the local criteria of flatness, a sheaf $F$ is flat over a closed point $p \in \CC$, if and only if 
\[L\iota^*_pF= \iota_p^*F,\]
where $\iota_p \colon S\times \{p\} \hookrightarrow S\times \CC$ is the natural inclusion. The inclusion $\iota_p$ can be factored as follows
\[ \iota_p=\iota_b \circ \iota'_p \colon S\times \{p\}  \xhookrightarrow{\iota'_p}  S\times \CC_b \xhookrightarrow{\iota_b} S\times \CC,\]
where $\CC_b$ is the fiber in which $p$ is contained, hence  
\[L\iota^*_pF= L\iota'^*_p (L\iota^*_bF).\]
Since $\CC$ and $F$ are flat over $B$, we conclude that $F_b$ is flat over $\CC_b$ for all closed points $b\in B$, if and only if $F$ is flat over $\CC$. 
\qed 
\section{Algebraicity and properness} \label{Sectionalg}


\subsection{Regularity of sheaves} \label{sheaves}
We recall the definition of $m$-regularity of sheaves, which will be necessary for proving quasi-compactness of our moduli spaces.  
\begin{defn} Fix a very ample line bundle $\CO_X(1)$ on a variety $X$.  Let $m$ be an integer, a sheaf $F$ on $X$ is $m$-regular, if 
	\[H^i(X, F(m-i))=0 \quad \text{for all }i>0.\]
	The Mumford--Castelnuovo regularity of a sheaf $F$ is a number 
	\[\mathrm{reg}(F):=\mathrm{inf}\{m \in \BZ \mid F \text{ is }  m\text{-regular}\}.\]
	\end{defn}

 By \cite[Lemma 1.7.6]{HL}, bounds on the Mumford--Castelnuovo regularity allow to conclude boundedness of families of sheaves. We plan to use this property to prove quasi-compactness of moduli spaces $\epsilon$-stable sheaves. However, we firstly need to show that the family of possible curves $C$ is bounded. This follows directly from Proposition \ref{positive}.

\begin{lemma}  \label{projection} The forgetful map to the moduli stack of marked nodal curves, 
	\[Q^{\epsilon}_{g,N}(M(\mathbf{v}),\beta)) \rightarrow \FM_{g,N},\]
factors through a substack of finite type. 
\end{lemma}
\textit{Proof.}	 The restriction of a stable quasimap to an unstable component (a rational bridge or a rational tail) must be non-constant by stability and it must pair positively with $\CL_{\beta}$ by Proposition \ref{positive}. Therefore the number of unstable components of the domain curve of a stable quasimap is bounded in terms of $\beta$. Hence the projection $Q^{\epsilon}_{g,N}(M(\mathbf{v}),\beta) \rightarrow \FM_{g,N}$ factors through a substack of finite type. 
\qed
\\

 By Proposition \ref{positive}, for each stable quasimap $f\colon C \rightarrow \rCoh(S,\mathbf{v})$ of a fixed degree $\beta$, we obtain an ample bundle $f^* \CL_\beta$ on $C$.  By Lemma \ref{projection},   for a big enough $k \gg 0$, the line bundle 
\[\CO_{S\times C}(1):= \CO_S(1)\boxtimes f^* \CL_\beta^k \]
is very ample on $S\times C$ for all $f\colon C \rightarrow \rCoh(S,\mathbf{v})$. In the next lemma, regularity and Hilbert polynomials of $\epsilon$-stable sheaves on $S\times C$ are defined with respect to $\CO_{S\times C}(1)$.  



\begin{lemma}\label{betac}
	Fix a class $\beta \in \Eff(M(\mathbf{v}), \rCoh(S,\mathbf{v}))$ and $\epsilon\in \BR_{>0}$.
	  Then there exists $A \in \BZ$, such that $\epsilon$-stable sheaves with the Chern character $(\ch(\mathbf{v}),\check{\beta})$ satisfy	
	\[\mathrm{reg}(F)\leq A,\]
	and the set of Hilbert polynomials of such sheaves is finite. 
	\end{lemma}
\textit{Proof.} 
 Let $F^0$ be an $\epsilon$-stable sheaf with the associated quasimap  $ f \colon C \rightarrow \Coh(S,\mathbf{v})$. The semistable reduction applied to all base points at once gives a sequence of short exact sequences,
\begin{align*}
	0 \rightarrow F^{1}\rightarrow &F^0 \rightarrow Q^{1} \rightarrow  0, \\
	&\vdots \\
	0\rightarrow F^{k} \rightarrow &F^{k-1} \rightarrow Q^{k} \rightarrow 0, 
\end{align*}
such that $F^k$ defines a map $f^k\colon C \rightarrow M(\mathbf{v})$. By the associated long exact sequence of cohomologies, $F^0$ is $m$-regular, if and only if $F^k$ and $Q^i$ are $m$-regular. Using the fact that a family of sheaves is bounded, if and only if its Mumford--Castelnuovo regularity is bounded and the set of Hilbert polynomials is finite \cite[Lemma 1.7.6]{HL}, we have to show:
\begin{itemize}
	\item[\textbf{(i)}] the number of steps $k$ is bounded,
	\item[\textbf{(ii)}] the family of possible $Q^i$ is bounded,
	\item[\textbf{(iii)}] the family of possible $F^k$ is bounded. \\
\end{itemize}

\noindent \textbf{(i)} By the proof of Proposition \ref{positive}, at each step of Langton's semistable reduction, the difference  between $\CL_\beta\cdot_{f^i}C$ and $\CL_\beta\cdot_{f^{i-1}}C$  is a strictly positive integer, and $\CL_\beta\cdot_{f^k}C\geq0$. Hence the number of steps $k$ must be bounded by the intersection number of $C$ with the line bundle $\CL_\beta$,
\[k\leq  \CL_\beta \cdot_f C=\beta(\CL_\beta).\]

\noindent \textbf{(ii)} We will show that the family of $Q^i$ is bounded by bounding the Hilbert polynomial and using the bound on $\mu_{\mathrm{max}}(Q^i)$. By the Grothendieck--Riemann--Roch theorem, to bound the Hilbert polynomial, we need to bound the rank, the degree and the Euler characteristics of $Q^i$.

  The rank of $Q^i$ is bounded by $\rk(\mathbf{v})$. The degree of $Q^{i}$ is bounded by Lemma \ref{bounds}. The Euler characteristics is bounded from above by the bounds from the proof of Proposition \ref{positive} and by Lemma \ref{bounds}. The bound from below for the Euler characteristics also follows from  Proposition \ref{positive} and its proof as we will now demonstrate. Firstly, by $\mathbf{(i)}$, the number of steps $k$ is bounded. 
 Moreover, the family of possible $f^k \colon C \rightarrow M(\mathbf{v})$ is also bounded, since its degree with respect to an ample line bundle $\CL_\beta$ is bounded, and the moduli space of maps from a bounded family of curves of bounded degree is bounded\footnote{This follows from the boundedness of Hilbert schemes of subschemes with bounded Hilbert polynomials.} (boundedness of $F^k$ requires an extra argument, but for this part boundedness of $f^k$ suffices, because intersection numbers with $\CL_0$ and $\CL_1$ are invariant with respect to tensoring families of sheaves with line bundles, see the proof of the part (\textbf{iii})).  This implies that $\CL_0 \cdot_{f^k} C$ is also bounded. By Lemma \ref{key}, we have the following expression for the difference between intersection numbers with the line bundle $\CL_0$, 
\begin{equation} \label{diff}
\CL_0 \cdot_{f^i} C-\CL_0 \cdot_{f^{i-1}} C=\chi(\mathbf{v})\mathrm{rk}(Q^i)-\mathrm{rk}(\mathbf{v})\chi(Q^i),
\end{equation}
  it must be bounded from below, because $\chi(Q^i)$ is bounded from above.   We conclude that the difference (\ref{diff}) must be also bounded from above, because $\beta(\CL_0)$ is fixed, there are bounded number of steps, and $\CL_0 \cdot_{f^k} C$ is bounded. This bounds the Euler characteristics of $Q^i$ from below.  We thereby bounded the Hilbert polynomial of $Q^i$. Finally, $\mu_{\mathrm{max}}(Q^i)$ is bounded by Lemma \ref{bounds}. Using \cite[Theorem 3.3.7]{HL}, we conclude that the family of possible $Q^i$ is bounded.
\\

\noindent \textbf{(iii)} For this part, we need the boundedness of corresponding maps $f^k\colon C \rightarrow M(\mathbf{v})$ and the boundedness of line bundles $L$, such that 
\[ \det(p_{C*}(p_{S}^{*}\mathbf{u}\otimes F^k))\cong L,\]
the latter is needed because the semistable reduction does not preserve the condition $\det(p_{C*}(p_{S}^{*}\mathbf{u}\otimes F))\cong \CO_{C}$. On the other hand, the sheaf $F^k\otimes p_C^*L^{-1}$ satisfies $\det(p_{C*}(p_{S}^{*}\mathbf{u}\otimes F^k\otimes p^*_CL^{-1}))\cong \CO_{C}$, hence the associated quasimap is a lift of a quasimap from $M(\mathbf{v})$ by the section $s_{\mathbf{u}}$. We refer to Remark \ref{liftsrem} for more on how different lifts are related. 

  The boundedness of $f^k\colon C \rightarrow M(\mathbf{v})$ was established in the proof of the part $(\mathbf{ii})$ above. It remains to  deal with line bundles $L$. By the exact sequences above and the boundedness of $Q^i$ from the part $(\mathbf{ii})$, the boundedness of possible $L$ follows from the boundedness of the number of steps $k$, which was proved in $(\mathbf{i})$. 
  
  Here, we crucially rely on the fact that our target is the rigidified stack $\rCoh(S,\mathbf{v})$, which means that  we start with a sheaf that satisfies the condition $\det(p_{C*}(p_{S}^{*}\mathbf{u}\otimes F^0))\cong \CO_C$. Otherwise, the line bundle   $\det(p_{C*}(p_{S}^{*}\mathbf{u}\otimes F^0))$ can be arbitrary, the family of possible $L$ is therefore unbounded. \\

  Overall, we conclude that the number of steps $k$ is bounded, and the families of possible $Q^i$ and $F^k$ are also bounded. This implies that the Mumford--Castelnuovo regularity of $\epsilon$-stable  sheaves with the given Chern character is bounded by some number $A$, 
  \[\mathrm{reg}(F)\leq A, \]
  and the set of Hilbert polynomials of such sheaves is finite. 
\qed

\subsection{Algebraicity and quasi-compactness}

We will show that the condition of $\epsilon$-stability is open in the stack of all maps  to $\Coh_r(S,\mathbf{v})$. The difficult part is to prove that the condition $\mathbf{(ii)}$ of $\epsilon$-stability is open. The reason is that it is not well-adapted for families, as it requires stabilisations of quasimaps. We therefore firstly show that it is locally  constructible and then use the valuative criterion\footnote{We refer to the result which states that (locally) constructable sets are open, if and only if they are stable under generalisations of points, see \cite[Section 0060]{stacks-project}; stability of generalisations of points can be check on discrete valuation rings. } of openess for locally constructible subsets. 

Let \[ \mathrm{Map}_{\FM_{g,N}}(\FC_{g,N}, \Coh_r(S,\mathbf{v})\times \FM_{g,N}) \colon (Sch/ \BC)^{\circ} \rightarrow Grpd \] be the stack of all maps from nodal curves to $\Coh_r(S,\mathbf{v})$.  Under our assumptions, we have
\[\Coh(S,\mathbf{v})\cong \rCoh(S,\mathbf{v})\times B\BC^*,\]
hence $\rCoh(S,\mathbf{v})$ must have affine stabilizers (see \cite{HR} for the definition), because $\Coh(S,\mathbf{v})$  have affine stabilizers. Moreover, $\Coh(S,\mathbf{v})$ is a quasi-separated algebraic stack locally of finite presentation by \cite[Section 0DLX]{stacks-project}. Hence the same applies to $\rCoh(S,\mathbf{v})$. By the representability theorem for mapping stacks of \cite{HR}, we  obtain that  $\mathrm{Map}_{\FM_{g,N}}(\FC_{g,N}, \Coh_r(S,\mathbf{v})\times \FM_{g,N})$ is a quasi-separated algebraic  stack locally of finite presentation over $\FM_{g,N}$.  

\begin{lemma} \label{const} The space $Q^{\epsilon}_{g,N}(M(\mathbf{v}),\beta)$ is locally constructible in the mapping stack  $\mathrm{Map}_{\FM_{g,N}}(\FC_{g,N}, \Coh_r(S,\mathbf{v})\times \FM_{g,N})$. 
	\end{lemma}
\textit{Proof.} The condition of prestability of quasimaps can be easily seen to be open. The condition $\mathbf{(i)}$ of $\epsilon$-stability is also open, since  ampleness of line bundles is open.  We therefore have to deal with the condition $\mathbf{(ii)}$.  Let 
\[\CU \subset \mathrm{Map}_{\FM_{g,N}}(\FC_{g,N}, \Coh_r(S,\mathbf{v})\times \FM_{g,N}) \] be an open quasi-compact substack.  We may assume that $\CU$ contains only prestable quasimaps. Let $\CC_{\CU}$ be the universal curve over $\CU$ and let
\[ \mathbb{f} \colon \CC_{\CU} \rightarrow \rCoh(S,\mathbf{v})\] 
be the universal map. We want to stratify the space $\CU$ in the way that $\mathbb{f}$ can be stabilized at base points to be able to evaluate their length. Firstly, we stratify $\CU$ according to the number of base points. Namely, let $\CU_k$ be moduli space of maps from prestable nodal curves  to $\rCoh(S,\mathbf{v})$ with $N+k$ markings, such that:
\begin{itemize} 
	\item there are exactly $k$ base points,
	\item all $k$ base points are marked by the last $k$ markings,
	\item the underlying map with $N$ markings is contained in $\CU$. 
	\end{itemize} 
These are locally closed substacks of mapping stacks with $N+k$ markings. Since $\CU$ is of finite type, the number of base points is bounded from below, hence a finite disjoint union $\bigcup_{k} \CU_k$  covers $\CU$, i.e.\ there exists a surjection,
\[p\colon \bigcup_{k} \CU_k \rightarrow \CU,\]
given by forgetting the last $k$ makings. 
We now further stratify $\CU_k$. Consider the universal maps, 
\[ \mathbb{f}_k \colon \CC_{\CU_k} \rightarrow \rCoh(S,\mathbf{v}),\]
since all base points are contained in the regular locus of the universal curve, the maps $\mathbb{f}_k$ can be stabilized at base points over the generic points of irreducible components of $\CU_k$ (we can reduce the components, if necessary,  to talk about generic points, as this does not affect the topology). Note that stabilisation just removes indeterminacies of  the associated rational map 
\[ \mathbb{f}_k \colon \CC_{\CU_k} \dashrightarrow M(\mathbf{v}).\]
 In particular, there exists a dense open subset of $\CU_k$, such that the quasimap $\mathbb{f}_k$ can be stabilised at $i^{\mathrm{th}}$ base point. By passing  to the complement of that open subset, repeating the argument and using the fact that $\CU_k$ is quasicompact, we obtain a finite stratification of $\CU_k$ into locally closed subsets, 
\[\bigcup_{j,i} \CU^j_{k,i}=\CU_k, \] 
such that the quasimap map $\mathbb{f}_k$ can be stabilized at  $i^{\mathrm{th}}$ base point over each stratum $\CU^j_{k,i}$. We now can evaluate the length of base points at each stratum  $\CU^j_{k,i}$. The condition $\mathbf{(ii)}$ of $\epsilon$-stability involves bounding a degree of line bundles from above, hence this condition just picks some connected components (or non at all) of the strata  $\CU^j_{k,i}$, which we denote by $\CU^{j'}_{k,i}$. Since $\CU^j_{k,i}$ are locally closed subsets, $\CU^{j'}_{k,i}$ are also locally closed and, in particular, constructible. The locus of $\epsilon$-stable quasimaps $\CU^{\epsilon}_k$ in $\CU_k$ is given by intersections of $\cup_{j'}\CU^{j'}_{k,i}$ for all $i$, 
\[\CU^{\epsilon}_k :=\bigcap_i \left( \bigcup_{j'}\CU^{j'}_{k,i}\right), \] 
which is constructible, since it involves  finite intersections and finite unions. 
The locus of $\epsilon$-stable quasimaps in $\CU$ is the image of the union $\cup_k\CU^{\epsilon}_k$, i.e.\ on the level of geometric points, we have 
\[|Q^{\epsilon}_{g,N}(M(\mathbf{v}),\beta)\cap \CU|=\lvert p\left(\bigcup_k\CU^{\epsilon}_k\right)\rvert.\]
Since $p$ is a map between stacks of finite type, the image of a constructible set is constructible. We conclude that $Q^{\epsilon}_{g,N}(M(\mathbf{v}),\beta)$ is locally constructible. \qed

\begin{prop} \label{openess}  The space $Q^{\epsilon}_{g,N}(M(\mathbf{v}),\beta)$ is open in the mapping stack $\mathrm{Map}_{\FM_{g,N}}(\FC_{g,N}, \Coh_r(S,\mathbf{v})\times \FM_{g,N})$.
\end{prop}
\textit{Proof.} As we said before, the only non-trivial part of the claim is the condition $\mathbf{(ii)}$ of $\epsilon$-stability. By Lemma \ref{const}, we can use the valuative criterion of openess for  (locally) constructible sets.  We therefore have to show that if a family of quasimaps $f \colon \CC \rightarrow \rCoh(S,\mathbf{v})$ over the spectrum of a discrete valuation ring $\Delta$ does not satisfy the condition $\mathbf{(ii)}$ at the generic fiber, it also does not satisfy it at the special fiber.

 If the generic fiber $f^\circ \colon \CC^\circ \rightarrow \rCoh(S,\mathbf{v})$ does not satisfy the condition $\mathbf{(ii)}$, then there exists a base points $b^\circ \in \CC^\circ$ which violates it, we can stabilize it  to obtain another family $f'^\circ \colon \CC^\circ \rightarrow \rCoh(S,\mathbf{v})$ over $\Delta^\circ$, which does not have a base point at $b^\circ$. It extends to the family $f' \colon \CC \rightarrow \rCoh(S,\mathbf{v})$ over $\Delta$ by Lemma \ref{extension0}. By definition, the length of the base point $b^\circ$  is given as follows,
\[\ell(b^\circ)=\CL_{\beta} \cdot_{f^\circ} \CC^\circ-\CL_{\beta} \cdot_{f'^\circ} \CC^\circ.\]
Let $b_0 \in \CC_0$ be the limit of $b^\circ$ in the central fiber $\CC_0$, then 
\[\ell(b_0)\geq \CL_{\beta} \cdot_{f_0} \CC_0-\CL_{\beta} \cdot_{f'_0} \CC_0=\CL_{\beta} \cdot_{f^\circ} \CC^\circ-\CL_{\beta} \cdot_{f'^\circ} \CC^\circ,\]
we use the sign of inequality, because $b_0$ might still be a base point of $f'_0$ whose length with respect to $f'_0$ is smaller than the length with respect to $f_0$.   We therefore conclude that $f_0$ also does not satisfy the condition $\mathbf{(ii)}$. This implies that the complement of  $Q^{\epsilon}_{g,N}(M(\mathbf{v}), \beta)$ is closed, hence $Q^{\epsilon}_{g,N}(M(\mathbf{v}), \beta)$ is open.  
\qed

\begin{cor}  \label{algebraicity} The space
	$Q^{\epsilon}_{g,N}(M(\mathbf{v}), \beta)$ is a quasi-separated algebraic stack of finite presentation. 
	\end{cor}
\textit{Proof}. Firstly, $Q^{\epsilon}_{g,N}(M(\mathbf{v}), \beta)$ is algebraic, quasi-separated and locally of finite presentation by Proposition \ref{openess}, since mapping stacks have these properties. Secondly,  By Lemma \ref{betac}, Lemma \ref{projection} and \cite[Lemma 1.7.6]{HL},  the  stack $M_{\mathbf{v},\check{\beta}}^{\epsilon}(S\times C_{g,N})$  is quasi-compact.  Hence by  Proposition \ref{positive}, the stack $Q^{\epsilon}_{g,N}(M(\mathbf{v}), \beta)$ is also quasi-compact, and therefore of finite presentation. 
\qed

\subsection{Properness}

In this section, we prove that $Q^{\epsilon}_{g,N}(M(\mathbf{v}), \beta)$ is a proper Deligne--Mumford stack. To this end, we need a few auxiliary results about sheaves on $S\times C$.  The central among which is  \textit{Hartogs' property} for families of nodal curves over a spectrum of a discrete valuation ring $\Delta$, established in  Lemma \ref{extension0}.

\begin{lemma}\label{torsion}
	Let $F$ be a sheaf on $S \times C$ flat over $C$, such that $F_{p}$ is torsion free for a general $p \in C$, then $F$ is torsion free. Conversely, if $F$ is torsion free  and flat over nodes of $C$, then $F$ is flat over $C$ and $F_p$ is torsion free  for a general $p \in C$. 
\end{lemma} 

\textit{Proof.}
Assume $F$ is flat, such that $F_{p}$ is torsion free for a general $p \in C$. Let $T(F) \subset F$ be the maximal torsion subsheaf.  Firstly, $T(F)\neq F$, because $\rk(F)> 0$ (a general fiber of $F$ is torsion free and therefore of non-zero rank). It also cannot be supported on fibers of $S\times C \rightarrow C$ due to flatness of $F$ over $C$, therefore Supp$(T(F))$ intersects a general fiber. Since $F/T(F)$ is generically flat, restricting $T(F) \subset F$ to a general fiber $p\in C$, we get a torsion subsheaf of $F_{p}$ for a general $p \in C$, which is zero, therefore $T(F)$ is zero. The converse follows from the fact that torsion freeness implies flatness over a discrete valuation ring, and torsion in a general fiber $F_p$ would produce torsion in $F$. 
\qed

\begin{lemma}\label{locfree} Let $\CU$ be a regular surface,  let $F_1$ and $F_2$ be sheaves on $ S\times \CU$ flat over $\CU$. Then the sheaf $\CH om_{p_{\CU}} (F_1,F_2)$ is locally free for the natural projection $p_{\CU} \colon S\times \CU \rightarrow \CU$.
\end{lemma}

\begin{proof} It is enough to show that $\CH om_{p_\CU} (F_1,F_2)=p_{\CU*}\CH om(F_1,F_2)$ is reflexive, because being  reflexive is equivalent to being locally free on regular surfaces. 
	
	Since $F_1$ is flat over $\CU$, by the construction from the proof of \cite[Proposition 2.1.10]{HL}, there is an exact sequence 
	\[ H_1 \rightarrow H_0 \rightarrow F_1 \rightarrow 0,\]
	such that $H_i=p_{\CU}^*H_i' \otimes \CO_S(-k_i)$ for a locally free sheaf $H_i'$ on $\CU$ and for a sufficiently big integer $k_i$. We apply $\CH om_{p_{\CU}}(-,F_2)$ to obtain an exact sequence 
	\[ 0\rightarrow \CH om_{p_{\CU}}(F_1,F_2) \rightarrow \CH om_{p_{\CU}}(H_0,F_2) \rightarrow \CH om_{p_{\CU}}(H_1,F_2).\]
	By using the identification
	\[
	\CH om(H_i,F_2)=H_i^*\otimes F_2=p_{\CU}^*H_i'^* \otimes F_2 \otimes \CO_S(k_i)
	\]
	and choosing big enough $k_i$, the flatness of $F_2$ and  \cite[Proposition 2.1.2]{HL} imply that 
	\[ \CH om_{p_{\CU}}(H_i,F_2)=H_i'^* \otimes p_{\CU*} (F_2 \otimes \CO_S(k_i))   \]
	is a locally free sheaf. Hence by \cite[Lemma 0EB8]{stacks-project}, we obtain that $ \CH om_{p_{\CU}}(F_1,F_2)$ is reflexive. 
	
\end{proof}

\begin{lemma} \label{extension0}
	Let $\CC \rightarrow \Delta$ be a family of nodal projective curves over the spectrum of a discrete valuation ring, and let $\{p_1, \dots, p_m\}\subset \CC$ be a set of finitely many closed points in the regular locus of the central fiber. Then  any quasimap 
	$\tilde{f}\colon \tilde{\CC}=\CC \setminus \{p_1,\dots, p_m\} \rightarrow \rCoh(S,\mathbf{v})$ extends to $f\colon\CC \rightarrow \rCoh(S,\mathbf{v})$, which is unique up to a unique isomorphism.
\end{lemma}
\textit{Proof.} Let $\tilde{F}$ be the family on $S \times \tilde{\CC}$ corresponding to the lift of $\tilde{f}$ by $s_{\mathbf{u}}$, we then extend $\tilde{F}$ to a coherent sheaf $F$ on $S\times \CC$, quotienting our the torsion, if necessary. The sheaf $F$ is therefore flat over $\Delta$. Moreover, by Lemma \ref{flatness}, the sheaf $F$ is flat over $\CC$, if and only if the central fiber $F_{0}$ is flat over $\CC_0$. By Lemma \ref{torsion}, the central fiber $F_{0}$ is flat over $\CC_0$,  if and only if it is torsion free, because $\CC_{0}$ is regular at $p_i$. If $F_{0}$ is not torsion free, we can remove the torsion inductively as follows. Let $F^{0}=F$ and  $F^{i}$ be defined by short exact sequences, 
\[0 \rightarrow F^{i} \rightarrow F^{i-1} \rightarrow Q^{i} \rightarrow 0,\]
such that $Q^{i}$ is the quotient of $F^{i-1}_{0}$ by the maximal torsion subsheaf. At each step, the torsion of $F_{0}^{i}$ is supported at slices $S\times\{ p_{i}\}$, therefore all $F^i$ are isomorphic to $F^0$ over $S\times \tilde{\CC}$. By the argument from \cite[Theorem 2.B.1]{HL}, this process terminates, i.e.\ $F^i=F^{i+1}$ and $F^i_{0}$ is torsion free for $i\gg 0$. Let us  redefine the sheaf $F$, we set $F=F^i$ for some $i\gg0$, then the sheaf $F$ induces a quasimap to $\Coh(S,\mathbf{v})$, and composing it with the projection to $\rCoh(S,\mathbf{v})$, we thereby obtain an extension $f\colon \CC \rightarrow \rCoh(S,\mathbf{v})$ of $\tilde{f}$.

Consider now another extension $f'\colon \CC \rightarrow \rCoh(S,\mathbf{v})$, we lift both $f$ and $f'$ to $\Coh(S,\mathbf{v})$ with $s_{\mathbf{u}}$. Let $F$ and $F'$ be the corresponding families on $S\times \CC$, which satisfy 
\[ F_{|\tilde{\CC}}=F'_{|\tilde{\CC}}.\]
We will show that $F$ and $F'$ are in fact isomorphic, using the argument from \cite[Lemma 3.25]{Ku}. 
Let $\CU\subseteq \CC$ be a regular open neighbourhood of $\{p_1, \dots, p_m\}\subset \CC$. By Lemma \ref{locfree}, $\CH om_{p_{\CU}}(F,F')$ is locally free. Moreover, since $F$ and $F'$ are equal away from $\{p_1 ,\dots, p_m\}$ and a general fiber is simple, $\CH om_{p_{\CU}}(F,F')$ must be a line bundle. It has a non-vanishing section given by the identity morphism, defined away from $\{p_1 ,\dots , p_m\}$. Since $\{p_1 ,\dots , p_m\}$ is of codimension 2 and $\CH om_{p_{\CU}}(F,F')$ is locally free, the section extends to the whole $\CU$, providing a trivialisation $\CH om_{p_{\CU}}(F,F')\cong \CO_\CU$. Using the tautological morphism\footnote{It is defined as the composition $F_{|\CU} \otimes p_{\CU}^*(\CH om_{p_{\CU}}(F,F')) \rightarrow F_{|\CU} \otimes \CH om_{\CU}(F,F') \rightarrow  F'_{|\CU}$, where the first map is given by the adjunction.},
\[
F_{|\CU} \otimes p_{\CU}^*(\CH om_{p_{\CU}}(F,F')) \rightarrow F'_{|\CU},
\]
and the trivialisaiton above, we therefore obtain a morphism 
\[ F_{|\CU} \rightarrow F'_{|\CU},\]
which is an identity away  $\{p_1 ,\dots, p_m\}$. Gluing it with the identity morphism over $\tilde{\CC}$, we obtain a morphism defined over the entire $\CC$, 
\begin{equation} \label{themor}
 F \rightarrow F'.
 \end{equation} 
It is injective and its cokernel $Q$ is supported over $\{p_1 ,\dots, p_m\}$. By pulling back to the closed fiber $\CC_0$ of $\CC$, we obtain an exact sequence 
\[0 \rightarrow F_{|\CC_0} \rightarrow F'_{|\CC_0} \rightarrow Q_{|\CC_0}\rightarrow0,\]
note that it is exact from the left, because $F_{|\CC_0}$ is torsion free by Lemma \ref{torsion}. Since $F_{|\CC_0}$ and $F'_{|\CC_0}$ have the same Chern character (equal to the one of the generic fiber), we obtain that $Q$ must be zero, hence the morphism (\ref{themor}) is an isomorphism. 
Finally, for a fixed isomorphism of $F$ and $F'$ over $\tilde{\CC}$, the isomorphism over $\CC$ is  unique, as the isomorphism over $\tilde{\CC}$ determines a section of  $\CH om_{p_{\CU}}(F,F')$ away from $\{p_1 ,\dots, p_m\}$, which uniquely extends to the whole $\CU$. 
\qed

\begin{rmk} It is not clear if for a general smooth surface $\CU$ a map $f\colon \CU \setminus \{p_1, \dots,p_m\} \rightarrow \Coh_r(S,\mathbf{v})$ extends to the whole $\CU$, like in the case of GIT quotients in \cite[Lemma 4.3.2]{CFKM}. Hence the assumption that our surface is given by a family of  curves $\CC \rightarrow \Delta$ might be necessary. This form of Hartogs' property is good enough for proving Theorem \ref{proper} in the spirit of \cite[Section 4]{CFKM}.
\end{rmk}

\begin{lemma} \label{simple} Let $F$ be the sheaf on $S\times C$ associated to a quasimap $f \colon C \rightarrow \Coh_r(S,\mathbf{v})$, then $F$ is simple. 
\end{lemma}

\begin{proof} We have to show that 
	\[\Hom(F,F)=\BC.\]
	By passing to the normalisation of $C$, we may assume that $C$ is smooth. The sheaf  $\CH om_{p_{C}}(F,F)=p_{C*}\CH om(F,F)$ is torsion free on $C$ and a general fiber of $F$ over $C$ is simple, hence $\CH om_{p_{C}}(F,F)$ is a line bundle.   Moreover, it has a non-vanishing section given by the identity morphism $\CO_C \rightarrow   \CH om_{p_{C}}(F,F)$. We conclude that $\CH om_{p_{C}}(F,F)$ is a trivial line bundle, and therefore 
	\[\Hom(F,F)=H^0(S\times C, \CH om(F,F))=H^0(C, \CH om_{p_{C}}(F,F))=\BC, \]
	which concludes the proof. 
\end{proof}

\begin{thm}\label{proper} The space $Q^{\epsilon}_{g,N}(M(\mathbf{v}),\beta)$ is a proper Deligne--Mumford stack. 
\end{thm}
\textit{Proof.}  We use identification from Theorem \ref{mapssheaves}. By Corollary \ref{algebraicity} and Lemma \ref{simple}, the forgetful map
\[ M_{\mathbf{v},\check{\beta}}^{\epsilon}(S\times C_{g,N}) \rightarrow \FM_{g,N}\]
is representable by algebraic spaces. By $\epsilon$-stability, $\epsilon$-stable sheaves are fixed by at most a finite discrete subgroup of the automorphism group of a curve. We therefore conclude that $M_{\mathbf{v},\check{\beta}}^{\epsilon}(S\times C_{g,N})$ is a quasi-separated Deligne--Mumford stack of finite type. Using the valuative criteria of properness for quasi-separated Deligne--Mumford stacks of finite type \cite[Section  0CLY]{stacks-project} and Lemma \ref{extension0}, the proof of properness then proceeds as in the GIT case \cite[Proposition 4.3.1]{CFKM}.
\qed

\subsection{Sheaves of fixed determinant} \label{stablesheaves}

Assume $C$ is smooth.  We define 
\[M_{\mathbf{v},\check{\beta}}^{0^+}(S\times C)\]
 to be the fiber of $ M_{\mathbf{v},\check{\beta}}^{0^+}(S\times C_{g,0})$ over the curve $C \in \FM_{g,0}(\BC)$. Unpacking the definition of $\epsilon$-stability for sheaves, we obtain that the space $M_{\mathbf{v},\check{\beta}}^{0^+}(S\times C)$ parametrises sheaves $F$ on $S\times C$ subject to the following conditions:
\begin{itemize}
	\item $\ch(F)=(\ch (\mathbf v), \check \beta)$,
	\item $F$ is torsion free,
	\item a general fiber of $F$ is stable,
	\item the group $\{g \in \mathrm{Aut}(C)\mid g^*F\cong F\}$ is finite,
	\item $\det(p_{C*}(p_{S}^{*}\mathbf{u}\otimes F))\cong \CO_{C}$. 
	\end{itemize}

If  $g(C)\geq 1$ and $\beta\neq 0$, then the group of automorphisms of $C$ fixing $F$ is automatically finite.  All of the conditions above are standard except the last one, as one usually fixes the determinant of $F$. Let 
\[\widetilde{M}_{\mathbf{v},\check{\beta}}(S \times C)\] 
be the moduli space of sheaves subject to all conditions above, except the last one, instead we require the following:
\begin{itemize}
	\item $\det(F)\cong L$ for a fixed line bundle $L$, such that $\mathrm{c}_1(L)=\mathrm{c}_1(F)$.
	\end{itemize}
There exists a map that relates two moduli spaces, 
\begin{align*} 
	p\colon \widetilde M_{\mathbf{v},\check{\beta}}(S \times C) &\rightarrow M_{\mathbf{v},\check{\beta}}^{0^+}(S\times C), \\
	F &\mapsto F\otimes p_C^*\det(p_{C*}(p_{S}^{*}\mathbf{u}\otimes F))^{-1}.
\end{align*}
  It is \'etale and surjective, as is shown in the following lemma.
\begin{lemma} \label{det} The map $p$ is \'etale and surjective of degree $\rk(\mathbf{v})^{2g}$.
\end{lemma}

\textit{Proof.}
In fact, we will show that $p$ is a finite-group torsor.  Let  
\[\Gamma:=\Pic_0(C)[\rk(\mathbf{v})]\] be the $\rk(\mathbf{v})$-torsion points of $\Pic_0(C)$. The group $\Gamma$ acts on the moduli space $\widetilde M_{\mathbf{v},\check{\beta}}(S \times C)$ by tensoring with line bundles, because $\det(F\otimes A)=\det(F)\otimes A^{\otimes\rk(\mathbf{v})}$ for a line bundle $A \in \mathrm{Pic}(S\times C)$. The action is free, because  the following holds,
\[\det(p_{C*}(p_{S}^{*}\mathbf{u}\otimes F\otimes p_C^* A)\cong \det(p_{C*}(p_{S}^{*}\mathbf{u}\otimes F)\otimes A,\]
which is due to $\chi( \mathbf{v} \cdot u)=1$. The map $p$  is $\Gamma$-invariant.  In particular, we have the induced map 
\begin{equation} \label{isotor}
\widetilde M_{\mathbf{v},\check{\beta}}(S \times C) /\Gamma \rightarrow M_{\mathbf{v},\check{\beta}}^{0^+}(S\times C). 
\end{equation}

Conversely,  the inverse of this map is given by a $\Gamma$-torsor on $M_{\mathbf{v},\check{\beta}}^{0^+}(S\times C)$ which parametrises sheaves $F$  in  $M_{\mathbf{v},\check{\beta}}^{0^+}(S\times C)$ together with a $\rk(\mathbf{v})^{\mathrm{th}}$-root of $\det(F)\otimes L^{-1}$. Using the assumption $h^1(S)=0$ and the identification $\Pic_0(S\times C)\cong\Pic_0(C)$, this torsor is constructed as the fiber product, 

\begin{equation*}
	\begin{tikzcd}[row sep=scriptsize, column sep = large]
		& M_{\mathbf{v},\check{\beta}}^{0^+}(S\times C)^{\frac{1}{\rk(\mathbf{v})}} \arrow[d]  \arrow[r] &  \Pic_0(C) \arrow[d,"\cdot \rk(\mathbf{v})"]\\
		&  M_{\mathbf{v},\check{\beta}}^{0^+}(S\times C) \arrow[r,"\det(F)\otimes L^{-1}"] &   \Pic_0(C)
	\end{tikzcd}
\end{equation*}

There exists a $\Gamma$-equivariant  map, 
\[M_{\mathbf{v},\check{\beta}}^{0^+}(S\times C)^{\frac{1}{\rk(\mathbf{v})}}  \rightarrow \widetilde M_{\mathbf{v},\check{\beta}}(S \times C), \quad F \mapsto F\otimes(\det(F)^{-1}\otimes L)^{\frac{1}{\rk(\mathbf{v})}},\]
which induces  the inverse of (\ref{isotor}). We conclude
\[\widetilde M_{\mathbf{v},\check{\beta}}(S \times C) /\Gamma \cong M_{\mathbf{v},\check{\beta}}^{0^+}(S\times C),\]
in particular, $\widetilde M_{\mathbf{v},\check{\beta}}(S \times C) $   is a $\Gamma$-torsor over $M_{\mathbf{v},\check{\beta}}^{0^+}(S\times C)$. A finite-group torsor is \'etale over the base. The cardinality of $\Gamma$ is $\rk(\mathbf{v})^{2g}$. 
\qed

\section{Obstruction theory} \label{Sectionobs}

\subsection{Obstruction theories}\label{obstruction}
Let $\CF_{r}$ be the universal family on $S\times \rCoh(S,\mathbf{v})$, i.e.\ $\CF_{r}$ is the descend of $\CF\otimes p^*_{\Coh(S,\mathbf{v})}\lambda(\mathbf{u})^{-1}$ from $S\times \Coh(S,\mathbf{v})$ to $S\times \rCoh(S,\mathbf{v})$. Let 
\[ \tr\colon R\CH om_{\pi_{\rCoh(S,\mathbf{v})}}(\CF_{r}, \CF_r) \rightarrow H^*(S,\CO_S)\otimes \CO_{\rCoh(S,\mathbf{v})} \]
be the universal trace map. We define  
\[R\CH om_{\pi_{\rCoh(S,\mathbf{v})}}(\CF_{r}, \CF_r)_0:=\mathrm{Cone}(\tr)[-1].\]
 The complex above defines a perfect obstruction on $\rCoh(S,\mathbf{v})$,
\[(\BT^{\mathrm{vir}})^{\vee}:=(R\CH om_{\pi_{\rCoh(S,\mathbf{v})}}(\CF_{r}, \CF_{r})_{0}[1])^{\vee}\rightarrow \BL_{\rCoh(S,\mathbf{v})},\]
e.g.\ see \cite{STV} for more details. Note that the complex $(\BT^{\mathrm{vir}})^{\vee}$ is of amplitude [-1,1] due to the presence of non-discrete automorphisms of the unstable part of $\rCoh(S,\mathbf{v})$.
Let 
\begin{align*}
	\pi_{1}&\colon  \CC_{g,N} \rightarrow Q^{\epsilon}_{g,N}(M(\mathbf{v}), \beta), \\
	\mathbb{f}&\colon \CC_{g,N}\rightarrow \rCoh(S,\mathbf{v}),
\end{align*}
be the canonical projection from the universal curve and the universal map. On the other hand, let 
\begin{align*}\pi_2&\colon S\times \CC_{g,N}  \rightarrow M^{\epsilon}_{\mathbf{v}, \check{\beta}}(S \times C_{g,N}), \\
	\BF &\in \Cohc(S\times \CC_{g,N}),
\end{align*}
be the canonical projection from the universal threefold and the universal sheaf. We have two naturally defined obstruction-theory complexes:
\begin{itemize} 
	\item a map-theoretic obstruction complex $\pi_{1*}\mathbb{f}^*\BT^{\mathrm{vir}}$,
	\item a sheaf-theoretic obstruction complex $R\CH om_{\pi_2}(\BF,\BF)_{0}[1]$. 
\end{itemize}
In fact, these two complexes are  isomorphic. 

\begin{prop}\label{comparison}
	There is a natural identification of complexes, 
	\[R\CH om_{\pi_2}(\BF,\BF)_{0}[1] \cong \pi_{1*}\mathbb{f}^*\BT^{\mathrm{vir}}.\]
\end{prop}

\textit{Proof.}  Consider the following diagram
\[
\begin{tikzcd}
	S\times \CC_{g,N} \arrow[r,"\mathrm{id}\times \mathbb{f}"] \arrow[d] \arrow[dd, bend right=60,"\pi_2" near start]
	& S \times  \rCoh(S,\mathbf{v}) \arrow[d,"\pi_{ \rCoh(S,\mathbf{v})}"] \\
	\CC_{g,N}  \arrow[r,"\mathbb{f}"] \arrow[d, "\pi_{1}"] 
	&  \rCoh(S,\mathbf{v}) \\
	Q^{\epsilon}_{g,N}(M(\mathbf{v}),\beta)
\end{tikzcd}
\]
The trace map 
\[\mathrm{tr}\colon  R\CH om(\CF_{r},\CF_{r}) \rightarrow \CO_{\rCoh(S,\mathbf{v})}\]
 has a section given by the inclusion of the identity morphism $\CO_{\rCoh(S,\mathbf{v})}\rightarrow R\CH om(\CF_{r}, \CF_{r})$, therefore 
\[R\CH om(\CF_{r},\CF_{r})\cong R\CH om(\CF_{r},\CF_{r})_{0}\oplus \CO_{\rCoh(S,\mathbf{v})},\]
and by the moduli interpretation of $\rCoh(S,\mathbf{v})$, we get that
\[(\mathbb{f}\times \mathrm{id})^{*}\CF_{r}\cong\BF.\]
Hence by functoriality of the trace and the splitting above, we obtain that 
\[(\mathbb{f}\times \mathrm{id})^{*}R\CH om(\CF_{r},\CF_{r})_{0}\cong R\CH om(\BF,\BF)_{0}.\]
By the base change theorem,
\[R\CH om_{\pi_2}(\BF,\BF)_{0}\cong\pi_{1*}\mathbb{f}^{*}R\CH om_{\pi_{\rCoh(S,\mathbf{v})}}(\CF_{r},\CF_{r})_{0},\]
this proves the claim.
\qed
\begin{prop}\label{perf2}
	Assume $\BT^{\mathrm{vir}}$ is a locally free sheaf in degree 0 over the stable locus $M(\mathbf{v})$, then the complex $\pi_{1*}\mathbb{f}^*\BT^{\mathrm{vir}}$  is of perfect amplitude $[0,1]$. In particular, this holds for moduli spaces of sheaves on $K3$ surfaces and del Pezzo surfaces. 
\end{prop}
\textit{Proof.} A sheaf $F \in M^{\epsilon}_{\mathbf{v}, \check{\beta}}(S \times C_{g,N})(\BC)$ is perfect, since it is a family of sheaves on a smooth surface $S$ over $C$. A threefold $S\times C$ has at worst normal crossing singularities, its dualizing sheaf is therefore a line bundle.  Moreover, by Lemma \ref{simple}, the sheaf  $F$ is simple. Using Serre's duality and  simplicity of $F$,  we conclude that \[\Ext^i(F,F)_0=0, \quad  i\not\in  \{1,2,3\}.\]
Hence $R\CH om_{\pi_2}(\BF,\BF)_{0}[1]$ is  at most of amplitude $[0,2]$. By Proposition \ref{comparison}, the complex $\pi_{1*}\mathbb{f}^*\BT^{\mathrm{vir}}$  is therefore also at most of amplitude $[0,2]$. We need to show that it is zero in degree $2$.
    It is enough to check it over a point $[f\colon (C,\mathbf{p})\rightarrow \rCoh(S,\mathbf{v})]\in Q^{\epsilon}_{g,n}(M(\mathbf{v}), \beta)(\BC)$.  Consider the distinguished triangle 
\begin{equation}  \label{triangle}
	\tau_{\leq 0}f^*\BT^{\mathrm{vir}} \rightarrow f^*\BT^{\mathrm{vir}} \rightarrow \tau_{\geq 1} f^*\BT^{\mathrm{vir}}\rightarrow,
\end{equation}
where $\tau_{\dots}$ is the truncation of complex with respect to the standard $t$-structure. Taking the long exact sequence of cohomologies associated to (\ref{triangle}), we obtain
\[\dots \rightarrow H^2(C,\tau_{\leq 0}f^*\BT^{\mathrm{vir}}) \rightarrow H^2(C,f^*\BT^{\mathrm{vir}}) \rightarrow H^2(C,\tau_{\geq 1} f^*\BT^{\mathrm{vir}}) \rightarrow \dots\]
Let us now analyse the terms in the long exact sequence.  Firstly,  $f$ maps generically to the stable locus $M(\mathbf{v})$. By the assumption, $\BT^{\mathrm{vir}}_{|M(\mathbf{v})}$ is a locally free sheaf concentrated in degree 0. Hence the object $f^*\BT^{\mathrm{vir}}$ is  generically a locally free sheaf concentrated in degree 0. Moreover, $\BT^{\mathrm{vir}}$ is of perfect amplitude $[-1,1]$. These two facts imply that  $\tau_{\geq 1} f^*\BT^{\mathrm{vir}}$ is a $0$-dimensional sheaf in degree 1 supported on the base points of $f$, therefore 
\[ H^2(C,\tau_{\geq 1} f^*\BT^{\mathrm{vir}})=0.\]
Since $C$ is a curve, 
\[ H^2(C,\tau_{\leq 0}f^*\BT^{\mathrm{vir}}) =0.\] 
By the long exact sequence, we therefore obtain that $H^2(C,f^*\BT^{\mathrm{vir}})=0$. \qed

\begin{thm} \label{obsthe} There exists an obstruction-theory morphism
	\[\phi\colon (\pi_{1*}\mathbb{f}^*\BT^{\mathrm{vir}})^{\vee} \rightarrow  \BL_{Q^{\epsilon}_{g,N}(M(\mathbf{v}), \beta) / \FM_{g,N}}.\]
	Moreover, if $\epsilon=0^+$, under the assumption of Proposition \ref{perf2},   the corresponding virtual fundamental classes coincide with those of Donaldson--Thomas theory. 
\end{thm}

\textit{Proof.}  
By  \cite{TV, TV2, Lur} and \cite{STV}, there exists a derived algebraic stack $\BR \rCoh(S,\mathbf{v})$ whose classical truncation is $\rCoh(S,\mathbf{v})$, 
\[ t_0(\BR \rCoh(S,\mathbf{v}))\cong \rCoh(S,\mathbf{v}),\]
and 
\[j^* \BT_{\BR \rCoh(S,\mathbf{v})} \cong \BT^{\mathrm{vir}}\]
with respect to the natural inclusion  $j\colon \rCoh(S,\mathbf{v}) \hookrightarrow \BR \rCoh(S,\mathbf{v})$. 

A derived enhancement gives rise to an obstruction theory of the underlying classical stack, e.g.\ see \cite[Section 1]{STV} for more details. By Proposition \ref{openess}, the stack $Q^{\epsilon}_{g,N}(M(\mathbf{v}), \beta)$ is open in the mapping stack $\mathrm{Map}_{\FM_{g,N}}(\FC_{g,N}, \Coh_r(S,\mathbf{v})\times \FM_{g,N})$. Hence  the obstruction-theory morphism
\[\phi\colon (\pi_{1*}\mathbb{f}^*\BT^{\mathrm{vir}})^{\vee} \rightarrow  \BL_{Q^{\epsilon}_{g,N}(M(\mathbf{v}), \beta) / \FM_{g,N}}\]
can  be given by the derived mapping stack 
\[\BR\mathrm{Map}_{\FM_{g,N}}(\FC_{g,N}, \BR\Coh_r(S,\mathbf{v})\times \FM_{g,N}),\]
which is algebraic by \cite{TV2, Lur}. We refer to \cite[Theorem 5.1.1]{HLP} for a more recent treatment of derived mapping stacks; the stack $\BR \rCoh(S,\mathbf{v})$ satisfies the requirements of \cite[Theorem 5.1.1]{HLP}.  The obstruction-theory complex is perfect by Proposition \ref{perf2}.

By \cite{S}, a virtual fundamental class depends only on Chern characters of the corresponding obstruction-theory complex. The second part of the claim therefore follows from Proposition \ref{comparison}.\qed
\\

Let 
\[[Q^{\epsilon}_{g,N}(M(\mathbf{v}), \beta)]^{\mathrm{vir}} \in H_{*}(Q^{\epsilon}_{g,N}(M(\mathbf{v}), \beta))\]
be the associated virtual fundamental class. Invoking the identification presented in Proposition \ref{comparison}, the relative virtual dimension of the moduli space over $\FM_{g,N}$ can be computed via the virtual dimension of the relative moduli space of sheaves, 

\begin{align*}
	\mathrm{rel. \  vdim}&=\sum_i (-1)^{i+1}\mathrm{ext}^i(F,F)_0\\
	&= -\int_{S\times C} (\ch(F)\cdot \ch(F^{\vee})
	-1)\cdot\td_{S\times C}\\
	&=-\int_{S\times C} ((\ch(\mathbf{v}),\check{\beta})\cdot (\ch(\mathbf{v})^*,-\check{\beta}^*) 
	-1)\cdot\td_{S\times C} \\
	&=\check{\beta}_0\cdot (\mathrm{c}_{1}(\mathbf{v})\cdot\mathrm{c}_1(S))-\rk(\mathbf{v})\cdot (\check{\beta}_1\cdot \mathrm{c}_1(S))-(g-1)\dim(M(\mathbf{v})),
\end{align*}
where $\check{\beta}_0$ and $\check{\beta}_1$ are the components of $\check{\beta}\in \Lambda$ of cohomogical degrees 0 and 2, respectively. The class with the superscript ``*" stands for a class of the dual object. On the other hand, the Riemann--Roch formula gives us that 
\[ \mathrm{rel. \ vdim}=\deg (f^*\BT^{\mathrm{vir}})-(g-1)\dim(M(\mathbf{v})).\] 
Hence we conclude that
\begin{equation} \label{virt}
	\deg (f^*\BT^{\mathrm{vir}})=\check{\beta}_0\cdot (\mathrm{c}_{1}(\mathbf{v})\cdot\mathrm{c}_1(S))-\rk(\mathbf{v})\cdot (\check{\beta}_1\cdot \mathrm{c}_1(S)). 
	\end{equation}




\subsection{Invariants}\label{invariants}
Moduli spaces $Q^{\epsilon}_{g,N}(M(\mathbf{v}),\beta)$ have the usual structures to define the enumerative invariants in the spirit of Gromov--Witten theory:
\begin{itemize}
	\item evaluation maps at marked points
	\[\mathrm{ev}_{i}\colon  Q^{\epsilon}_{g,N}(M(\mathbf{v}),\beta) \rightarrow M(\mathbf{v}), \quad i=1, \dots, N,\]
	\item cotangent line bundles 
	\[L_{i}: =s^{*}_{i}(\omega_{\CC_{g,N}/Q^{\epsilon}_{g,N}(M(\mathbf{v}),\beta)}), \quad i=1, \dots, N,\]
\end{itemize}
where $s_{i}\colon  Q^{\epsilon}_{g,N}(M(\mathbf{v}),\beta) \rightarrow \CC_{g,N}$ are universal markings, and $\CC_{g,N}$ is the universal curve over $Q^{\epsilon}_{g,N}(M(\mathbf{v}),\beta)$.  We denote 
\[\psi_{i}:=\mathrm{c}_{1}(L_{i}), \quad i=1,\dots, N.\]

\begin{defn} \label{defninv} Under the assumption of Proposition \ref{perf2}, the \textit{descendent} $\epsilon$-\textit{invariants} are 
	\[\langle \gamma_{1}\psi^{k_1}, \dots, \gamma_{N}\psi^{k_N} \rangle^{\epsilon}_{g,\beta}:= \int_{[Q^{\epsilon}_{g,N}(M(\mathbf{v}),\beta)]^{\mathrm{vir}}}\prod^{i=N}_{i=1}\ev^{*}_{i}(\gamma_{i})\psi_{i}^{k_{i}},\]
	where $\gamma_{1}, \dots, \gamma_{N} \in H^{*}(M(\mathbf{v}))$ and $k_{1}, \dots, k_{N}$ are non-negative integers. 
\end{defn}

\begin{rmk} \label{absrel}We can also define another kind of invariants by the identification of quasimaps with the relative moduli space of sheaves from Proposition \ref{sheavesqms}.  Consider
	\[
	\begin{tikzcd}[row sep=small, column sep = small]
		&  S \times \CC_{g,N} \arrow{dl}[swap]{\pi_1} \arrow{dr}{\pi_{2}} \\
		S\times C_{g,N} && Q^{\epsilon}_{g,N}(M(\mathbf{v}), \beta)
	\end{tikzcd}
	\]
	where $\CC_{g,N}$ is the univeral curve over $Q^{\epsilon}_{g,N}(M(\mathbf{v}), \beta)$,  $C_{g,N}$ is the universal curve over $\Mbar_{g,N}$, and  $\pi_{1}$ is the morphism given by stabilisation of curves. For the unstable values of $g$ and $N$ we set the product $S\times \Mbar_{g,N+1}$ to be $S$. For a class $\tilde {\gamma} \in H^l(S\times C_{g,N})$, define the following operation on cohomology, 
	\[\ch_{k+2}(\tilde{\gamma})\colon  H_{*}(Q^{\epsilon}_{g,N}(M(\mathbf{v}), \beta)) \rightarrow H_{*-2k+2-l}(Q^{\epsilon}_{g,N}(M(\mathbf{v}), \beta)),\]
	\[\ch_{k+2}(\tilde{\gamma})(\xi)=\pi_{2*}\left(\ch_{k+2}(\BF)\cdot \pi_{1}^*(\tilde{\gamma})\cap \pi_2^*(\xi) \right).\]
Descendent Donaldson--Thomas  invariants are then defined by
	\begin{multline*} 
		\langle  \tilde{\tau}_{k_{1}}(\tilde{\gamma}_{1}), \dots, \tilde{\tau}_{k_{r}}(\tilde{\gamma}_{r}) \rangle^{\epsilon}_{g,n,\beta}\\
		=(-1)^{k_{1}+1}\ch_{k_{1}+2}(\tilde{\gamma}_{1}) \ \circ \ \dots \ \circ \ (-1)^{k_{r}+1} \ch_{k_{r}+2}(\tilde{\gamma}_{r}) \left([Q^{\epsilon}_{g,N}(M(\mathbf{v}), \beta)]^{\mathrm{vir}}\right).
	\end{multline*}
	We can also combine descendent quasimap invariants and descendent Donaldson--Thomas  invariants, 
	\[\langle   \gamma_{1}\psi^{k_1}, \dots, \gamma_{N}\psi^{k_N}  \mid \tilde{\tau}_{k_{1}}(\tilde{\gamma}_{1}), \dots, \tilde{\tau}_{k_{r}}(\tilde{\gamma}_{r})\rangle^{\epsilon}_{g,N,\beta},\]
	which are essentially a combination of relative and absolute Donaldson--Thomas invariants of the relative geometry \[ S \times C_{g,N} \rightarrow \Mbar_{g,N}\]
	for different $\epsilon$-stabilities. However, we will not be concerned with the absolute  Donaldson--Thomas invariants defined above in the present work.
\end{rmk}

The discussion in \cite[Section 6]{CFKM} also applies to $\epsilon$-invariants in our setting. In particular, $\epsilon$-invariants satisfy an analogue of the Splitting Axiom in Gromov--Witten theory, and there exists a projection to the moduli of stable nodal curves
\[p\colon Q^{\epsilon}_{g,N}(M(\mathbf{v}),\beta) \rightarrow \Mbar_{g,N}\]
by taking the stabilisation of the domain of a quasimap, so that the classes
\[p_{*}(\prod^{i=N}_{i=1} \mathrm{ev}^{*}_{i}(\gamma_{i})\psi_{i}^{k_{i}}) \in H^*(\Mbar_{g,N})\]
gives rise to a Cohomological field theory without a unit\footnote{In  the quasimap theory for $\epsilon\leq 1$, forgetting a marking is not permitted in general, because it might turn a rational bridge into a rational tail  which does not satisfy $\epsilon$-stability; hence a string equation might not hold true.} on $H^{*}(M(\mathbf{v}))$.

\section{Hilbert schemes} 
\subsection{Ideal sheaves} \label{relativeHilbert}We now restrict to the case of Hilbert schemes of points on a surface $S^{[n]}$. Hence, in this section, we assume that
\[ \ch(\mathbf{v})=(1,0,-n).\]  Punctorial Hilbert schemes are special, because there exists a canonical trivialisation of $\tau \colon \Coh(S,{\mathbf{v}}) \rightarrow \rCoh(S,{\mathbf{v}})$.  It is given by the determinant of the universal sheaf $\CF$ on $S \times \Coh(S, {\mathbf{v}})$,
\[\det(\CF) \in \Pic(S\times \Coh(S,{\mathbf{v}})).\]
Indeed, it is a line bundle of $\BC^*$-weight 1, because $\CF$ is of rank 1. Hence the family $\CF\otimes \mathrm{det}(\CF)^{-1}$ descends to $S\times \rCoh(S,{\mathbf{v}})$, giving the section 
\[s_{\det}\colon \rCoh(S,{\mathbf{v}}) \rightarrow \Coh(S,{\mathbf{v}}).\] 

Applying the same analysis as in Section \ref{sectionrel}, we obtain that there exists a canonical identification, 
\[ \phi \colon \det(\CF\otimes \mathrm{det}(\CF)^{-1})\xrightarrow{\sim} \CO_{ S\times\Coh(S,{\mathbf{v}})}.\]
Hence for any quasimap $f \colon C \rightarrow \rCoh(S,\mathbf{v})$, the determinant of the sheaf $F$ on $S\times C$ associated to the lift $s_{\det} \circ f$ can be canonically trivialised,
\[\phi \colon \det(F)\xrightarrow{\sim} \CO_{S\times C}.\] 

If $C$ is smooth, then since $F$ is a torsion free sheaf of rank 1, its double dual  $F^{\vee \vee}$ is locally free and therefore isomorphic to $\det(F)$. This is proved for families in \cite[Lemma 6.13]{Koll}. If $C$ is not smooth, then because $F$ is flat over nodes of $C$, we arrive at the same conclusion by passing to the normalisation of $C$. More generally, for families of nodal curves we can also apply  \cite[Lemma 6.13]{Koll} (which is stated only for the smooth case) to conclude that $F^{\vee \vee}$ is locally free as follows. Firstly,  since the question is local, we can split $\pi\colon  \CC \rightarrow B$ into the $\pi$-smooth locus $U_1$, and a Zariski neighbourhood $U_2$ of the $\pi$-singular locus. We require $U_2$ to be disjoint from the base locus of quasimaps, so that $F_{|S\times \CU_2}$  is flat  and torsion free over $\CU_2$ by the prestabiity of quasimaps. We then apply \cite[Lemma 6.13]{Koll} to the smooth maps $S\times U_1 \rightarrow B$ and $S\times U_2 \rightarrow U_2$. 

Overall, $\phi$ gives us a canonical identification, $F^{\vee\vee}\cong \det(F)\cong \CO_{S\times C}$, and since $F$ admits an embedding into its double dual, we obtain that $F$ canonically embeds into the trivial sheaf,
\[F \hookrightarrow F^{\vee \vee}\cong \CO_{S\times C}.\]
Finally, $F$ is  a family of 0-dimensional subschemes on $S$ at a general point of $C$, we conclude that $F$ is an ideal sheaf $I$ of a 1-dimensional subscheme,
\[F \cong I, \]
such that if the associated quasimap $f$ is of degree $\beta$, then
\[\ch(I)=((1,0,-n),\check{\beta})\in \Lambda\oplus \Lambda(-1).\]
We now will  translate $\epsilon$-stability of quasimaps into an explicit condition that involves only ideal sheaves and the associated subschemes. 

\subsection{Explicit $\epsilon$-stability for ideal sheaves} \label{explicit}   Let $\Gamma$ be a 1-dimensional subscheme on $S\times C$, and let $I$ be the associated ideal sheaf. The fiber $I_{p}$ of $I$ over a regular point $p\in C$ is stable\footnote{A sheaf of rank 1 is stable, if and only if it is torsion free.}, if and only if  $\Gamma$ is flat over $C$ at $p$, cf.\ Remark \ref{predeformable}. On the other hand,  a subscheme  $\Gamma$ is flat over $p \in C$, if and only if it does not have embedded points or non-dominant components over $p$. The latter requirement is equivalent to the condition that the structure sheaf $\CO_\Gamma$ of $\Gamma$ does not have $\CO_C$-torsion over $p$. 

Let
\[T^p \subseteq \CO_\Gamma \]
be the maximal $\CO_C$-torsion subsheaf of $\CO_\Gamma$ supported over a regular point $p\in C$. The sheaf $\CO_\Gamma/T^p$ is then flat over $p$.  By using the composition \[\CO_{S\times C} \rightarrow \CO_\Gamma \rightarrow \CO_\Gamma/T^p,\] we conclude that  $\CO_\Gamma/T^p$ is the structure sheaf of a subscheme, which we denote by $\Gamma^p$. By the discussion above, the ideal sheaf $I_{\Gamma^p}$ of $\Gamma^p$ is associated to the stabilisation of the quasimap $f$ at the point $p$ in the sense of Section \ref{stabilisation}. By construction, we have the following identity in the $K$-group of $S\times C$, 
\[[I]-[I_{\Gamma^p}]=-[T^p].\]
By Lemma \ref{key}, we obtain an expression for the length of $p$ in terms of $T^p$,
\[\ell(p)=m\cdot \deg(T^p)+\chi(T^p),\] where we used that $T^p$ is of rank 0 after the projection to $S$; the integer $m$ is chosen such that Proposition \ref{positive} holds, and the degree $\deg(T^p)$ is defined as $\mathrm{ch}(T^p)\cdot \mathrm{c}_1(\CO_S(1))$.  Hence the condition $\mathbf{(ii)}$ of Definition \ref{stabilityqm} translates into the following, 
\[ m\cdot \deg(T^p)+\chi(T^p)\leq 1/\epsilon.   \]

The condition $\mathbf{(i)}$ of Definition \ref{stabilityqm} can also be translated into the statement about the degree and the Euler characteristics of $\Gamma$ over rational bridges and rational tails.  Firstly, given a rational bridge $B\subseteq C$, then by using Lemma \ref{key}, we obtain that the following must be satisfied,
\[-m\cdot \mathrm{deg}(p_{S*}I_{|B})+\chi(\mathbf{v})\mathrm{rk}(p_{S*}I_{|B})- \chi(p_{S*}I_{|B})>0,\]
this simplifies to 
\[m\cdot \deg (\CO_{\Gamma|B})+\chi(\CO_{\Gamma|B})-n>0,\]
where $\deg (\CO_{\Gamma|B})$ is defined as $\ch(\CO_{\Gamma|B})_\mathrm{d}\cdot \mathrm{c}_1(\CO_S(1))$.  The term $n$ is there to cancel the fiber contribution\footnote{For example, if the quasimap is constant on a rational bridge $B$, then the associated $\CO_\Gamma$ satisfies the following identity $m\cdot \deg (\CO_{\Gamma|B})+\chi(\CO_{\Gamma|B})-n=0$, i.e.\ it is unstable with respect to $\epsilon$-stability as expected.}
 $\ch(\CO_{\Gamma|B})_\mathrm{f}$ to $\chi(\CO_{\Gamma|B})$. 
 
Similarly, given a rational tail  $R\subseteq C$, then the following must be satisfied,
\[m\cdot \deg (\CO_{\Gamma|R})+\chi(\CO_{\Gamma|R})-n>1/\epsilon.\]
Finally, by prestability of quasimaps,  $\CO_\Gamma$ has to be flat over nodes and marked points.   Summing up the discussion above, and using the main results of Section \ref{Sectionalg} and \ref{Sectionobs}, we obtain the following.  
\begin{cor} \label{Hilb} The space $Q^{\epsilon}_{g,N}(S^{[n]}, \beta)$ is a proper Deligne--Mumford stack with a perfect obstruction theory. There exists a natural identification of  moduli spaces, which respects perfect obstruction theories,  
	\[Q^{\epsilon}_{g,N}(S^{[n]}, \beta) \cong \Hilb_{n,\check{\beta}}^{\epsilon}(S \times C_{g,N}),\]
	where \[ \Hilb_{n,\check{\beta}}^{\epsilon}(S \times C_{g,N}) \colon (Sch/ \BC)^{\circ} \rightarrow Grpd \] 
	is the moduli space of triples $(C, \mathbf{p}, I)$,  satisfying the following properties for some fixed integer $m\gg0$:
	\begin{itemize}
		\item $I$ is the ideal sheaf of a 1-dimensional subscheme $\Gamma$ on $S\times C$, such that $\CO_{\Gamma}$ is flat over nodes and marked points of $C$,
		\item $\ch(\CO_\Gamma)=((0,0,n),-\check{\beta})\in \Lambda\oplus \Lambda(-1)$, 
		\item $m\cdot \deg(T^p)+\chi(T^p)\leq 1/\epsilon$ for all points $p \in C$,
		\item $m\cdot \deg (\CO_{\Gamma|R})+\chi(\CO_{\Gamma|R})-n>1/\epsilon$ for all rational tails $R\subseteq C$,
			\item $m\cdot \deg (\CO_{\Gamma|B})+\chi(\CO_{\Gamma|B})-n>0$ for all rational bridges $B\subseteq C$.
			\end{itemize}
\end{cor}
\textit{Proof.} Using the discussion in Section \ref{relativeHilbert} and in the beginning of Section \ref{explicit},  the identification of the moduli stacks follows from the same arguments as in Theorem \ref{mapssheaves}. The rest follows from the main results of Section \ref{Sectionalg} and \ref{Sectionobs}. 
\qed

\begin{rmk}\label{predeformable}
	Let $I$ be the ideal sheaf associated to a quasimap $f$, whose target is a Hilbert scheme of points. The fiber $I_{s}$ over a node $s \in C$ is stable, if and only if it is torsion free, which is equivalent to injectivity on the left of the exact sequence    
	\[I_{s} \rightarrow\CO_{S\times \{s\}}\rightarrow \CO_{\Gamma_{s}}\rightarrow 0,\]
	which in turn is equivalent to $\Tor_{S\times C}^{1}(\CO_{\Gamma}, \CO_{S\times \{s\}})=0$. In relative Donaldson--Thomas theory, the latter condition is referred to as \textit{admissibility}, it is one of the stability conditions of ideal sheaves on threefolds with normal crossing singularities. 
\end{rmk}

\begin{rmk}
	Fix a smooth curve $C$ with $g\geq1$ and $\beta \neq 0$.   By Corollary \ref{Hilb}, the moduli space of $0^+$-stable quasimaps from $C$ can be identified with a moduli space of 1-dimensional subschemes on $S\times C$,
	\[Q_{C}(S^{[n]},\beta)\cong \Hilb_{n,\check{\beta}}(S\times C).\]
 On the other hand, for a fixed smooth curve with one marking $(C,p)$, we obtain a moduli space of 1-dimensional subschemes relative to the divisor $S\times \{p\}\subset S\times C$, 
	\[Q_{(C,p)}(S^{[n]},\beta) \cong \Hilb_{n,\check{\beta}}(S\times C/S_{p}).\] 
	Moreover, pulling back a class from $S^{[n]}$ by the evaluation map associated to a marking on the left is equivalent to pulling back the class from a relative divisor on the right.
\end{rmk}

\subsection{Changing the t-structure}\label{perverse}
Consider the following \textit{torsion pair} in the abelian category $\mathrm{Coh}(S)$,
\begin{align*}
	\CT&=\{A\in \Cohc(S)\mid  \dim(A)=0\}, \\
	\CT^{\perp}&=\{A'\in \Cohc(S)\mid  \Hom(A,A')=0, \forall A\in \CT \}.
\end{align*}
Let \[\Cohc^{\sharp}(S)=\langle \CT^{\perp}, \CT[-1] \rangle\]
be the corresponding  perverse abelian heart, we refer to \cite{HRS} for the construction of abelian hearts associated to torsion pairs. A 2-term complex $A^\bullet=[A_0 \xrightarrow{d}A_1]$ belongs to $\Cohc^{\sharp}(S)$, if and only if
\[\ker(d)\in \CT^{\perp}, \quad \mathrm{coker}(d)\in \CT.\] 
A family of objects in $\Cohc^{\sharp}(S)$ over a base scheme $B$ is a perfect object
\[F\in \mathrm{D}_{\mathrm{perf}}(S\times B),\] such that for all closed points $b \in B$, we have 
\[F_{b}:=L \iota_{b}^{*}F\in \Cohc^{\sharp}(S), \]
where $\iota_b  \colon S\times \{b\} \hookrightarrow  S \times B$. We will refer to $F_{b}$ as fibers of $F$. 
Punctorial Hilbert schemes sit inside the rigidification of the corresponding moduli stack, 
\[S^{[n]} \subset \rpCoh(S, \mathbf{v}):=\pCoh(S, \mathbf{v}) \thickslash \BC^*,\]
which is constructed, for example, in \cite{Li}.
As in the case of the standard heart, we have a canonical section given by the determinant of the universal sheaf,
\begin{equation*}
	s_{\det}: \rpCoh(S, \mathbf{v})\rightarrow \pCoh(S, \mathbf{v}),
\end{equation*}
which exists  because the universal object is also of rank 1.  

Recall that a \textit{stable pair} on a scheme $X$ is a 1-dimensional sheaf $G$ with a section $g \in H^0(X,G)$, viewed as a 2-term complex,
\[I^\bullet =[\CO_X \xrightarrow{g} G],\]
such that $G$ is pure and $\mathrm{coker}(g)$ is 0-dimensional. We refer to \cite{PT} for the theory of stable pairs in the context of enumerative geometry. 

In what follows, we will also use the Abramovich--Polishchuk heart, denoted by  $\Cohc^{\sharp}(S\times C)$, on $S\times C$  associated to the perverse heart on $S$ for a smooth curve $C$. By definition \cite[Section 1.1]{AP}, an object $F\in \mathrm{D^b}(S\times C)$ is in $\Cohc^{\sharp}(S\times C)$, if and only if 
\[p_{S*}(F \otimes \CO_{C}(m)) \in \Cohc^{\sharp}(S), \quad \text{for all }  m\gg0.\] 
We now prove the following.

\begin{prop}\label{stablepair}
Let $f\colon C \rightarrow \rpCoh(S, \mathbf{v})$ be a prestable quasimap and 	let $F$ be the family on $S \times C$ associated to $s_{\det}\circ f$.  Then $F$ is a stable pair.  Conversely, given a stable pair $I^\bullet$ on $S\times C$ whose fibers are ideal sheaves over nodes of $C$, then it is a family of objects in $\Cohc^{\sharp}(S)$. 
\end{prop}

\textit{Proof.} The object $F$ is of rank 1, and $\det(F)\cong \CO_{S\times C}$ by the choice of the section $s_{\det}$. Hence by \cite[Lemma 3.11]{To}, in order to show that $F$ is a stable pair, we have to establish the following properties:
\begin{itemize}
	\item[$\mathbf{(i)}$] $\CH^{i}(F)=0$, for $i\neq0,1$, 
	\item[$\mathbf{(ii)}$] $\CH^{0}(F)$ is a rank 1 torsion free sheaf and $\CH^{1}(F)$ is 0-dimensional,
	\item[$\mathbf{(iii)}$] $\Hom(Q[-1], F)=0$ for any 0-dimensional sheaf $Q$,
\end{itemize}
where $\CH^{i}(F)$ is the cohomology of a complex with respect to the standard $t$-structure. Note that the proof of \cite[Lemma 3.11 (ii)]{To} easily extends to the case of a threefold with normal crossing singularities, as long as $F$ satisfies the first two properties above and its restriction  to the singular locus is an ideal sheaf (i.e.\ it is admissible over the  singular locus). Indeed, in this case, $\CH^{0}(F)$ is an ideal sheaf by the discussion in Section \ref{relativeHilbert}, $\CH^1(F)$ is supported away from the singularity,  and $F$ fits into the distinguished triangle, 
\[\CH^{0}(F) \rightarrow F \rightarrow \CH^1(F)[-1] \rightarrow,\]
then \cite[Lemma 3.11 (ii)]{To} follows from the same arguments as in the smooth case. 
\\

\noindent $\mathbf{(i)}$ Consider two distinguished triangles that the object $F$ fits in,
\[\tau_{\leq-1}F\rightarrow F\rightarrow \tau_{\geq0}F \rightarrow,\] 
\[ \tau_{\leq1}F\rightarrow F \rightarrow \tau_{\geq2}F \rightarrow,\]
where the truncation is taken with respect to the standard $t$-structure. Firstly, over a dense open subsect $U \subseteq C$, which contains all nodes, $F$ is a family of rank 1 sheaves. Hence
\[(\tau_{\leq -1}F)_{|U}=0, \quad (\tau_{\geq 2}F)_{|U}=0.\] 
We now deal with $p \in C\setminus U$.  By using the distinguished triangles above and the associated long exact sequences of sheaf cohomologies (cohomologies with respect to the standard $t$-structure), we obtain the long exact sequence
\[ \dots \rightarrow \CH^{-2}((\tau_{\geq0}F)_p) \rightarrow \CH^{-1}((\tau_{\leq -1 }F)_p) \rightarrow \CH^{-1}(F_p) \rightarrow \dots \] 
Since $C\setminus U$ contains only regular points, fibers $(\tau_{\geq0}F)_p$ are of amplitude $[-1,a]$ for some $a$. On the other hand, fibers $F_p$ are of amplitude $[0,1]$. We conclude that $\CH^{-1}((\tau_{\leq -1 }F)_p)=0$, the same applies to cohomologies of lower degrees, hence
\[ (\tau_{\leq-1}F)_p=0, \quad p\in C\setminus U.\]
By the same reasoning, we obtain that the sheaf cohomology $\CH^i((\tau_{\geq 2}F)_p)=0$ for $i\geq 2$. Since $\tau_{\geq 2}F$ is of amplitude $[2,b]$ for some $b$, by truncating $\tau_{\geq 2}F$ and repeating the argument, we obtain that  the underived fiber   $\iota_p^*\tau_{\leq2}(\tau_{\geq 2}F)$ is trivial.  Together with the vanishing of the higher sheaf cohomology, we conclude
\[ (\tau_{\geq 2}F)_p=0, \quad p\in C\setminus U.\]
Overall, we obtain that $\tau_{\leq-1}F=0$ and $\tau_{\geq 2}F=0$, since an object is trivial, if and only if all of its fibers are.  Hence $F$ must be of amplitude $[0,1]$. 
\\

\noindent $\mathbf{(ii)}$ Since $F_{p}$ is an ideal sheaf for a general $p \in C$ and $F_{p} \in \Cohc^{\sharp}(S)$ for all $p \in C$, the sheaf $\CH^{1}(F)$ must be 0-dimensional by the definition of $\Cohc^{\sharp}(S)$. Consider now the distinguished triangle,
\begin{equation} \label{trian}
 \CH^{0}(F) \rightarrow F \rightarrow \CH^{1}(F)[-1] \rightarrow,
 \end{equation}
and assume there is a non-zero torsion subsheaf (in the standard heart),
\[T(\CH^{0}(F)) \hookrightarrow \CH^{0}(F).\] Firstly, because a general fiber of $F$ is an ideal sheaf, $T(\CH^{0}(F))$ is supported on fibers over $C$. In particular, $\CH^{0}(F)$ is not flat over some point $p\in C$. Hence
$L\iota^*_p\CH^{0}(F)$ is strictly of amplitude $[-1,0]$ and $L\iota^*_p\CH^{1}(F)[-1]$ is at most of amplitude $[0,1]$,  because $C$ is 1-dimensional and $p$ can be assumed to be regular (since $F$ is an ideal sheaf over nodes).  This contradicts the assumption that $F$ is a family of objects in $\Cohc^{\sharp}(S)$, which are at most of amplitude $[0,1]$. Hence $\CH^{0}(F)$ must be torsion free. 
\\

\noindent$\mathbf{(iii)}$ Assume $\Hom(Q[-1], F)\neq0$ for a 0-dimensional sheaf $Q$. Note that around nodes $F$ is a family of rank 1 sheaves, hence the 0-dimensional sheaf $Q[-1]$ might be supported only away from nodes.  Since $Q[-1]$ is supported away from nodes, we can apply the results of  \cite{AP}, assuming that $C$ is smooth. In particular, the object $Q[-1]$ is in the heart $\Cohc^{\sharp}(S\times C)$ on $S\times C$ by the definition of the heart \cite[Section 1.1.1]{AP}. The family $F$ is also in $\Cohc^{\sharp}(S\times C)$ by \cite[Corollary 3.3.3]{AP}. By passing to the image of a morphism, it is enough to consider $Q[-1]$ as a torsion subobject of $F$.  Let $T(F) \subset F$ be the maximal torsion subobject, then $F':=F/T(F)$ is a torsion free object, hence it is flat by \cite[Corollary 3.1.3]{AP}. Restricting to a fiber over some $p \in C$ we obtain an exact sequence in $\Cohc^{\sharp}(S\times C)$,
\[0\rightarrow T(F)_{p} \rightarrow F_{p} \rightarrow F'_{p} \rightarrow 0,\]
because $F'$ is flat. We get that $T(F)_{p}\in \Cohc^{\sharp}(S)$. Since both $F$ and $F'$ are families of perverse objects, which are generically isomorphic (generically they are families of rank 1 sheaves), we have $\ch(F_{p})=\ch(F'_{p})$ for all $p \in C$. By  exactness of the sequence above, we conclude that $\ch(T(F)_{p})=0$ for all $p \in C$. It is not difficult to check that an object $A^\bullet \in  \Cohc^{\sharp}(S)$ is trivial, if and only if $\ch(A^\bullet)=0$ (this follows from the same result in the standard heart $\Cohc(S)$). Hence $T(F)_{p}=0$ for all $p \in C$, which implies that $T(F)=0$. We conclude that $Q=0$. 
\\

Conversely, given now a stable pair $I^{\bullet}$. To show that it is a family of objects in $\Cohc^{\sharp}(S)$, we must show that $I^\bullet_p \in \Cohc^{\sharp}(S)$ for all $p \in C$. We may assume that $C$ is smooth by passing to the normalisation, since it does not change the fibers, and fibers of $I^{\bullet}$  over nodes are ideal sheaves.  By definition,  $I^{\bullet}$ sits in a distinguished triangle
\[ \ker(g) \rightarrow I^{\bullet} \rightarrow \mathrm{coker}(g)[-1] \rightarrow.\] 
Applying $p_{S*}(-\otimes \CO_{C}(m))$ for $m \gg 0$ to the triangle, we obtain that $p_{S*}(\ker(g) \otimes \CO_{C}(m))$ is a torsion free sheaf and $p_{S*}(\mathrm{coker}(g) \otimes \CO_{C}(m))$ is a $0$-dimensional sheaf. Hence $p_{S*}(I^{\bullet} \otimes \CO_{C}(m))\in \Cohc^{\sharp}(S)$ for $m \gg 0$.  By the definition \cite[Section 1.1]{AP}, we obtain that $I^{\bullet}$ is contained in the heart $\Cohc^{\sharp}(S\times C)$ on $S\times C$. Moreover, by the assumption that $G$ is pure, $I^{\bullet}$ is torsion free in the sense of \cite[Defintion 3.1.1]{AP}. Hence by \cite[Corollary 3.1.3]{AP}, we conclude that it is a family of objects in $\Cohc^{\sharp}(S)$. 
\qed 
\\

The determinant-line-bundle construction in this setting also provides us with the map 
\[\lambda\colon  K_{0}(S) \rightarrow \Pic(\pCoh(S, \mathbf{v})).\]
 The line bundles $\CL_{0}$ and $\CL_ {1}$ satisfy the same properties as in the case of the standard heart. 
\begin{lemma}\label{positive2}
	Let $f\colon  C \rightarrow\pCoh(S, \mathbf{v})$ be a prestable quasimap of degree $\beta$. Then there exists $m_{0} \in \BN$ which depends only on $\beta$, $\mathbf{v}$ and $\CO_S(1)$\footnote{While Hilbert schemes of points are independent of stability, we still need $\CO_S(1)$ for the definition of $\CL_1$. }, such that for all $m \geq m_{0}$, the quasimap is non-constant, if and only if
	\[ \CL_{0}\otimes \CL_{1}^{m} \cdot_{f} C>0.\]
	This also holds true for all subcurves $C'\subseteq C$ and the induced maps for the same choice of $m$.
\end{lemma}  
\textit{Proof.} The proof is similar to the one of Proposition \ref{positive}. We use Proposition \ref{stablepair} to identify objects associated to perverse quasimaps with stable pairs. We then consider the distinguished triangle 
\[\ker(g) \rightarrow I^\bullet \rightarrow \mathrm{coker}(g)[-1] \rightarrow, \] 
where $\ker(g)$ is an ideal sheaf, and $\mathrm{coker}(g)$ is a $0$-dimensional sheaf. Let $f'$ be the quasimap associated to $\ker(g)$. Note that it is a quasimap to the pair $(S^{[n]}, \rpCoh(S, \mathbf{v}))$, because the subscheme associated to $\ker(g)$ does not have embedded points.  By Lemma \ref{key}, which can be proved in this setting using the same arguments as in the standard heart, we obtain that 
\begin{align*} 
\CL_{0}\otimes\CL_{1}^{m} \cdot_{f} C&=\CL_{0}\otimes\CL_{1}^{m} \cdot_{f'} C-\chi(\mathrm{coker}(g)[-1]) \\
&=\CL_{0}\otimes\CL_{1}^{m} \cdot_{f'} C+\chi(\mathrm{coker}(g)).
\end{align*}
The quantity $\chi(\mathrm{coker}(g))$ is positive, if and only if $\mathrm{coker}(g)$ is non-trivial. We then apply the arguments of Proposition \ref{positive} to the quasimap $f'$. 
\qed
\\

Fixing a positive line bundle $\CL_{\beta}$ from  Lemma \ref{positive2} once and for ever for all $\beta \in \Eff(S^{[n]}, \rpCoh(S, \mathbf{v}))$, we can define the notion of $\epsilon$-stability  as in Definition \ref{stabilityqm}. 
\begin{defn} Let
\[Q^{\epsilon}_{g,N}(S^{[n]}, \beta)^{\sharp}\colon  (Sch/ \BC)^{\circ} \rightarrow Grpd\]
be the moduli space of $\epsilon$-stable quasimaps of genus $g$ and degree $\beta$ with $N$ marked points associated to a pair $(S^{[n]}, \rpCoh(S, \mathbf{v}))$. 
\end{defn}

 For properness of these moduli spaces, we will need a version of  Proposition \ref{extension0} for the perverse heart. 

\begin{lemma} \label{extension2}
	Let $\CC \rightarrow \Delta$ be a family of nodal projective curves, and  let $\{p_1, \dots, p_m\} \subset \CC$ be a set of finitely many closed points in the regular locus of the central fiber. Then any quasimap 
	$\tilde{f}\colon  \tilde{\CC}=\CC\setminus \{p_1, \dots, p_m\} \rightarrow \rpCoh(S, \mathbf{v})$ extends to $f\colon \CC \rightarrow \rpCoh(S, \mathbf{v})$, which is unique up to unique isomorphism.
\end{lemma}
\textit{Proof.} Employing the same proof as the one of Lemma \ref{extension0} is problematic in this case, as we do not know how to extend objects in the derived category (unlike sheaves), so we follow a different strategy. 

Restricting $\tilde{f}$ to the generic fiber $\CC^{\circ}$ of $\CC$ over $\Delta$, we obtain a stable pair $I^{\bullet,\circ}$ on $S\times\CC^{\circ}$, which is flat over $\CC^{\circ}$, hence it is flat over $\Delta$ by Lemma \ref{flatness} (the argument of Lemma \ref{flatness} easily extends to complexes). By the properness of the relative moduli space of stable pairs $\mathrm{P}(S\times\CC/\Delta)$, it  can be completed to a pair $I^\bullet$ on $S\times \CC$. We refer to \cite[Section 5]{LiWu} for the properness of stable pairs on threefolds with normal crossing singularities. By Lemma \ref{torsion}, the central fiber may be non-flat only over nodes. Hence by Lemma \ref{torsion} and Lemma \ref{flatness}, it defines a quasimap with indeterminacies  (i.e.\ a non-everywhere defined quasimap, or a rational quasimap),
 \[f\colon  \CC \dashrightarrow \rpCoh(S, \mathbf{v}),\]  such that indeterminacies are possibly only at the nodes of the central fiber. We will show that it does not have any indeterminacies. 
 
 Firstly,  $f$ also defines a rational map, 
\[f_{\mathrm{rat}}\colon  \CC \dashrightarrow S^{[n]},\]
and so does $\tilde{f}$, 
\[\tilde{f}_{\mathrm{rat}}\colon  \CC \dashrightarrow S^{[n]}.\] 
 The corresponding graphs in $\Hilb(S^{[n]}\times \CC)$ agree generically, therefore by the separatedness of Hilbert schemes they are equal, i.e.\ 
 \[f_{\mathrm{rat}}=\tilde{f}_{\mathrm{rat}}.\] 
If $p_{i}$ is not a limit of base points of $\tilde{f}$, then there is a neighbourhood $U \subset \CC$ around $p_{i}$, where \[\tilde{f}_{|U\setminus \{p_{i}\}}=\tilde{f}_{\mathrm{rat}|U\setminus \{p_{i}\}}=f_{\mathrm{rat}|U\setminus \{p_{i}\}}=f_{|U\setminus \{p_{i}\}},\]
we then define $\tilde{f}_{|U}=f_{|U}$ ($f$ is defined at $p_{i}$, because $p_{i}$ is in the regular locus). Since quasimaps to $\rpCoh(S, \mathbf{v})$ are given by stable pairs, which do not have automorphisms other than the $\BC^*$-scaling (which is removed, because $\rpCoh(S,\mathbf{v})$ is rigidified), we can glue maps in a unique way, thereby extending $\tilde{f}$ to $p_{i}$. 

If $p_{i}$ is a limit base point of $\tilde{f}$, let $B_{i} \subset \CC$ be the  section corresponding to these base points, then there is some neighbourhood $U$ around $B_{i}$, such that
\[\tilde{f}_{|U\setminus B_{i}}=\tilde{f}_{\mathrm{rat}|U\setminus B_{i}}=f_{\mathrm{rat}|U\setminus B_{i}}=f_{|U\setminus B_{i}},\]
but since $\tilde{f}_{|\CC^{\circ}}=f_{|\CC^{\circ}}$, we conclude that $\tilde{f}_{U\setminus \{p_{i}\}}=f_{U\setminus \{p_{i}\}}$. We then proceed as before. Let 
\[f'\colon \CC \rightarrow \rpCoh(S, \mathbf{v})\] be the resulting extension and ${I^{\bullet}}'$ be the associated family. By the separatedness of relative moduli spaces of stable pairs (see e.g.\  \cite[Section 5]{LiWu}), we get that ${I^{\bullet}}'  \cong I^\bullet$, and that the extension is unique. In particular, $f=f'$. 
\qed
\\

Using the discussion above, we can adopt the proofs of the main results for the standard heart to the case of the perverse heart. For the expression of $\epsilon$-stability in terms of stable pairs $I^\bullet=[\CO_{S\times C} \xrightarrow{g} G]$ like in Corollary \ref{Hilb}, we define $[T^p]$ as a $K$-theoretic sum of the $\CO_C$-torsion of $G$ over $p$, denoted by $T'^p$, and the restriction of $\mathrm{coker}(g)$ over $p$, 
\[[T^p]:=[T'^p]+[\mathrm{coker}(g)_{p}],\]
while instead of $\CO_\Gamma$, we take the sheaf $G$, cf. the proof of  Lemma \ref{positive2}. 
\begin{thm}\label{PT}  The space $Q^{\epsilon}_{g}(S^{[n]}, \beta)^{\sharp}$ is a proper Deligne--Mumford stack with a perfect obstruction theory. There exists a natural identification of moduli spaces, which respects perfect obstruction theories, 
	\[Q^{\epsilon}_{g,N}(S^{[n]}, \beta)^{\sharp} \cong \mathrm{P}_{n,\check{\beta}}^{\epsilon}(S \times C_{g,N}),\]
	where
	\[\mathrm{P}_{n,\check{\beta}}^{\epsilon}(S \times C_{g,N}) \colon  (Sch/ \BC)^{\circ} \rightarrow Grpd\]
is the moduli space of triples $(C,\mathbf{p},I^\bullet)$,  satisfying the following properties for some fixed integer $m\gg0$\text{:}
	\begin{itemize}
		\item $I^\bullet=[\CO_{S\times C} \xrightarrow{g} G]$ is a stable pair, such that $G$ is flat over nodes and marked points of $C$, and $\mathrm{coker}(g)$ is supported away from nodes and marked points, 
		\item $\ch(G)=((0,0,n),-\check{\beta})\in \Lambda\oplus \Lambda(-1)$,
		\item $m\cdot \deg([T^{p}])+\chi([T^{p}])\leq 1/\epsilon$ for all points $p \in C$,
		\item $m\cdot \deg (G_{|R})+\chi(G_{|R})-n>1/\epsilon$ for all rational tails $R\subseteq C$,
		\item $m\cdot \deg (G_{|B})+\chi(G_{|B})-n>0$ for all rational bridges $B\subseteq C$.
		
	\end{itemize}
\end{thm} 
\textit{Proof.}  Firstly, as in the case of the standard heart, Corollary \ref{algebraicity}, we deduce that $Q^{\epsilon}_{g,N}(S^{[n]}, \beta)^{\sharp}$ is a quasi-separated algebraic stack locally of finite presentation. Indeed, this follows from the representability of mapping stacks \cite{HR}, since the stack $\pCoh(S,\mathbf{v})$ is algebraic, locally of finite presentation and has affine stabilizers (given by invertible elements in $\Hom(A^\bullet, A^\bullet)$) by \cite{Li}. 
On the other hand, $\mathrm{P}_{n,\check{\beta}}^{\epsilon}(S \times C_{g,N})$ is a quasi-separated algebraic stack locally of finite presentation as a relative moduli space of coherent systems by arguments from \cite[Section 4]{LiWu} (see also \cite{LeP}). 

Consider now  another moduli space $\mathrm{P}_{n,\check{\beta}}^{\epsilon}(S \times C_{g,N})'$, whose $B$-valued points are given by quadruples $(\CC,\mathbf{p}, F, \phi)$, where $F$ is an object in $\mathrm{D_{perf}}(S\times \CC)$ with a trivialisation of the  determinant $\phi\colon \det(F)\xrightarrow{\sim} \CO_{S\times \CC}$,  $(\CC,\mathbf{p})$ is a family of nodal curves over $B$, and over geometric points $b\in B$ objects $F_{|S\times \CC_b}$ define $\epsilon$-stable quasimaps to $\rpCoh(S, \mathbf{v})$. This is indeed a well-defined 2-functor by the construction of $\pCoh(S,\mathbf{v})$ in \cite{Li}. By arguments of Theorem \ref{mapssheaves}, we have an identification of $2$-functors, 
\begin{equation} \label{identpt}
 Q^{\epsilon}_{g,N}(S^{[n]}, \beta)^{\sharp} \cong \mathrm{P}_{n,\check{\beta}}^{\epsilon}(S \times C_{g,N})'.
 \end{equation}
  We obtain that $\mathrm{P}_{n,\check{\beta}}^{\epsilon}(S \times C_{g,N})'$ is  algebraic.  
 Since the universal stable pair on $\mathrm{P}_{n,\check{\beta}}^{\epsilon}(S \times C_{g,N})$ defines a family of objects in $\mathrm{P}_{n,\check{\beta}}^{\epsilon}(S \times C_{g,N})'$, there is a natural map between algebraic stacks, 
\begin{equation} \label{mapstacks}
	\mathrm{P}_{n,\check{\beta}}^{\epsilon}(S \times C_{g,N}) \rightarrow \mathrm{P}_{n,\check{\beta}}^{\epsilon}(S \times C_{g,N})', 
	\end{equation}
 which is an equivalence on geometric points\footnote{In fact, it is already enough to conclude that the virtual intersection theories of these moduli spaces are equivalent by the results of \cite{MaPush}, because the relative obstruction theories of both spaces are given by $R \Hom(I^\bullet, I^\bullet )_0$, which follows from the first-order version of \cite[Theorem 2.7]{PT}.} by Proposition \ref{stablepair} (note that the proposition extends to curves over a field of characteristic 0). 
 
  One can readily verify  that  \cite[Theorem 2.7]{PT} extends to the case of relative threefolds $S\times \CC \rightarrow B$ with nodal fibers, as long as stable pairs  are admissible over nodes. More precisely, in the terminology of \cite{PT}, the result holds\footnote{The proof of \cite[Theorem 2.7]{PT} involves projective dimension and Serre's duality arguments, which extend to a moving threefold; if stable pairs are admissible over nodes, the same arguments apply to a singular threefold $S\times C$; see also Section \ref{relativeHilbert} for how to prove \cite[Lemma 6.13]{Koll} in the singular case. } for a nilpotent thickening given by an extension $S\times \CC_0 \subset S\times \CC$ flat over $B_0 \subset B$.   It implies that (\ref{mapstacks}) is formally \'etale, and since both stacks are locally of finite presentation, it is \'etale. Moreover, it is an equivalence on geometric points, inducing an isomorphism on stabilisers. Hence by \cite[Section 0DUD]{stacks-project}, (\ref{mapstacks}) is an isomorphism. Using (\ref{identpt}), we therefore obtain the desired identification,
  \[  Q^{\epsilon}_{g,N}(S^{[n]}, \beta)^{\sharp} \cong \mathrm{P}_{n,\check{\beta}}^{\epsilon}(S \times C_{g,N}).\]

We now deal with properness of $Q^{\epsilon}_{g,N}(S^{[n]}, \beta)^{\sharp}$. It  is quasi-compact, and therefore of finite presentation, by Lemma \ref{positive2} and quasi-compactness of moduli spaces of stable pairs on threefolds with normal crossing singularities (see \cite[Section 5.25]{LiWu}).  The stability of stable pairs and $\epsilon$-stability imply that the stack is Deligne--Mumford. For the properness we use Lemma \ref{extension2} and the proof presented in \cite[Proposition 4.3.1]{CFKM}.  The existence and the compatibility of the perfect obstruction theories follows from the same arguments as in Proposition \ref{comparison} and \ref{perf2}, we refer to \cite{TV}  and  \cite{STV}  for the derived enhancement of a moduli space of objects with a fixed determinant in $\mathrm{D^b}(S)$.  
\qed

\begin{defn}
Theorem \ref{PT} allows us to define  \textit{perverse descendent} $\epsilon$-\textit{invariants},
\[ \langle \gamma_{1}\psi^{k_1}, \dots, \gamma_{N}\psi^{k_N} \rangle^{\sharp,\epsilon}_{g,\beta}:=\int_{[Q^{\epsilon}_{g,N}(S^{[n]},\beta)^{\sharp}]^{\mathrm{vir}}}\prod^{i=N}_{i=1}\ev^{*}_{i}(\gamma_{i})\psi_{i}^{k_{i}}, \] 
where all terms have the same meaning as  in Definition \ref{defninv}.
\end{defn}

\begin{rmk}
	As in the case of the standard heart, for a fixed smooth curve $C$ with $g\geq1$ and $\beta \neq 0$,  by Theorem \ref{PT} we have an identification with the moduli space of stable pairs, 
	\[Q_{C}(S^{[n]},\beta)^{\sharp}\cong \mathrm{P}_{n,\check{\beta}}(S\times C). \]
	 The same applies to a curve with a marking, 
	\[Q_{(C,p)}(S^{[n]},\beta)^{\sharp} \cong \mathrm{P}_{n,\check{\beta}}(S\times C/S_{p}).\] 
	
\end{rmk}

\subsection{Affine plane} \label{affineplane}
A punctorial Hilbert scheme of the affine plane  $\BC^2$ admits two equivalent descriptions, one is a Nakajima quiver variety, which is a GIT construction, 
\[(\BC^2)^{[n]}=[\mu^{-1}(0)/\mathrm{GL}_n]^{s}\subset [\mu^{-1}(0)/\mathrm{GL}_n],\]
we refer to \cite{G} for the notation. Another description is provided by a moduli space of framed sheaves on $\BP^{2}$. Both descriptions sit in some bigger stacks.  To match the unstable loci, one has to consider a non-standard heart of $\mathrm{D^{b}}(\BP^{2})$, namely, the perverse heart $\Cohc^{\sharp}(\BP^2)$. Let 
\[\mathfrak{Coh}_{\mathrm{fr}}^{\sharp}(\BP^2, \mathbf{v}),\] 
be the stack of objects in $\Cohc^{\sharp}(\BP^2)$ with framing at the line at infinity $\ell_\infty \subset \p^2$. Framing kills $\BC^*$-automorphisms, hence rigidification is not necessary.  Then by \cite[Theorem 5.7]{BFG}, we have a canonical isomorphism 
\[[\mu^{-1}(0)/\mathrm{GL}_n]\cong \mathfrak{Coh}_{\mathrm{fr}}^{\sharp}(\BP^2, \mathbf{v}),\]
which identifies stable loci on both sides. Hence moduli spaces of GIT quasimaps \cite{CFKM} and our moduli spaces of quasimaps to $(\BC^2)^{[n]}$ are isomorphic, 
\[Q^{0^+}_{g,N}((\BC^2)^{[n]}, \beta)^{\text{GIT}} \cong Q^{0^+}_{g,N}((\BC^2)^{[n]}, \beta)^{\sharp}.\]
Moreover, since $[\mu^{-1}(0)/\mathrm{GL}_{n}]$ is a local complete intersection, an easy virtual-dimension calculation shows that the obstruction-theory map, 
\[R\CH om_{\pi} (\CF, \CF)_{0}[1]^{\vee}\rightarrow \BL_{\mathfrak{Coh}_{\mathrm{fr}}^{\sharp}(\BP^2, \mathbf{v})},\]
is an isomorphism. Hence obstruction theories of both quasimap theories also match. We conclude that GIT and moduli-of-sheaves quasimap theories are equivalent in this case. 



\section{Wall-crossing} \label{wallcrossingg}
\subsection{$I$-function} We fix a parametrized projective line $\p^1$ with a $\BC^*$-action, 
\[t[x:y]=[tx:y], \quad t\in \BC^{*},\]
such that $0:=[0:1]$ and $\infty:=[1:0]$. By convention we set 
\[z:=e_{\BC^*}(\BC_{\mathrm{std}}),\] 
where $\BC_{\mathrm{std}}$ is the weight 1 representation of $\BC^*$. 
We now define a \textit{Vertex} space, see Remark \ref{VI} for the origin of the notation. Let
\[V(M(\mathbf{v}),\beta) \colon (Sch/ \BC)^{\circ} \rightarrow Set\]
be the space of quasimaps  $f \colon \p^1 \rightarrow \rCoh(S,\mathbf{v})$ subject to two conditions:
\begin{itemize}
	\item $f$ is of degree $\beta$,
	\item  $f(\infty) \in M(\mathbf{v})\subset \rCoh(S,\mathbf{v})$. 
	\end{itemize}
Note that we do not impose stability on quasimaps, i.e.\ they are allowed to have infinitely many automorphisms, and we do not identify maps by automorphisms of $\p^1$.  

Using the argument from Theorem \ref{mapssheaves}, the space $V(M(\mathbf{v}),\beta)$ can be viewed a moduli space of sheavs $F$ on $S\times \p^1$ subject to the following conditions:
\begin{itemize}
	\item $\ch(F)=(\ch(\mathbf{v}),\check{\beta})$, 
	\item $F$ is torsion free,
	\item the fiber $F_\infty$  is stable, 
	\item $\det(p_{\p^1*}(p_{S}^{*}\mathbf{u}\otimes F))\cong \CO_{\p^1}$.
\end{itemize}
 By construction, there is a natural evaluation map, 
\begin{align*}\ev\colon V(M(\mathbf{v}),\beta) \rightarrow M(\mathbf{v}), \quad f \mapsto f(\infty).
	\end{align*}
Moreover, by acting on the source of quasimaps, we obtain a $\BC^*$-action on $V(M(\mathbf{v}),\beta)$. Under the assumption of Proposition \ref{perf2}, the space $V(M(\mathbf{v}),\beta)$ has a perfect obstruction theory. However, it is not proper, only its $\BC^*$-fixed locus 
\[V(M(\mathbf{v}),\beta)^{\BC^*}\] is proper. Indeed, this follows from properness of the space of all prestable quasimaps from $\p^1$ and the fact that $V(M(\mathbf{v}),\beta)^{\BC^*}$ is just a connected  component of its $\BC^*$-fixed locus. We use the virtual localisation of \cite{GP} to construct its virtual fundamental class,
\[[V(M(\mathbf{v}),\beta)]^{\mathrm{vir}}:=\frac{[V(M(\mathbf{v}),\beta)^{\BC^*}]^{\mathrm{vir}}}{e_{\BC^*}(\CN^{\mathrm{vir}})} \in H_*(V(M(\mathbf{v}),\beta)^{\BC^*})[z^{\pm}],\]
where $\CN^{\mathrm{vir}}$ is the virtual normal complex of $V(M(\mathbf{v}),\beta)^{\BC^*}$ inside $V(M(\mathbf{v}),\beta)$, and $z$ is the equivariant parameter. We are now ready to define Givental's $I$-function, introduced in the context of GIT quasimaps in \cite{CFKM}.

\begin{defn}We define 
	\begin{align*}
		I_\beta(z)&:= \ev_* [V(M(\mathbf{v}),\beta)]^{\mathrm{vir}} \in H^*(M(\mathbf{v}))[z^\pm], \\
		\mu_{\beta}(z)&:=[zI_\beta(z)]_{z^{\geq 0}}  \in H^*(M(\mathbf{v}))[z].
	\end{align*}
$I$-function and its truncation are then defined as the following generating series,
\begin{align*}
	I(q,z)&:=1+ \sum_{\beta \neq0} I_\beta(z)q^\beta \in H^*(M(\mathbf{v}))[z^\pm][\![q^\beta]\!], \\
	\mu(q,z)&:=\sum_{\beta \neq0} \mu_{\beta}(z)q^\beta \in H^*(M(\mathbf{v}))[z][\![q^\beta]\!].
\end{align*}
\end{defn}

\begin{rmk}\label{VI} In the context of Donaldson--Thomas theory, $I$-functions are also known as \textit{Vertex} functions. The latter terminology originates from \cite{MNOP1,MNOP}. To pay tribute to both theories, we denote $I$-functions by ``$I$", while the spaces that are used to define them by ``$V$".
\end{rmk}
These functions will govern the quasimap wall-crossing formulas. 
\subsection{Master space and wall-crossing} \label{masterHilb}
For the material discussed in this section we refer the reader to \cite{YZ}. Here, we just glide over the machinery developed in \cite{YZ}, adjusting some minor details to our needs. 

The space $\BR_{>0} \cup \{0^+,\infty\}$ of $\epsilon$-stabilities is divided into chambers, in which the moduli space $Q^{\epsilon}_{g,N}(M(\mathbf{v}),\beta)$ stays the same, and as $\epsilon$ crosses a wall between chambers, the moduli space changes discontinuously. Given a class $\beta \in \Eff(M(\mathbf{v}), \rCoh(S,\mathbf{v}))$,  for a class $\beta'\in \Eff(M(\mathbf{v}), \rCoh(S,\mathbf{v}))$ which is a summand of $\beta$, we define 
\[\deg(\beta'):=\beta'(\CL_\beta),\]
and, in this section, we define $\epsilon$-stability of quasimaps of degree $\beta'$ with respect to the line bundle $\CL_\beta$. 

 Let $\epsilon_{0}=1/d_{0} \in \BR_{>0}$ be a wall for $\beta \in \Eff(M(\mathbf{v}), \rCoh(S,\mathbf{v}))$ and $\epsilon_-$, $\epsilon_+$ be some values that are close to $\epsilon_{0}$ from left and right of the wall, respectively. Assuming \[2g-2+N+\epsilon_0\deg(\beta)>0,\]
  i.e.\ $\epsilon$-stability is defined\footnote{If $2g-2+N+\epsilon_0\deg(\beta)\leq 0$, then the moduli space of $\epsilon^-$-stable quasimaps is empty.} on both sides of the wall,  we define
\[MQ^{\epsilon_{0}}_{g,N}(M(\mathbf{v}),\beta) \rightarrow M\widetilde{\FM}_{g,N,d}\]
to be the master space with the projection to a moduli space of curves with calibrated tails, constructed in \cite{YZ}. The construction applies to our case verbatim. The space $ M\widetilde{\FM}_{g,N,d}$ is a $\p^1$-bundle  over $\widetilde{\FM}_{g,N,d}$. The latter space  is obtained by a series of blow-ups of a moduli space of semistable curves weighted by the degree $\FM^{\mathrm{ss}}_{g,N,d}$, such that the total degree is $d=\deg(\beta)$.  As in the GIT case, we have the following result.
\begin{thm} The space $MQ^{\epsilon_{0}}_{g,N}(M(\mathbf{v}),\beta)$ is a proper Deligne--Mumford stack. 
\end{thm}
\textit{Proof.} With Lemma \ref{extension0} the proof is exactly the same as in GIT case, we therefore refer to \cite[Section 5]{YZ}. 
\qed 
\\

The master space also carries a perfect obstruction theory, which is obtained in the same way as the one for $Q^{\epsilon}_{g,N}(M(\mathbf{v}),\beta)$. Let 
\begin{align*}\mathbb{f} &\colon  \CC_{g,N} \rightarrow \rCoh(S,\mathbf{v}), \\
\pi&\colon  \CC_{g,N} \rightarrow MQ^{\epsilon_0}_{g,N}(M(\mathbf{v}), \beta)
\end{align*}
be the universal quasimap and the canonical projection.
\begin{thm} There exists an obstruction-theory morphism
	\[\phi\colon (\pi_*\mathbb{f}^*\BT^{\mathrm{vir}})^{\vee} \rightarrow  \BL_{ MQ^{\epsilon_0}_{g,N}(M(\mathbf{v}), \beta) / M\widetilde{\FM}_{g,N,d}}.\]
	If $\BT^{\mathrm{vir}}_{|M(\mathbf{v})}$ is a locally free sheaf in degree 0, then the complex $\pi_{*}\mathbb{f}^*\BT^{\mathrm{vir}}$  is of perfect amplitude $[0,1]$.
\end{thm}
\textit{Proof.} The first claim follows from the same arguments as in Theorem \ref{obsthe}. The second claim follows from the arguments of Proposition \ref{perf2}.
\qed 
\\

Using the master space, we can establish the wall-crossing formula. 
\begin{thm} \label{wallcrossingHilb}Given a wall $\epsilon_0\in \BR_{>0}$. Assuming $2g-2+N+\epsilon_0\deg(\beta)>0$, we have
	\begin{multline*}
		\langle \lambda_1 \psi^{k_1}_1, \dots,\lambda_N \psi^{k_N}_N \rangle^{\epsilon_-}_{g,\beta} = \langle \lambda_1 \psi^{k_1}_1, \dots,\lambda_N \psi^{k_N}_N \rangle^{\epsilon_+}_{g,\beta} \\
		+\sum_{\underline{\beta}} \langle \lambda_1 \psi^{k_1}_1, \dots,\lambda_N \psi^{k_N}_N, \mu_{\beta_1}(-\psi_{N+1}), \dots, \mu_{\beta_k}(-\psi_{N+k}) \rangle^{\epsilon_+}_{g,\beta_0}/k!,
	\end{multline*}
	where $\underline{\beta}$ runs through all the $(k+1)$-tuples of effective  classes 
	\[\underline{\beta}=(\beta_0, \beta_{1}, \dots, \beta_{k}),\]
	such that $\beta=\sum^{i=k}_{i=0}\beta_i$, and $\deg(\beta_{i})=d_{0}$ for all $i\in \{1, \dots, k\}$. The $\epsilon_{+}$-stability for the class $\beta_0$ is given by $\CL_{\beta}$. The same holds for perverse quasimap invariants 	$\langle \lambda_1 \psi^{k_1}_1, \dots,\lambda_N \psi^{k_N}_N \rangle^{\sharp, \epsilon}_{g,\beta}$. 
\end{thm}
\textit{Sketch of Proof.} Here we will only sketch the proof,  for all the details we refer to \cite[Section 6]{YZ}, as the proof in our case is exactly the same as the one for GIT quasimaps. 

The master space  $MQ^{\epsilon_{0}}_{g,N}(M(\mathbf{v}),\beta)$ carries a natural $\BC^*$-action, such that the $\BC^*$-fixed locus is a union of the following three types of spaces (up to finite coverings):
\begin{itemize}
	\item $Q^{\epsilon^-}_{g,N}(M(\mathbf{v}),\beta)$,
	\item $\widetilde{Q}^{\epsilon^+}_{g,N}(M(\mathbf{v}),\beta)$, the  base change of $Q^{\epsilon^+}_{g,N}(M(\mathbf{v}),\beta)$  from $\FM_{g,N,d}$ to $\widetilde{\FM}_{g,N,d}$,
	\item $Y\times_{M(\mathbf{v})^{k}} \prod_{i=1}^k  V(M(\mathbf{v}),\beta_i)^{\BC^*}$,  where $Y$ is a finite-group gerbe  over the space  $\widetilde{Q}^{\epsilon^+}_{g,N+k}(M(\mathbf{v}),\beta_0)$.
\end{itemize}
Applying the virtual localisation formula and  taking the equivariant residue, we obtain certain relations between the classes associated to the spaces above. Projecting everything to a point, we get the wall-crossing formula. 

All efforts are aimed at the careful construction of the master space and the analysis of moving and fixed parts of the obstruction theories at fixed loci. For more details we refer the reader to \cite[Section 6]{YZ}. \qed 
\begin{rmk} In the GIT set-up, there are naturally defined maps $[W/G]\rightarrow [\BC^{n+1}/\BC^*]$, which induce 
	$Q^{\epsilon}_{g,N}(W/G,\beta)\rightarrow Q^{\epsilon}_{g,N}(\p^n,d)$. This allows to give a more refined class-valued wall-crossing by pushforwarding the classes on $MQ^{\epsilon_0}_{g,N}(W/G, \beta)$ to $Q^{\epsilon^{-}}_{g,N}(\p^n,d)$ instead of a point. In our case, this seems to be less natural. Even though $\rCoh(S,\mathbf{v})$ is Zariski locally a GIT stack, we do not have these naturally defined maps, because it is unclear, if a line bundle $\CL_{\beta}$ is ample on any of the GIT loci through which the universal quasimap factors. Moreover, these loci change as we change $\beta$. 
	
	It is also possible to pushforward the classes to $\Mbar_{g,N}$ instead of $Q^{\epsilon}_{g,N}(\p^n,d)$. The problem with this approach is that the projection \[Q^{\epsilon}_{g,N+k}(M(\mathbf{v}),\beta) \rightarrow \Mbar_{g,N}\] involves stabilisation of a curve, which implies that $\psi$-classes do not pullback to $\psi$-classes. Consequently, the wall-crossing formula becomes inefficient to state.  
\end{rmk}
Crossing all walls from $\epsilon=0^+$ to $\epsilon=\infty$ and applying Theorem \ref{wallcrossingHilb} inductively, we obtain a statement that relates two extremal values of $\epsilon$. 

\begin{cor} \label{wallcrossingHilb2} Assuming $(g, N)\neq (0,1)$, we have
	\begin{multline*}
		\langle \lambda_1 \psi^{k_1}_1, \dots,\lambda_N \psi^{k_N}_N \rangle^{0^+}_{g,\beta} = \langle \lambda_1 \psi^{k_1}_1, \dots,\lambda_N \psi^{k_N}_N \rangle^{\infty}_{g,\beta} \\
		+\sum_{\underline{\beta}} \langle \lambda_1 \psi^{k_1}_1, \dots,\lambda_N \psi^{k_N}_N, \mu_{\beta_1}(-\psi_{N+1}), \dots, \mu_{\beta_k}(-\psi_{N+k}) \rangle^{\infty}_{g,\beta_0}/k!,
	\end{multline*}
	where $\underline{\beta}$ runs through all the $(k+1)$-tuples of effective quasimap classes 
	\[\underline{\beta}=(\beta_0, \beta_{1}, \dots, \beta_{k}),\]
	such that $\beta=\sum^{i=k}_{i=0}\beta_i$, and $\beta_{i}\neq 0$ for all $i\in \{1, \dots, k\}$. The same holds for perverse quasimap invariants 	$\langle \lambda_1 \psi^{k_1}_1, \dots,\lambda_N \psi^{k_N}_N \rangle^{\sharp, \epsilon}_{g,\beta}$. 
\end{cor}
 
The result above can be restated as a change-of-variables relation of the following generating series. Let $\epsilon\in \{0^+, \infty\}$, we define 
\[F^{\epsilon}_{g}(\mathbf{t}(z)):=\sum^{\infty}_{N=0}\sum_{\beta\geq0}\frac{q^{\beta}}{N!}\langle \mathbf{t}(\psi), \dots, \mathbf{t}(\psi) \rangle^{\epsilon}_{g,N,\beta},\]
where $\mathbf{t}(z) \in H^{*}(M)[\![z]\!]$ is a generic element, and the unstable terms are set to be zero. 
\begin{cor} \label{changeofvariableHilb} For all $g\geq 1$, we have
	\[F^{0^+}_{g}(\mathbf{t}(z))=F^{\infty}_{g}(\mathbf{t}(z)+\mu(-z)).\]
	For $g=0$, the same equation holds true modulo constant and linear terms in $\mathbf{t}$. 	
\end{cor}

\subsection{The genus-zero case}If 
\[2g-2+N+\epsilon_0\deg(\beta)\leq 0,\]
then the moduli space $Q^{\epsilon^-}_{0,1}(M(\mathbf{v}),\beta)$ is empty, and the wall-crossing formula  takes a different form.

\begin{thm} \label{unstable}
For $\epsilon \in (\frac{1}{\deg(\beta)},\frac{1}{\deg(\beta)-1})$, we have 
\[\ev_*\left( \frac{[Q^{\epsilon}_{0,1}(M(\mathbf{v}),\beta)]^{\mathrm{vir}}}{z(z-\psi_1)}\right)=[I_\beta(z)]_{z^{\leq-2}},\]
where $[\dots]_{z^{\leq-2}}$ means that we take a truncation up to $z^{-2}$.
\end{thm}

\textit{Proof.} See \cite[Lemma 7.2.1]{YZ}. \qed 
\\

To express the wall-crossing formula above in terms of a change of variables, we do the following. Let $\{B^{i}\}$ be a basis of $H^{*}(M(\mathbf{v}))$ and $\{B_{i}\}$ be its dual basis with respect to the intersection pairing. Let 
\begin{multline*}
J^{0^+}(\mathbf{t}(z),q,z):=\frac{\mathbf{t}(-z)}{z}+I(q,z)\\
+\sum_{\beta\geq0,N\geq0}\frac{q^\beta}{N!}\sum_{i}B_{i}\langle \frac{B^{i}}{z(z-\psi)}, \mathbf{t}(\psi), \dots, \mathbf{t}(\psi)\rangle^{0^+}_{0,1+N,\beta},
\end{multline*}
where unstable terms are set to be zero, and let
\begin{multline*} J^{\infty}(\mathbf{t}(z),q,z):=\frac{\mathbf{t}(-z)}{z}+1\\
+\sum_{\beta\geq0,N\geq0}\frac{q^\beta}{N!}\sum_{i}B_{i}\langle \frac{B^{i}}{z(z-\psi)}, \mathbf{t}(\psi), \dots, \mathbf{t}(\psi)\rangle^{\infty}_{0,1+N,\beta},
\end{multline*}
then Theorem \ref{wallcrossingHilb} and Theorem \ref{unstable} imply the following. 
\begin{cor} \label{unstable2}
We have 
\[J^{\infty}(\mathbf{t}(z)+\mu(-z))=J^{0^+}(\mathbf{t}(z)).\]
\end{cor}
\textit{Proof.} We again refer to \cite[Section 7.4]{YZ}. \qed

\section{Semi-positive targets} \label{semi-positiv}

\subsection{$I$-function.} The complexity of $I$-functions largely depends on virtual dimensions of corresponding Vertex spaces. There is a big class of moduli spaces for which $I$-functions can be given a more explicit form. 
\begin{defn}
A pair $(M(\mathbf{v}),\rCoh(S,\mathbf{v}))$ is \textit{semi-positive}, if for all prestable quasimaps $f \colon C \rightarrow \rCoh(S,\mathbf{v})$, the following holds 
\[\deg (f^*\BT^{\mathrm{vir}})\geq 0.\]
\end{defn}
An example of a semi-positive target would be a Hilbert scheme of points $S^{[n]}$ on a del Pezzo surface, e.g.\ $\p^2$. Indeed, by the ampleness of the anti-canonical bundle, Corollary \ref{Hilb} and (\ref{virt}), we obtain that
\begin{align*}\deg(f^*\BT^{\mathrm{vir}})&=\check{\beta}_0\cdot (\mathrm{c}_{1}(\mathbf{v})\cdot\mathrm{c}_1(S))-\rk(\mathbf{v})\cdot (\check{\beta}_1\cdot \mathrm{c}_1(S))\\
	&=-\check{\beta}_1\cdot \mathrm{c}_1(S)\geq 0,
	\end{align*}
where by Corollary \ref{Hilb}, the class   $-\check{\beta}_1$ is an effective curve class on $S$. In particular,  $\check{\beta}_1\cdot \mathrm{c}_1(S)=0$, if and only if $\check{\beta}_1=0$. 
\\

Consider now  the expansion 
\[\mu(q,z)=I_{1}(q)+(I_{0}(q)-1)z+I_{-1}(q)z^2+I_{-2}(q)z^3+\dots\]
We will show that all terms $I_{k}$ with $k\leq-1$  vanish for a semi-positive target. 

\begin{lemma}\label{semipos} For a semi-positive target, we have
	\[I_{k}=0, \quad \text{if }k\leq -1.\]
	\end{lemma}
\textit{Proof.} The virtual dimension of $V(M(\mathbf{v}),\beta)$ is equal to 
\[\deg(f^*\BT^{\mathrm{vir}})+\dim(M(\mathbf{v})).\]
 Hence cohomological  degrees of  classes
\[\ev_* [V(M(\mathbf{v}),\beta)]^{\mathrm{vir}}\in H^*(M(\mathbf{v}))[z^\pm]\]
are equal to 
\[-2\deg(f^*\BT^{\mathrm{vir}}),\]
such that $z$ is of cohomological degree 2. Cohomological degrees of $I_k$ are therefore of the following form
\[2k-2\deg(f^*\BT^{\mathrm{vir}}).\]
The assumption of semi-positivity and the fact that $H^*(M(\mathbf{v}))$ is non-negatively graded imply the claim. 
\qed 
\\

The terms $I_0(q)$ and $I_1(q)$ have the following useful expression, which we, however, will not use in this work. 

\begin{prop} \label{semi-positive}
For a semi-positive target, we have
\begin{itemize}
	\item[$\mathbf{(i)}$]
	\begin{equation*} I_{0}(q)^{-1}=1+\sum_{\beta\neq0}q^{\beta}\langle [\mathrm{pt}],\mathbbm{1},\mathbbm{1} \rangle_{0,\beta}^{0^+},
	\end{equation*}
	
	\item[$\mathbf{(ii)}$] \[I_{1}(q)=f_{0}(q)\mathbbm{1}+\sum_{j}f_{j}(q)D_{j},\]
\end{itemize}
where $\{D_{j}\}$ is a basis of $H^2(M(\mathbf{v}))$, and 
\[\frac{f_{0}(q)}{I_{0}(q)}=\sum_{\beta\neq0}q^{\beta}\langle[\mathrm{pt}], \mathbbm{1} \rangle_{0,\beta}^{0^+} \quad \frac{f_{j}(q)}{I_{0}(q)}=\sum_{\beta\neq0} \sum_j q^{\beta}\langle D^{j}, \mathbbm{1} \rangle_{0,\beta}^{0^+}.\]
\end{prop}

\textit{Proof.} The result follows from the same arguments as  \cite[Corollary 5.5.4]{CFK}. However, it can be obtained more directly via the wall-crossing formula, Corollary \ref{wallcrossingHilb2}, by plugging in the  insertions above, and  using the string and divisor equations on the Gromov--Witten side. 
\qed

\subsection{Computations}
We will now explicitly  compute $I$-functions associated to a perverse pair $(S^{[n]}, \rpCoh(S, \mathbf{v}))$ for a del Pezzo surface $S$.  We start with some notational preparations. To derive expressions for $I$-functions, we firstly need to remind the reader how the cohomology of $S^{[n]}$ looks like. By the virtual dimension constraints, we only need the degree 2 cohomology, or, dually, the degree 2 homology.  The homology  $H_2(S^{[n]},\BZ)$ admits a Nakajima--Grojnowski basis,
\begin{align*}
	H_2(S,\BZ) \oplus \BZ &\xrightarrow{\sim} H_2(S^{[n]},\BZ), \\
	(\gamma, k) &\mapsto \gamma_n+kA,
\end{align*}
where the classes above are defined in terms of Nakajima--Grojnowski operators as follows, 
\[\gamma_n=\mathfrak{q}_{-1}(\gamma)\mathfrak{q}_{-1}([\mathrm{pt}])^{n-1}1_S, \quad  A=\mathfrak{q}_{-2}([\mathrm{pt}])\mathfrak{q}_{-1}([\mathrm{pt}])^{n-2}1_S,\]
we refer to \cite[Section 1]{Ob18} for the notation and the definition of Nakajima--Grojnowski operators in the similar context. 
In more geometric terms, if a class $\gamma$ is represented by a curve $\Gamma \subset S$, then the class $\gamma_n$ is represented by the curve $\Gamma_n \subset S^{[n]}$ which is given by letting one point move along $\Gamma$ and keeping $n-1$ other distinct points fixed. The class $A$ is given by the locus of  length 2 non-reduced structures on a fixed point $p\in S$, keeping other $n-2$ reduced points fixed.  We will apply the same notation for the dual cohomology classes in $H^2(S^{[n]},\BZ)$.

Let us now discuss degrees of quasimaps. By Corollary \ref{Hilb}, the sheaves associated to quasimaps are ideals of 1-dimensional subschemes, hence a class $-\check{\beta}$ is of the following form, 

\begin{equation} \label{cohidentification}
-\check{\beta}=(0,\gamma, k)\in H^{0}(S)\oplus H^{1,1}(S)\oplus H^{4}(S),
\end{equation}
where the negative sign amounts to considering Chern characters of subschemes rather than their ideals. For simplicity we will denote a class $-\check{\beta}$ just by $(\gamma, k)$. The decomposition above and the one given by Nakajima--Grojnowski basis are related. Given a map $f \colon C \rightarrow S^{[n]}$, then with respect to the identification above, we have
\[f_*[C]=\gamma_n+mA\]
By \cite[Lemma 2]{Ob18}, the associated Chern character  is of the following form,
\begin{equation} \label{km}
-\check{\beta}=(\gamma, k)=\left( \gamma, m-\mathrm{c}_1(S)\cdot \gamma/2 \right).
\end{equation}
Using this formula, we can write degrees of quasimaps in terms of curve classes in $H_2(S^{[n]},\BZ)$. In fact, this is the most convenient basis to express $I$- functions, because $\check\beta$ is not necessarily an integral class. We therefore will use the following identification
\begin{align*}
\BQ[\![q^\beta]\!]\xrightarrow{\sim}& \BQ[\![q^\gamma,  y]\!], \quad q^\beta=q^\gamma \cdot y^{m}, \\
\beta \mapsto& (\gamma,m).
\end{align*}
\begin{rmk} \label{rmk}Another expression of degrees of quasimaps in the case of $S^{[n]}$ involves Euler characteristics of the associated subscheme $\Gamma$, which is related to $m$ as follows, 
	\[ \chi(\Gamma)=m+n(1-g(C)),\]
	this is useful to keep in mind for Lemma \ref{IfunctionFano}.
	\end{rmk}
With this notation, we will compute $I$-functions in Proposition \ref{IfunctionFano}, a result which was kindly communicated to the author by Georg Oberdieck. However, before doing that, we need the following lemma about exceptional classes on del Pezzo surfaces. 
\begin{lemma} \label{curves} Let $S$ be a del Pezzo surface and $\gamma \in H^2(S,\BZ)$ be an effective curve class, such that $\gamma \cdot \mathrm{c}_1(S)=1$. If $\mathrm{c}_1(S)^2\neq1$, then $\gamma$ is an exceptional class, i.e.\  $\gamma=[\p^1]$ and $\gamma^{2}=-1$. If $\mathrm{c}_1(S)^2=1$, then $\gamma$ is an exceptional class or $\gamma=\mathrm{c}_1(S)$.
\end{lemma}	
	
\textit{Proof.} Assume $\mathrm{c}_1(S)^2\neq1$, then by \cite[Proposition 3.4]{Koll} the divisor $\mathrm{c}_1(S)$ is base-point free. Hence a curve $\Gamma$ in the class  $\gamma$ is irreducible and is mapped with degree 1 onto a line via the anti-canonical map. This implies that $\Gamma\cong \p^1$. By the adjunction formula, we get that $\gamma^2=-1$. 

Assume $\mathrm{c}_1(S)^2=1$, then $\mathrm{c}_1(S)$ is not base-point free. In this case, $S$ is a blow-up of $\p^2$ in 8 generic points. We use a direct lattice computation to show that  $\gamma=\mathrm{c}_1(S)$ or $\gamma$ is exceptional. The cohomology of $S$ can be described as follows, 
\begin{align*}
&H^2(S,\BZ)=\BZ\langle H, E_1, \dots, E_8\rangle ,\\
&H^2=1, \quad E_i^2=-1, \quad  H\cdot E_i=0, \quad E_i\cdot E_j=0\  \text{ for }i \ \neq j,
\end{align*}
such that the anti-canonical class has the following expression, 
\[\mathrm{c}_1(S)=3H-\sum^{8}_{i=1} E_i.\]
Consider now a class 
\[\gamma= a_0H + \sum^{8}_{i=1}a_iE_i.\]
By the assumption, 
\begin{equation} \label{boundv}
	\gamma \cdot \mathrm{c}_1(S)= 3a_0+\sum^{8}_{i=1}a_i=1.
	\end{equation}
By the Cauchy--Schwarz inequality applied to the standard Euclidean pairing of  vectors $(a_1,\dots,a_8)$ and $(1,\dots ,1)$, and (\ref{boundv}),  we obtain the following bound,  
\begin{equation} \label{bound2}
	\frac{(1-3a_0)^2}{8}\leq \sum^{8}_{i=1}a^2_i,
	\end{equation}
such that the equality holds, if and only if $(a_1,\dots,a_8)=(a,\dots,a)$.
Consider now the self-intersection product, 
\[\gamma^2=a_0^2 - \sum^{8}_{i=1}a_i^2,\]
using the bound (\ref{bound2}), we obtain 
\[\gamma^2\leq -\frac{a^2_0-6a_0+1}{8}.\]
Using the fact that $\gamma^2$ is an integer, we conclude  $\gamma^2 \leq 1$, such that the equality holds, if and only if $(a_1,\dots,a_8)=(-1,\dots ,-1)$ and $a_0=3$. In other words, $\gamma^2 =1$, if and only if  $\gamma=\mathrm{c}_1(S)$. Assume now $\gamma^2 \leq 0$, then by the adjunction formula, we obtain that a curve $\Gamma$ in the class $\gamma$ satisfies 
\[ 2g_{a}(\Gamma)-2=\gamma^2-\gamma \cdot \mathrm{c}_1(S)\leq -1.\]
Moreover,  $\Gamma$ is irreducible by $\gamma \cdot \mathrm{c}_1(S)=1$ and  the ampleness of $\mathrm{c}_1(S)$. Hence $g_{a}(\Gamma)=0$ and $\gamma^2=-1$. Since $\Gamma$ is Gorenstein, we also get that  $\Gamma\cong \p^1$. \qed 
\begin{prop}[Georg Oberdieck] \label{IfunctionFano} Assume $S$ is a del Pezzo surface, such that $\mathrm{c}_1(S)^2\geq 2$. Then for a perverse pair $(S^{[n]}, \rpCoh(S, \mathbf{v}))$, we have
\begin{align*}
	I_0(q,y) & = 1, \\
	I_1(q,y) & = \log( 1 + y) \mathrm{c}_1(S)_n +\frac{1}{1+y}\left( \sum_{\gamma\cdot \mathrm{c}_1(S)=1}q^\gamma\right) \mathbbm{1}.
\end{align*}
If $\mathrm{c}_1(S)^2=1$, the same holds with the exception that the second term in $I_1(q,y)$ acquires a contribution from $\gamma=\mathrm{c}_1(S)$. 
\end{prop}
\textit{Proof.}
The pair $(S^{[n]}, \rpCoh(S, \mathbf{v}))$ is semi-positive. Hence by Lemma \ref{semipos} and its proof, and (\ref{virt}), we only need to consider classes $(\gamma, m)$, such that 
\[\gamma \cdot \mathrm{c}_1(S)\leq 1.  \]
Since $\mathrm{c}_1(S)$ is ample, we have $\gamma \cdot \mathrm{c}_1(S)=0$, if and only if $\gamma=0$. We therefore start our analysis with classes of the form 
\[-\check \beta =(0, m).\] In this case, 
\[I_0(q) \in H^0(S^{[n]})[\![q^\beta]\!], \quad  I_1(q)\in H^2(S^{[n]})[\![q^\beta]\!].\]
\noindent To evaluate $I_0(q)$, it is enough to compute its pairing  $\langle I_0(q), [P] \rangle$ with a point class $[P]\in H^{4n}(S^{[n]})$.
Moreover, using \cite[Proposition 5.10]{BeFa}\footnote{We use the fact  the evaluation map can be upgraded to a map of derived schemes, such that $S^{[n]}$ is a classical smooth scheme, while the fiber is taken in the category of derived schemes; this gives us the compatibility of obstruction theories.} and a pushforward-pullback theorem for the fiber product with respect to the embedding of a point $P \hookrightarrow S^{[n]}$, we obtain that this pairing can be computed as the virtual degree of the fiber, 
\[\langle I_{(0,m),0}, [P] \rangle=\deg [\ev^{-1}(P)]^{\mathrm{vir}},\]
where the equivariant localisation is used to define the quantity on the right.
   Let $P \in S^{[n]}$ be represented  by a collection of $n$ distinct points on $S$. We identify $P$ with the corresponding subset of points $P \subset S$. Then the fiber above parametrizes stable pairs supported on $P \times \p^1 \subset S\times \p^1$, which are admissible  over $\infty \in \p^1$. In other words, 
\[\ev^{-1}(P)\cong \coprod_{\underline{m}} \prod^{n}_{i=1}\Sym^{m_i}\BC \subset \coprod_{\underline{m}}\prod^{n}_{i=1} \Sym^{m_i}\p^1,\]
where the union is taken over $n$-tuples of positive integers $(m_1,\dots, m_n) $, such that $m=\sum^{i=n}_{i=1}m_i$.  
The deformation theory of such pairs inside a threefold $S \times \p^1$  is completely local, hence it can be assumed that $S=\BC^2$ (see \cite[Section 4.2]{PT} for the deformation theory). In this case, the local model is provided by $\BC^2\times  \Sym^{m_i}\BC$, and taking the fiber removes deformations coming from the first factor.  Since $\BC^2$ carries a symplectic form, the obstruction theory acquires a surjective cosection which remains intact after taking the fiber. Hence the contribution of $\prod^{n}_{i=1}\Sym^{m_i}\BC$ is zero, if at least one $m_i$ is non-zero. We conclude that all classes of the form $(0, m)$  with non-zero $m$ do not contribute to $I_0(q)$. Hence $\langle I_0(q), [P] \rangle = 1$, which implies that $I_0(q)=1$.

We now consider the term $I_1(q)$. The class $A$ can be represented by a smooh curve (this is needed to apply  \cite[Proposition 5.10]{BeFa}), namely, $\p^1$.  By the same arguments as above, we therefore obtain that 
\[\langle I_1(q), A \rangle = 0,\]
since this is again a completely local question with the difference that a stable pair  acquires a varying non-reduced structure of multiplicity 2 along one of the components of $P\times \p^1$.
 
 Now let us  evaluate $I_1(q)$ at the classes in $H_2(S,\BZ) \subseteq H_2(S^{[n]},\BZ)$.  Since $S$ is a del Pezzo surface,  $H_2(S,\BZ)$ is generated by exceptional curve classes and a hyperplane class.  Hence we may assume that  $\gamma_n$ is represented by a smooth curve of the form $\Gamma_n$.  The fiber product with respect to $\Gamma_n \hookrightarrow S^{[n]}$ has the following form
\[\ev^{-1}(\Gamma_n)=\coprod_{\underline{m}} \left( \left(\Gamma \times \Sym^{m_1}\BC\right) \times \prod^{n}_{i=2}\Sym^{m_i}\BC \right).\]
As before, the only non-vanishing contribution comes from the case when $m_{1}=m$. In this case,  each factor $\Sym^{0}\BC$
 contributes 1 to the product formula, while the contribution of $\Gamma \times \Sym^{m} \BC$ is independent of $n$. We conclude that
\begin{equation} \label{reduction}
\langle I^{S^{[n]}}_1(q), \gamma_n \rangle = \langle I^{ S}_1(q), \gamma \rangle.
\end{equation}
Hence we may assume $n=1$.
In this case, the fixed locus 
\[V(S,(0,m))^{\BC^*}=\mathrm{P}_{1,(0,m)}(S\times \BC)^{\BC^*}\] is isomorphic to $S$, parametrizing stable pairs
of the form  
\[I^\bullet =[\CO_{\p^1_\mathrm{pt} } \to \CO_{\p^1_\mathrm{pt} }(D)],\]
where $\p^1_\mathrm{pt} = \p^1 \times \{\mathrm{pt} \}$ for a point $ \mathrm{pt} \in S$, and $D = m \cdot [0]$.
The obstruction theory was computed in \cite[Section 4.2]{PT},
\begin{align*}
\mathrm{Def}_{I^\bullet}&=H^0(\p^1, \CO_D(D))\oplus T_{S,\mathrm{pt} } \\ \mathrm{Obs}_{I^\bullet}&=H^0(\p^1, \CO_D(D)\otimes \omega_{S\times \p^1})^\vee=H^0(\p^1,\CO_D(D)\otimes \omega_{ \p^1})^\vee \otimes \omega_{S,\mathrm{pt} }^\vee.
\end{align*}
Recall that the $\BC^{\ast}$-action on $\p^1$ is given by $t \cdot [x:y] = [t x :y]$.
The coordinate function $X= x/y$ acquires the dual scaling,  $t \cdot X = t^{-1}X$, hence has weight $-z$. Let us analyse the $\BC^*$-equivariant structure of the obstruction theory. Firstly, 
\[ H^{0}(\p^1,\CO_D(D)) = (X^{-m}) \otimes \BC[X]/X^{m}
= \BC X^{-m} \oplus \BC X^{-m+1} \oplus \ldots \oplus \BC X^{-1}, \]
which therefore has weights $z, 2z, \ldots, mz$ as a $\BC^*$-representation.
Moreover, the fiber $\omega_{\p^1, 0} $ over $0\in \p^1$ has weight $-z$, so we get
that $H^0(\p^1,\CO_D(D) \otimes \omega_{\p^1})$
has weights $0, z, \ldots, (m-1) z$, therefore its dual has weights
$(-m+1)z, \ldots, -z , 0$. 	We therefore obtain the following
\begin{align*}
\ev_{*} [V(S,(0,m))]^{\mathrm{vir}}
& = p_{S*} \left( \frac{ e_{\BC^*}( \mathrm{Obs}_{I^\bullet}^{\text{mov}} ) }{ e_{\BC^*}( \mathrm{Def}_{I^\bullet}^{\text{mov}} ) }  \right)\\
& =\frac{ (-z + \mathrm{c}_1(S)) \cdots ((-m+1)z + \mathrm{c}_1(S))}{ z \cdot 2z \cdots m z } \cdot \mathrm{c}_1(S) \\
& =\frac{(-1)^{m-1} (m-1)! z^{m-1}}{ m! z^{m} } \cdot  \mathrm{c}_1(S) + ( \ldots ) \cdot \mathrm{c}_1(S)^2 \\
& = \frac{(-1)^{m-1}}{m z} \cdot  \mathrm{c}_1(S) +  ( \ldots ) \cdot \mathrm{c}_1(S)^2,
\end{align*}
by using (\ref{reduction}), we conclude that 
\[\ev_{*} [V(S^{[n]},(0,m))]^{\mathrm{vir}}= \frac{(-1)^{m-1}}{m z} \cdot  \mathrm{c}_1(S)_n+ O(1/z^2),\]
this computes the part of  $I(q,y)$ which is contained in $H^2(S^{[n]})$, 
\[[I_1(q,y)]_{\deg=2}=\log( 1 + y) \mathrm{c}_1(S)_n. \]

Assume now 
 \[-\check{\beta}=\left( \gamma, k \right)=\left( \gamma, m-\mathrm{c}_1(S)\cdot \gamma/2 \right),\]
such that $\gamma \cdot \mathrm{c}_1(S)=1$. By Lemma \ref{curves}, $\gamma$ is an exceptional curve class, or $\gamma=\mathrm{c}_1(S)$ on a degree 1 del Pezzo surface. We assume that $\gamma$ is the former.  
By the virtual dimension constraint, we have that 
\[I_0(q)=1, \quad  I_1(q)\in H^0(S^{[n]})[\![q^\beta]\!] .\]
As before, to evaluate $I_1(q)$ we need to compute its pairing with a point class $\langle I_1(q), [P] \rangle$. We take $P$ to be a collection of $n$ distinct points disjoint from $\p^1$ which represents $\gamma$. In this case, the fiber over $P$ consists of stable pairs  supported on the vertical $\Gamma$ inside $S$ in the class $\gamma$ and stable pairs supported on $n$ disjoint horizontal $\p^1$ inside $S\times \p^1$.  More precisely, the fiber has the following expression 
\[\ev^{-1}(P)=\coprod_{\underline{m}} \left( \left(\Sym^{m_0}\Gamma \times \BC\right) \times \prod^{n}_{i=1}\Sym^{m_i}\BC \right),\]
such that $m_0+1+\sum_im_i=m$ (we add 1, because $\Gamma$ is vertical, i.e.\ $n=0$, see Remark \ref{rmk} for more details). As before, the non-vanishing contribution is due to $m_0=m-1$. In this case, factors $\Sym^{0}\BC$ contribute 1 to the product formula.  The remaining factor $\Sym^{m-1}\Gamma \times \BC$ can be computed for $n=0$ via localisation on the vertex space, 
\[V(S^{[0]},(\gamma,m))=\mathrm{P}_{0, (\gamma,m)}(S\times \BC^1).\]
This is essentially the situation of a rigid smooth curve inside a threefold, because $\Gamma$ is rigid. Using the obstruction-theory computations from \cite[Section 4.2]{PT} again, we get the following expressions for a stable pair  $I^\bullet \in \mathrm{P}_{0, (\gamma,m)}(S\times \BC)^{\BC^*}$,
\begin{align*}
	\mathrm{Def}_{I^\bullet}&=H^0(\p^1,\CO_D(D))\oplus  \omega^\vee_{ \p^1,0}\\
	\mathrm{Obs}_{I^\bullet}&=H^0(\p^1,\CO_D(D)\otimes  \omega_{S} \otimes \omega_{\p^1})^\vee = H^0(\CO_D(D)\otimes  \omega_{S})^\vee \otimes  \omega_{ \p^1,0}^{\vee}.
	\end{align*}
 Using the notation of  \cite[Section 4.2]{PT}, we define $K_m$ to be the tautological rank $m$ vector bundle on $\Sym^m\Gamma$ associated to the restriction of $\omega_S$ to $ \Gamma$. The analysis in the end of \cite[Section 4.2]{PT} applies to our case verbatim with the difference that we acquire equivariant parameters. The only source of equivariance is $\omega_{ \p^1,0}^\vee$, which has weight $z$.  Since $\gamma\cdot \mathrm{c}_1(S)=1$,  we conclude that 
\begin{align*}
	\sum_{m}\deg [\mathrm{P}_{0, (\gamma,m)}(S\times \BC^1)]^{\mathrm{vir}}y^{m+1}&=\sum_m
\int_{\Sym^{m}\Gamma }\frac{ e_{\BC^*}( \mathrm{Obs}_{I^\bullet}^{\text{mov}} )}{ e_{\BC^*}( \mathrm{Def}_{I^\bullet}^{\text{mov}})}y^{m+1} \\
&= \sum_{m} \int_{\Sym^{m}\Gamma} \frac{\sum_i z^{m-i} \mathrm{c}_i(K_{m})}{z}y^{m+1} \\
&=\frac{y}{(1+y)z}.
\end{align*}
Overall, we obtain the part of $I_1(q,y)$ which is contained in $H^0(S^{[n]})$, if $S$ is not of degree 1, 
\[[I_1(q,y)]_{\deg=0}=\frac{y}{1+y}\left( \sum_{\gamma\cdot \mathrm{c}_1(S)=1}q^\gamma\right). \]
The same holds for a degree 1 del Pezzo surface with the difference that there is a term associated to $\mathrm{c}_1(S)$, which cannot be computed by the same methods, as curves in this class are not rigid. By \cite[Lemma 3.2.2]{Koll}, $h^0(S,\omega_S)\geq 2$, hence the contribution of $\mathrm{c}_1(S)$ is potentially non-trivial. 
\qed
\\

For Hilbert schemes of points, we define
\[\left\langle \gamma_1, \ldots, \gamma_N \right\rangle^{\sharp, \epsilon}_{g, \gamma}:=\sum_{m } \langle \gamma_1, \dots, \gamma_N \rangle_{g,(\gamma,m)}^{\sharp, \epsilon}y^{m}.\]
Using the wall-crossing formula from Theorem \ref{wallcrossingHilb}, Proposition \ref{IfunctionFano}, and string and divisor equations, we obtain the following result. Note that we  assume $N> 2$ to apply divisor and string equations for the degree-0 invariants  in the wall-crossing formula.

\begin{cor} \label{Delpezzowall} Assume $N>2$, then for Hilbert schemes of points $S^{[n]}$ on a del Pezzo surface $S$, we have
\[\left\langle \gamma_1, \ldots, \gamma_N \right\rangle^{\sharp, 0^+}_{g,\gamma} \\
=
(1+y)^{\mathrm{c}_{1}(S) \cdot \gamma} \cdot \left\langle \gamma_1, \ldots, \gamma_N \right\rangle^{\sharp, \infty}_{g, \gamma}.
\]
In particular,  3-point genus-0 Gromov--Witten invariants of $S^{[n]}$ are determined by Pandharipande--Thomas theory invariants of $S\times \p^1$ with relative insertions. 
\end{cor}

\appendix
\section{Stability of fibers} \label{stability-of-fibers}

\subsection{Stability of fibers versus stability} In this section, we will compare the stability of fibers of a sheaf $F$ on $S\times C$ with the stability of $F$ itself. The main results of the section is Corollary \ref{stability}, where we show that $\epsilon$-stable sheaves are in fact slope stable, if for $\mathbf{v}$ and $\CO_S(1)$ all semistable sheaves are slope stable on $S$. In Corollary \ref{converse}, we show its converse under the assumption that $\rk(\mathbf{v})=2$ and there are no strictly slope semistable sheaves. 

\subsection{Notation} Let $S\times C \rightarrow C$ be a trivial surface fibration over a connected nodal curve $C$, and let
\[\pi \colon \bigcup_i S\times C_{i} \rightarrow S\times C,\]
be its normalisation, such that $S\times C_{i}$ are its irreducible components.  For simplicity, we assume $h^1(S)=0$. 
 We  fix very ample line bundles $\CO_{S}(1)$ on $S$ and $\CO_C(1)$ on $C$.  We denote 
 \begin{align*}
 	L_k&= \CO_{S}(1)\boxtimes \CO_{C}(k), &d_i=\deg(\CO_C(1)_{|C_i}),	\\
 	d_S &=\mathrm{c}_1(\CO_S(1))^2, &d_C= \sum^m_{i=1}d_i,  \\
 	& & q_i=\frac{d_i}{d_C}. 
 	\end{align*}
 Given a sheaf $F$ on $S\times C$, let
 \[F_i:=\pi^*F_{|S\times C_i}.\] 
 Using the K\"unneth decomposition on $S\times C_i$,
 \[H^{2}(S\times C_i)=H^{2}(S)  \oplus \BQ,\]
 the first Chern class of a sheaf $F_i$ can be expressed accordingly
 \begin{equation} \label{decomps}
 \mathrm{c}_{1}(F_i)=\mathrm{c}_{1}(F_{i,p}) \oplus \beta_0(F), 
 \end{equation}
 where $F_{i,p}$ is a fiber of $F_i$ over a general point $p \in C_i$.  With respect to this decomposition, we have
 \begin{equation} \label{subs2}
 		\mathrm{c}_1(F_i)\cdot \mathrm{c}_1(L_k)^2= 2kd_i\cdot \mathrm{c}_{1}(F_{i,p})\cdot \mathrm{c}_1(\CO_S)+d_S \cdot \beta_0(F_i).
 \end{equation}
For latter, it is convenient to define the following quantities,  
\begin{align*}
\beta_0(F):&=\sum_i \beta_0(F_i), & \widetilde{\rk}(F):=\sum_i q_i \cdot \rk(F_i),\\ 
\deg(F_i)_{\mathrm{f}}:&=\mathrm{c}_{1}(F_{i,p})\cdot \mathrm{c}_1(\CO_S(1)),
\end{align*}
for more about the quantities on the left we refer to Section \ref{sheavesass}, the quantity on the right can be seen as a weighted rank of a sheaf. If $F$ has the same rank on all irreducible components, then $\widetilde{\rk}(F)=\rk(F_i)=\rk(F)$. 

 \subsection{Slope functions} \label{Sectionslopes}Recall that for a possibly singular and reducible variety $S\times C$, 
  the slope of a non-torsion sheaf $F$   with respect to a line bundle $L_k$ can be defined as follows,
\[\hat{\mu}_k(F)=\frac{a_{2}(F)}{a_{3}(F)},\]
where $a_{i}(F)$ are the coefficients of the Hilbert polynomial 
\[P(F,t)=\chi(F\otimes L^t_k)=\sum_n a_{n}(F) \frac{t^n}{n!}.\]

In the smooth case, we have 
\[\hat{\mu}_k(F)=\frac{\mu_k(F)}{a_{3}(\CO_{S\times C})}+\frac{a_{2}(\CO_{S\times C})}{a_{3}(\CO_{S\times C})},\]
where $\mu_k(F)$ is defined via the degree and the rank of $F$. A sheaf $F$ is slope semistable, if $\hat{\mu}_k(F)\geq \hat{\mu}_k(G)$ for all proper subsheaves $0 \neq G \subset F$. It is slope stable, if the inequality is strict. 

For us, a more convenient function to measure the difference of slopes of two  sheaves $G\subseteq F$ with respect to $L_k$ is
\[  \hat{\mu}_k(F,G)=\frac{a_{2}(F)\cdot a_{3}(G)-a_{2}(G)\cdot a_{3}(F)}{ a_{3}(\CO_{S\times C})}.\]
We will now derive a more explicit expression for $ \hat{\mu}_k(F,G)$. Firstly, consider the normalisation sequence 
\[0 \rightarrow F \rightarrow \bigoplus_i \pi_* F_i \rightarrow \bigoplus_{s \in \mathrm{nodes}}F_s\rightarrow 0,\] 
the sequence is exact on the left, because $F$ is torsion free. We get that 
\begin{align*}
	a_3(F)&=\sum_ia_3(F_i), \\
	a_2(F)&= \sum_i a_2(F_i)-\sum_sa_2(F_s),
\end{align*}
using the Grothendieck--Riemann--Roch theorem,  we also obtain 
\begin{equation} \label{subs}
\begin{aligned}
	a_3(F_i)&= \rk(F_i)  \cdot d_S\cdot kd_i, \\
	a_2(F_i)&=\mathrm{c}_1(F_i)\cdot \mathrm{c}_1(L_k)^2+\rk(F_i)\cdot a_2(\CO_{S\times C_i} ), \\
	a_2(\CO_{S\times C_{i}})&=d_S\cdot (1-g(C_i))+  kd_i\cdot \mathrm{c}_{1}(\CO_S(1))\cdot \mathrm{c}_1(S),\\
		a_2(F_s)&=  \rk(F_s) \cdot d_S .
\end{aligned}
\end{equation}

Assume   $F$ is flat, then the rank  and the degree of fibers of $F$ are constant, hence for any chosen $i$, we define
\[ \rk(F):=\rk(F_i), \quad \deg(F)_{\mathrm{f}}:=\deg(F_i)_{\mathrm{f}}. \]
 We then split $ \hat{\mu}_k(F,G)$ into two summands,
\begin{equation} \label{split}
 \hat{\mu}_k(F,G)=N_0+N_1, 
\end{equation}
where the first summand $N_0$ corresponds to the contributions of fibers over nodes, which admits a very simple expression after substitution of (\ref{subs}),
\begin{align*}
	N_0:&=-\frac{\left(\sum_sa_{2}(F_s)\right)\cdot a_{3}(G)-\left(\sum_sa_{2}(G_s)\right)\cdot a_{3}(F)}{ d_S\cdot kd_C} \\
	&=d_S \cdot \rk(F)\cdot \sum_s \left(\rk(G_s) -  \widetilde{\rk}(G) \right).
	\end{align*}
While the second summand $N_1$ corresponds to contributions from restrictions to irreducible components $F_i$, 
\[	N_1:=\frac{\left(\sum_ia_{2}(F_i)\right)\cdot a_{3}(G)-\left(\sum_ia_{2}(G_i)\right)\cdot a_{3}(F)}{ d_S\cdot kd_C},\]
substituting (\ref{subs2}) and (\ref{subs}), we obtain
\begin{multline*}N_1=\left(  d_S \cdot \beta_0(F)+2kd_C\cdot \deg(F)_{\mathrm{f}}+\rk(F)\cdot \sum_ia_2(\CO_{S\times C_i} ) \right) \cdot \widetilde{\rk}(G)\\
	-\left(d_S \cdot \beta_0(G)+ \sum_i \left(2kd_i\cdot \deg(G_i)_{\mathrm{f}}+\rk(G_i)\cdot a_2(\CO_{S\times C_i} ) \right) \right) \cdot  \rk(F).
\end{multline*}
Now let us rearrange the terms in $N_1$, cancelling some of them, and express it as 
\begin{equation} \label{split2}
N_1 =N_{1,0}+N_{1,1}+N_{1,2},
\end{equation}
where
\begin{align*}
N_{1,0}:&=\rk(F)\cdot \sum_i \left( \widetilde{\rk}(G) \cdot a_2(\CO_{S\times C_i} )  -\rk(G_i)\cdot a_2(\CO_{S\times C_i} )  \right), \\
&= d_S\cdot \rk(F)\cdot \sum_i \left( ( \widetilde{\rk}(G) -\rk(G_i) )\cdot (1-g(C_i))\right) \\
N_{1,1}:&= d_S \cdot \left(\widetilde{\rk}(G) \cdot \beta_0(F) -
\rk(F)\cdot \beta_0(G)\right), \\
N_{1,2}:&=\sum_{i} 2kd_i \cdot  \left( \rk(G_i)\deg(F)_\mathrm{f}- \rk(F)\deg(G_i)_{\mathrm{f}} \right).
\end{align*}
The terms $N_{1,1}$ and $N_{1,2}$ have obvious interpretations. The term $N_{1,1}$ is the slope difference associated to the ``degree" parts of the Chern characters of $F$ and $G$ in the sense of Section \ref{sheavesass}, while the term  $N_{1,2}$ is the slope difference associated to the fiber parts of $F$ and $G$. On the other hand, the summands $N_0$ and $N_{1,0}$ can be uniformly bounded from above and below, as is shown in the next lemma, which is also a summary of the preceding discussion.

\begin{lemma} \label{lemmaslopes} Assume $F$ is flat over $C$ and $\rk(F)>0$. Consider a subsheaf $G\subseteq F$, then 
\begin{multline*}	
	\sum^{j=2}_{j=1}N_{1,j} +	N(d_S, C, \rk(F))> 
	  \hat{\mu}_k(F,G)  >
	 \sum^{j=2}_{j=1}N_{1,j}-N(d_S, C, \rk(F)), 
	 \end{multline*}
where $N(d_S, C, \rk(F))$ is a number that depends only on $C$, $d_S$ and $\rk(F)$; $N_{1,1}$ and $N_{1,2}$ are defined as above. 
	\end{lemma}
\textit{Proof.} We use (\ref{split}) and (\ref{split2}). Since 
\begin{equation} \label{bound3}
  \rk(F) \geq \widetilde{\rk}(G)\geq 0, \quad \rk(F)\geq \rk(G_i)\geq 0,
  \end{equation}
we obtain that 
\begin{multline*}
d_S\cdot\rk(F)^2\cdot\left(\sum_{g(C_i)\neq0}(g(C_i)-1)+R\right)>N_{1,0}\\
>-d_S\cdot\rk(F)^2\cdot \left(\sum_{g(C_i)\neq0}(g(C_i)-1)+R\right),
\end{multline*}
where $R$ is the number of rational components in $C$.

For an upper bound of $N_0$, one has to do a little bit of work. For simplicity, assume $C$ has just one node.  Consider the sequence 
\[ 0 \rightarrow T(\pi^*G)\rightarrow \pi^*G \rightarrow \widetilde{G} \rightarrow0,\]
where $T(\pi^*G)$ is the maximal torsion subsheaf of the pullback $\pi^*G$. Since the normalisation of $C$ is smooth and of dimension 1, the sheaf  $\widetilde{G}$ is flat over the normalisation. Hence restricting the sequence to the fiber over a preimage  $s_i$ of the node $s\in C$, we obtain
\begin{equation*} 
0 \rightarrow T(\pi^*G)_{s_i}\rightarrow \pi^*G_{s_i} \rightarrow \widetilde{G}_{s_i} \rightarrow0, \quad i\in \{1,2\}.
\end{equation*}
Flatness of $\tilde{G}$ also implies that $\rk(G_i)=\rk(\widetilde{G}_{s_i})$, where we interpret  $G_i$ as the sheaf on the unique component, if $s$ is a non-separating node. Moreover, $\pi^*G_{s_i}=G_{s}$. Hence 
\begin{equation} \label{bound1}
 \rk(T(\pi^*G)_{s_i})+\rk(G_i)=\rk(G_s).
\end{equation}
Consider now the normalisation sequence,
\[0 \rightarrow G \rightarrow \pi_*\pi^*G \rightarrow G_s\rightarrow 0,\]  
since $G$ is torsion free, the pushforward $\pi_*(\oplus_i T(\pi^*G)_{s_i})$ must inject into $G_s$, 
\[\pi_*(\oplus_i T(\pi^*G)_{s_i}) \hookrightarrow G_s, \]
hence we obtain that 
\begin{equation}\label{bound22}
\sum^{i=2}_{i=1} \rk(T(\pi^*G)_{s_i})\leq \rk(G_s). 
\end{equation}
Combing (\ref{bound1}) and  (\ref{bound22}), we obtain 
\[ 0 \leq \rk(G_s)\leq \sum^{i=2}_{i=1}\rk(G_i).\]
Using (\ref{bound3}) and the bound above, we obtain a bound on $N_{0}$, 
\[ 2Bd_S \cdot \rk(F)^2>N_{1,0}>-2Bd_S \cdot \rk(F)^2,\]
where $B$ is the number of nodes in $C$. Combing bounds on $N_0$ and $N_{1,0}$, we obtain that 
\begin{align*}
	N(d_S, C, \rk(F))=d_S\cdot\rk(F)^2\cdot \left(\sum_{g(C_i)\neq0}(g(C_i)-1)+R+2B\right).
	\end{align*}
\qed 

\subsection{Stability of fibers implies stability}

For the next corollary,  following the discussion in Section \ref{sheaves}, we define
\[L_k:= \CO_S(1)\boxtimes f^* \CL_\beta^k.\]
We refer to Section \ref{Sectionslopes} for the definition of slopes for singular varieties. 

\begin{cor} \label{stability} Assume that for $\mathbf{v}$ and $\CO_S(1)$, all semistable sheaves are slope stable, or $\rk(\mathbf{v})=1$.  Fix a class $\beta \in \Eff(M(\mathbf{v}), \rCoh(S,\mathbf{v}))$ and $\epsilon \in \BR_{>0}\cup\{0^+,\infty\}$.  There exists $k_0 \in \BZ$, such that  $\epsilon$-stable sheaves with the Chern character $(\ch(\mathbf{v}),\check{\beta})$ are slope stable with respect to $L_k$ for all $k\geq k_0$. 
\end{cor}


\textit{Proof.} 
 If $\rk(\mathbf{v})=1$, the claim holds for a simple reason: on a smooth variety, a rank $1$ sheaf is stable, if and only if it is torsion free. More precisely, by tensoring sheaves with a line bundle we may assume $\mathbf{v}=(1,0,-n)$, then  if we choose a right section of $\Coh(S,\mathbf{v}) \rightarrow \Coh_r(S,\mathbf{v})$, sheaves $F$ associated to quasimaps are ideal sheaves of curves, as explained in Section \ref{relativeHilbert}. Ideal sheaves are stable with respect to all ample line bundles on $S \times C$. 

We therefore assume $\rk(\mathbf{v})\geq 2$. Given an $\epsilon$-stable sheaf $F$, a general fiber of $F$ over $C$ is stable, and, in particular, torsion free. Hence by Lemma \ref{torsion}, the sheaf $F$ is torsion free itself. Moreover, by the stability of a general fiber of $F$, a general fiber of any saturated subsheaf $G \subset F$ is non-destabilizing. Hence by Lemma \ref{lemmaslopes}, the difference of slopes,
\[\hat{\mu}_k(F)-\hat{\mu}_k(G),\]
with respect to a line bundle $L_k$ can be made positive for a large enough $k$. Indeed,
in the notation of Lemma \ref{lemmaslopes}, the term $N_{1,2}$, which depends just on Chern characters of general fibers, can be made arbitrarily big, while the term $N_{1,1}$ does not depend on $k$. Moreover, if $G$ is non-destabilizing for some $k'$, then it stays non-destabilizing for all $k\geq k'$ for the same reason.  

Now, for a fixed $k'$, the family of $L_{k'}$-destabilising subsheaves of $F$ is bounded by \cite[Lemma 1.7.9]{HL}. By the boundedness of this family and the discussion above,  there exists $k_0 \geq k'$, such that all subsheaves in this family become non-destabilizing. On the other hand, all other subsheaves remain non-destabilizing. Hence  $F$ is slope stable for all $k \geq k_0$.

By Corollary \ref{algebraicity}, the family of $\epsilon$-stable quasimaps  is bounded, therefore there exists a uniform choice of $k_{0}$, such that the above conclusion holds for all $\epsilon$-stable sheaves. 
\qed

\subsection{Stability implies stability of fibers}

\begin{prop}\label{stab} 
	Assume $C$ is smooth. Fix a class \[\alpha=(\alpha_0, \alpha_1, \alpha_2, \alpha_3) \in H^{*}(S\times C),\] such that $\alpha_0=2$ (i.e. $\rk=2$). There exists $k_{0}\in \BN$,  such that for all $k\geq k_{0}$ and for all torsion free sheaves $F$ with $\ch(F)=\alpha$, the following statement holds:  a saturated\footnote{Saturation is needed to ensure that fibers are indeed subsheaves $G_{p}\subset F_{p}$  for a general point $p$.} subsheaf $G \subset F$ is slope $L_{k}$-destabilizing,  if $G_{p} \subset F_{p}$ is slope destabilizing for a general $p \in C$.    
\end{prop}  

\textit{Proof.} 
We will prove the proposition by restricting to a hyperplane section and then applying \cite[Theorem 5.3.2]{HL} together with \cite[Remark 5.3.5]{HL}.

By (\ref{decomps}), we can decompose the first Chern class of $F$ as follows, 
\[ \mathrm{c}_{1}(F)=\mathrm{c}_{1}(F_{p}) \oplus \beta_0(F) \in H^{2}(S)  \oplus \BQ ,\]
where each summand is in the corresponding K\"unneth component, and $F_{p}$ is a general fiber of $F$ over $p\in C$. The intersection numbers with a curve class $\mathrm{c}_1(L_{n})\cdot \mathrm{c}_1(L_{m})$ take the following form
\begin{equation}\label{int1}
	\mathrm{c}_{1}(F)\cdot \mathrm{c}_1(L_{n})\cdot \mathrm{c}_1(L_{m})=(n+m)d_C\cdot \deg(F)_{\mathrm{f}} + d_S \cdot  \beta_0(F).
\end{equation} 
In particular, slope-stability with respect to a curve class $\mathrm{c}_1(L_{1})\cdot \mathrm{c}_1(L_{2k-1})$ coincides with slope-stability with respect to a curve class $\mathrm{c}_1(L_{k})\cdot \mathrm{c}_1(L_{k})$.

Consider now a general smooth hyperplane section 
\[H \in H^0(S\times C, L_1).\]
Let $2k_{0}-1$ be the smallest odd integer such that  \cite[Remark 5.3.5]{HL} and  \cite[Theorem 5.3.2]{HL} hold for the surface fibration $H\rightarrow C$, the class $\beta_{|H}$ and a polarisation $L_{2k_{0}-1|H}$. 

Consider now a saturated subsheaf $G \subset F$, such that  $G_{p} \subset F_{p}$ is destabilizing for a general $p \in C$. Then this also holds for a restriction $G_{|H} \subset F_{|H}$. Hence by   \cite[Theorem 5.3.2]{HL}, the subsheaf $G_{|H}$ destabilizes  $F_{|H}$ with respect to $ L_{2k_{0}-1|H}$. This implies that $G$ destabilises $F$ with respect to $\mathrm{c}_1(L_1)\cdot \mathrm{c}_1(L_{2k_{0}-1})$, and  by (\ref{int1})  with respect to $\mathrm{c}_1(L_{k_0})\cdot \mathrm{c}_1(L_{k_{0}})$. The same applies to all $k \geq k_0$. 
\qed

\begin{rmk} The proof of Proposition \ref{stab} is inspired by the proof of \cite[Proposition 4.2]{T}, which, however, contains a mistake.
	A restriction of a sheaf $F$ to a hyperplane section is stable with respect to the polarisation that defines the hyperplane section, which is not necessarily suitable in the sense of \cite[Theorem 5.3.2]{HL}. If one adds fiber classes to the polarisation to make it suitable, then one has to  take a hyperplane section of a bigger degree, for which suitable polarisation may be different.
\end{rmk}
\begin{rmk} We assume that $\rk=2$, because for fibred surfaces,  Proposition \ref{stab} does not work for $\rk>2$. Indeed, destabilizing sheaves might be of higher rank, hence \cite[(5.3)]{HL} is no longer true. For $\rk>2$, the correct approach is provided by \cite[Proposition 6.2]{Yo}, so in the end Corollary \ref{main1} still holds true for fibered surfaces. However, we do not know how to use the hyperplane-restriction tricks in the context of \cite{Yo} to extend the result to higher dimensions. 
\end{rmk}



\begin{cor} \label{main1} Let $C$ be a nodal curve. 
	Fix classes $\alpha_{i}\in H^{*}(S\times C_{i})$, such that $\alpha_{i,0}=2$. There exists $k_{0}\in \BN$, such that for all $k\geq k_{0}$ and for all sheaves $F$ flat over $C$ with $\ch(F_{i})=\alpha_{i}$, the following statement holds:   a saturated subsheaf $G \subset F$ is slope $L_{k}$-destabilizing,  if $G_{p} \subset F_{p}$ is slope destabilizing for a general $p \in C$.  
\end{cor}
\textit{Proof.} By (\ref{split}), the slop difference  $ \hat{\mu}_k(F_i,G_i)$ on a component of the normalisation $C_i$ can be expressed as follows, 
\begin{multline*}
 \hat{\mu}_k(F_i,G_i)= d_S \cdot \left(\rk(G_i) \cdot \beta_0(F_i) -
\rk(F)\cdot \beta_0(G_i)\right) \\
 +  2kd_i \cdot  \left( \rk(G_i)\deg(F)_\mathrm{f}- \rk(F)\deg(G_i)_{\mathrm{f}} \right).
\end{multline*}
Note that the terms $N_0$ and $N_{1,0}$ vanish, because $C_i$ is smooth and irreducible. 
 Lemma \ref{lemmaslopes} therefore gives us that
\begin{multline*}
 \hat{\mu}_k(F,G)< \sum_i  \hat{\mu}_k(F_i,G_i)+d_S\cdot \sum_i \left(\widetilde{\rk}(G) - \rk(G_i) \right) \cdot \beta_0(F_i)\\
+ N(d_S, C, \rk(F)), 
\end{multline*}
where we used $\beta_0(F)=\sum_i \beta_0(G_i)$ and $\beta_0(F)=\sum_i \beta_0(G_i)$. 
Applying Proposition \ref{stab} to each component $C_i$ and increasing  $k$ by a finite amount to cancel extra bounded contributions on the right of the inequality above, we can make the quantity $ \hat{\mu}_k(F,G)$ negative. Since the last two summands on the right can be bounded in terms of the given Chern-character data and degree of $\CO_C(1)$ on irreducible components, the resulting $k$ is uniform for all $F$ in the given class.   \qed

\begin{cor} \label{converse} Assume we are in the situation of Corollary \ref{main1}, and that there are no strictly slope semistable sheaves on fibers. There exists $k_{0}\in \BN$, such that for all $k\geq k_{0}$ and for all sheaves $F$ flat over $C$ with $\ch(F_{i})=\alpha_{i}$, the following statement holds:   a sheaf $F$ is slope $L_{k}$-stable,  only if $F_{p}$ is slope stable for a general $p \in C$.  	\end{cor}

\textit{Proof.} Assume $F_{p}$ is not slope stable  for all $p \in C$, i.e.\ it is slope unstable by the assumption. Then by \cite[Section 2.3]{HL}, we can construct a relative Harder--Narasimhan filtation over an open dense subset in $C$. We  extend it to the whole curve $C$ and take its lowest piece,  $G\subset  F$.  We may assume that $G$ is saturated,  because by construction the injection $G \hookrightarrow F$ restricts to an injection $G_p\hookrightarrow F_p$ for a general $p \in C$. Hence taking the saturation does not affect a general fiber of $G$, i.e.\ it remains destabilizing. The claim then follows from Corollary \ref{main1} applied to the subsheaf $G$. 
\qed

\section{Flatness}
\subsection{Flatness is an open condition} We will show that flatness of sheaves on $S\times C$ over $C$ is an open condition for sheaves on $S\times C$.   We managed to bypass this result in the main part of the article by working with quasimaps instead of sheaves in some places, e.g.\ Proposition \ref{openess}. 


   Let $\Coh_{\FM_{g,N}}(S\times \FC_{g,N})$ be the stack of all coherent sheaves $F$ on moving threefolds $S\times C$. Let 
\[ \Coh_{\FM_{g,N}}^{\mathrm{flat}}(S\times \FC_{g,N}) \subset \Coh_{\FM_{g,N}}(S\times \FC_{g,N})\]  be the locus of sheaves on $S\times C$ flat over $C$.

\begin{lemma}  \label{flatopen} The locus $\Coh_{\FM_{g,N}}^{\mathrm{flat}}(S\times \FC_{g,N})$ is an open substack. 
\end{lemma}
\textit{Proof.} Let $\pi \colon \CC \rightarrow B$ be family of curves, and let  $F$ be a family of sheaves on $S\times \CC$ flat over $B$. Assume that the sheaf $F_{b_0}$ on $S \times \CC_{b_0}$ is flat over $\CC_{b_0}$ for some closed point $b_0\in B$. We will show that there is a Zariski open neighbourhood  $U\subset B$ containing $b_0$, such that for all $b \in U$, the sheaf $F_b$ is flat over $\CC_{b}$.  

By Lemma \ref{flatness}, the family $F$ is flat at  $\CC_{b_0} \subset \CC$. Hence by the usual openness of flatness, there is a Zariski open neighbourhood $\CU' \subseteq  \CC$ containing $\CC_{b_0}$, over which $F$ is flat. We want to make $\CU'$ proper and flat over $B$, i.e.\ we want to get rid of fibers, which are not completely contained in $\CU'$. The set $\CU'$  intersects a fiber $\CC_b$ but does not  fully contain it, if and only if  $\CC_b$  intersects the complement $\CU'^c$. Hence we can get rid of such fibers by the following construction, 
\[\CU=\CU' \cap \pi^{-1}(\pi(\CU'^{c}))^c \subseteq \CC,\] 
since $\pi$ is a closed map, we obtain that $\CU$ is a Zariski open  neighbourhood, such that $\CU$  intersect a fiber $\CC_b$, if and only if $\CU$ contains it. Moreover, $\CU$ is not empty, because $\CC_{b_0}\subset \CU$. By construction, the image of $\CU$ in $B$, 
\[\pi(\CU) \subseteq B\] 
is the desired Zariski open neighbourhood. Hence there is a Zariski open neighbourhood around all sheaves $F \in \Coh_{\FM_{g,N}}^{\mathrm{flat}}(S\times \FC_{g,N})(\BC)$ inside the stack  $\Coh_{\FM_{g,N}}(S\times \FC_{g,N})$ which is also contained in $\Coh_{\FM_{g,N}}^{\mathrm{flat}}(S\times \FC_{g,N})$. We conclude that $\Coh_{\FM_{g,N}}^{\mathrm{flat}}(S\times \FC_{g,N})$ is open. 
\qed

\bibliographystyle{amsalpha}
\bibliography{QMs}
\end{document}